# Robust Sample Average Approximation

Dimitris Bertsimas · Vishal Gupta · Nathan Kallus

**Acknowledgements** This material is based upon work supported by the National Science Foundation Graduate Research Fellowship under Grant No. 1122374.

**Abstract** Sample average approximation (SAA) is a widely popular approach to data-driven decision-making under uncertainty. Under mild assumptions, SAA is both tractable and enjoys strong asymptotic performance guarantees. Similar guarantees, however, do not typically hold in finite samples. In this paper, we propose a modification of SAA, which we term Robust SAA, which retains SAA's tractability and asymptotic properties and, additionally, enjoys strong finite-sample performance guarantees. The key to our method is linking SAA, distributionally robust optimization, and hypothesis testing of goodness-of-fit. Beyond Robust SAA, this connection provides a unified perspective enabling us to characterize the finite sample and asymptotic guarantees of various other data-driven procedures that are based upon distributionally robust optimization. This analysis provides insight into the practical performance of these various methods in real applications. We present examples from inventory management and portfolio allocation, and demonstrate numerically that our approach outperforms other data-driven approaches in these applications.

## 1 Introduction

In this paper, we treat the stochastic optimization problem

$$z_{\text{stoch}} = \min_{x \in X} \mathbb{E}_F[c(x; \xi)], \tag{1}$$

where $c(x; \xi)$ is a given cost function depending on a random vector $\xi$ following distribution $F$ and a decision variable $x \in X \subseteq \mathbb{R}^{d_x}$. This is a widely used modeling paradigm in operations research, encompassing a number of applications [10, 46].

In real-world applications, however, the distribution $F$ is unknown. Rather, we are given data $\xi^1, \ldots, \xi^N$, which are typically assumed to be drawn IID from $F$. The most common approach in these settings is the sample average approximation (SAA). SAA approximates the true, unknown distribution $F$ by the

Dimitris Bertsimas
Sloan School of Management
Massachusetts Institute of Technology, Cambridge, MA 02139
E-mail: dbertsim@mit.edu

Vishal Gupta
Marshall School of Business
University of Southern California, Los Angeles, CA 90029
E-mail: guptavis@usc.edu

Nathan Kallus
School of Operations Research and Information Engineering and Cornell Tech
Cornell University, New York, NY 10011
E-mail: kallus@cornell.edu

empirical distribution $\hat{F}_N$, which places $1/N$ mass at each of the data points. In particular, the SAA approach approximates (1) by the problem

$$\hat{z}_{\text{SAA}} = \min_{x \in X} \frac{1}{N} \sum_{j=1}^{N} c(x; \xi^j). \tag{2}$$

Variants of the SAA approach in this and other contexts are ubiquitous throughout operations research, often used tacitly without necessarily being referred to by this name.

Under mild conditions on the cost function $c(x; \xi)$ and the sampling process, SAA enjoys two important properties:

**Asymptotic Convergence:** As the number of data points $N \to \infty$, both the optimal value $\hat{z}_{\text{SAA}}$ of (2) and an optimal solution $x_{\text{SAA}}$ converge to the optimal value $z_{\text{stoch}}$ of (1) and an optimal solution $x_{\text{stoch}}$ almost surely (e.g. [27, 29]).

**Tractability:** For many cost functions $c(x; \xi)$ and sets $X$, finding the optimal value of and an optimal solution to (2) is computationally tractable (e.g. [10]).

In our opinion, these two features – asymptotic convergence and tractability – underlie SAA's popularity and practical success in data-driven settings.

This is not to suggest, however, that SAA is without criticism. First, useful a priori performance guarantees do not exist for SAA with finite $N$, except in certain special cases (e.g. [29, 31]). Consequently some authors have suggested modifications of SAA to evaluate performance with a fixed, finite sample. One such example is a two-sample SAA approach (2-SAA) wherein SAA is applied to one half of the data to find a solution $x_{\text{2-SAA}}$ and the other half is used for a posteriori out-of-sample performance evaluation (see Section 1.3). Although this approach (and similar variants) can yield performance bounds under certain assumptions, the bound comes at the price of sacrificing half of the data. This sacrifice may be impractical for moderate $N$. Moreover, as we illustrate, the quality of this bound depends strongly on the unknown distribution $F$ and, even for non-pathological examples, the approach may fail dramatically (see Section 7).

A more subtle, and indeed more problematic, criticism is that in many settings the SAA solution $x_{\text{SAA}}$ and its true out-of-sample cost $\mathbb{E}_F[c(x_{\text{SAA}}; \xi)]$ may be highly unstable for finite, moderate $N$; small changes in the data or small changes in $N$ can yield large changes in the solution. This phenomenon is perhaps best documented for portfolio allocation problems [16, 32] and in ill-conditioned or under-determined regression [25, 50]. The 2-SAA method also suffers from this instability and only exacerbates the instability by halving the amount of data upon which the SAA solution is computed. See Section 7 for a numerical illustration.

In this paper, we seek to address these two criticisms and build upon SAA's practical success. We propose a novel approach to (1) in data-driven settings which we term *Robust SAA*. Robust SAA inherits SAA's favorable asymptotic convergence and tractability. Unlike SAA, however, Robust SAA enjoys a strong *finite sample* performance guarantee for a wide class of optimization problems, and we demonstrate that its solutions are stable, even for small to moderate $N$. The key idea of Robust SAA is to approximate (1) by a particular data-driven, distributionally robust optimization problem using ideas from statistical hypothesis testing.

More specifically, a distributionally robust optimization (DRO) problem is

$$\overline{z} = \min_{x \in X} \mathcal{C}(x; \mathcal{F}), \tag{3}$$

$$\text{where } \mathcal{C}(x; \mathcal{F}) = \sup_{F_0 \in \mathcal{F}} \mathbb{E}_{F_0}[c(x; \xi)], \tag{4}$$

where $\mathcal{F}$ is a set of potential distributions for $\xi$. We call such a set a *distributional uncertainty set* or DUS in what follows. Initial research (see literature review below) focused on DUSs $\mathcal{F}$ specified by fixing the first few moments of a distribution or other structural features, but did not explicitly consider the data-driven setting. Recently, the authors of [13, 15] took an important step forward proposing specific data-driven DRO formulations in which the DUS $\mathcal{F}$ is a function of the data, i.e., $\mathcal{F} = \mathcal{F}(\xi^1, \ldots, \xi^N)$, and showing that (3) remains tractable. Loosely speaking, their DUSs consist of distributions whose first few moments are close to the sample moments of the data. The authors show how to tailor these DUSs so that for any $0 \le \alpha \le 1$, the probability (with respect to data sample) that the true (unknown) distribution $F \in \mathcal{F}(\xi^1, \ldots, \xi^N)$ is at least $1 - \alpha$. Consequently, solutions to (3) based on these DUSs enjoy a distinct, finite-sample guarantee:

**Finite-Sample Performance Guarantee:** With probability at least $1 - \alpha$ with respect to the data sampling process, for any optimal solution $\overline{x}$ to (3), $\overline{z} \ge \mathbb{E}_F[c(\overline{x}; \xi)]$, where the expectation is taken with respect to the true, unknown distribution $F$.

In contrast to SAA, however, the methods of [13, 15] do not generally enjoy asymptotic convergence. (See Section 4.3.2).

Our approach, Robust SAA, is a particular type of data-driven DRO. Unlike existing approaches, however, our DUSs are not defined in terms of the sample moments of the data, but rather are specified as the confidence region of a goodness-of-fit (GoF) hypothesis test. Intuitively, our DUSs consist of all distributions which are "small" perturbations of the empirical distribution – hence motivating the name Robust SAA – where the precise notion of "small" is determined by the choice of GoF test. Different GoF tests yield different DUSs with different computational and statistical properties.

We prove that like other data-driven DRO proposals, Robust SAA also satisfies a finite-sample performance guarantee. Moreover, we prove that for a wide-range of cost functions $c(x; \xi)$, Robust SAA can be reformulated as a tractable convex optimization problem. In most cases of practical interest, we show that it reduces to a linear or second-order cone optimization problem suitable for off-the-shelf solvers. In general, it can be solved efficiently by cutting plane algorithms. (See [12] for a review of convex optimization and algorithms.) Unlike other data-driven DRO proposals, however – and this is key – we prove that Robust SAA also satisfies an asymptotic convergence property similar to SAA. In other words, Robust SAA combines the strengths of both the classical SAA and data-driven DRO. Computational experiments in inventory management and portfolio allocation confirm that these properties translate into higher quality solutions for these applications in both small and large sample contexts.

In addition to proposing Robust SAA as an approach to addressing (1) in data-driven settings, we highlight a connection between GoF hypothesis testing and data-driven DRO more generally. Specifically, we show that any DUS that enjoys a finite-sample performance guarantee, including the methods of [13, 15], can be recast as the confidence region of *some* statistical hypothesis test. Thus, hypothesis testing provides a unified viewpoint. Adopting this viewpoint, we characterize the finite-sample and asymptotic performance of DROs in terms of certain statistical properties of the underlying hypothesis test, namely significance and consistency. This characterization highlights an important, new connection between statistics and data-driven DRO. From a practical perspective, our results allow us to describe which DUSs are best suited to certain applications, providing important modeling guidance to practitioners. Moreover, this connection motivates the use of well-established statistical procedures like bootstrapping in the DRO context. Numerical experimentation confirms that these procedures can significantly improve upon existing algorithms and techniques.

To summarize our contributions:

1. We propose a new approach to optimization in data-driven settings, termed Robust SAA, which enjoys both finite sample and asymptotic performance guarantees for a wide-class of problems.
2. We develop new connections between SAA, DRO and statistical hypothesis testing. In particular, we characterize the finite-sample and asymptotic performance of data-driven DROs in terms of certain statistical properties of a corresponding hypothesis test, namely its significance and consistency.
3. Leveraging the above characterization, we shed new light on the finite sample and asymptotic performance of existing DRO methods and Robust SAA. In particular, we provide practical guidelines on designing appropriate DRO formulations for specific applications.
4. We prove that Robust SAA yields tractable optimization problems that are solvable in polynomial time for a wide class of cost functions. Moreover, for many cases of interest, including two-stage convex optimization with linear recourse, Robust SAA leads to single-level convex optimization formulations that can be solved using off-the-shelf software for linear or second-order optimization.
5. Through numerical experiments in inventory management and portfolio allocation, we illustrate that Robust SAA leads to better performance guarantees than existing data-driven DRO approaches, has performance similar to classical SAA in the large-sample regime, and is significantly more stable than SAA in the small-to-moderate-sample regime with comparable or better performance.
6. Finally, we show how Robust SAA can be used to obtain approximations to the "price of data" – the price one would be willing to pay in a data-driven setting for additional data.

The remainder of this paper is structured as follows. We next provide a brief literature review and describe the model setup. In Section 2, we illustrate the fundamental connection between DRO and the confidence regions of GoF tests and explicitly describe Robust SAA. Section 3 connects the significance of the hypothesis test to the finite-sample performance of a DRO. Section 4 connects the consistency of the hypothesis test to the asymptotic performance of the DRO. Section 5 proves that for the tests we consider, Robust SAA leads to a tractable optimization problem for many choices of cost function. Finally, Section 7 presents an empirical study and Section 8 concludes. All proofs except that for Proposition 1 are in the appendix.

## 1.1 Literature review

As summarized by [46, Ch. 5], there are two streams of SAA literature distinguished by how they regard the data: "[the] random sample $[\xi^1, \ldots, \xi^n]$ can be viewed as historical data of $N$ observations of $\xi$, or it can be generated in the computer by Monte Carlo sampling techniques," More explicitly, SAA can be used as a solution for (1) either when $F$ is not known, but data drawn IID from it is available, or when $F$ is known and random variates can easily be simulated from it.

Although given a particular sample the SAA procedure is equivalent in both contexts, the pratical usage of SAA differs in each setting. In a Monte Carlo setting, it is easy to generate additional data as needed (e.g., [35]). Thus, the key issues are understanding how many samples are needed and how to compute the SAA solution efficiently. In this context [34] have proposed methods to evaluate both upper and lower bounds on the out-of-sample performance of the SAA solution. By contrast, in the historical data setting, additional data is unavailable. Thus, the key issue is using the given data most efficiently. [1, 31] are both examples of SAA used in this context. In statistics, SAA is also used in this context under the name empirical risk minimization [51], an example of which is the usual ordinary least squares.

In this paper, we focus on SAA (and other solutions to (1)) in the historical data setting, where a single, fixed data sample of size $N$ is given.

DRO was first proposed in the operations research literature by the author in [42], where $\mathcal{F}$ is taken to be the set of distributions with a given mean and covariance in a specific inventory context. DRO has since received much attention, with many authors focusing on DUSs $\mathcal{F}$ defined by fixing the first few moments of the distribution [8, 11, 37, 38], although some also consider other structural information such as unimodality [19]. In [55], the authors characterized the computational tractability of (3) for a wide range of DUSs $\mathcal{F}$ by connecting tractability to the geometry of $\mathcal{F}$.

As mentioned, in [13, 15], the authors extended DRO to the data-driven setting. In [13], the authors studied chance constraints, but their results can easily be cast in the DRO setting. Both papers focus on tractability and the finite-sample guarantee of the resulting formulation. Neither considers asymptotic performance. In [26], the authors also propose a data-driven approach to chance constraints, but do not discuss either finite sample guarantees or asymptotic convergence. Using our hypothesis testing viewpoint, we are able to complement and extend these existing works and establish a unified set of conditions under which the above methods will enjoy a finite-sample guarantee and/or be asymptotically convergent.

Recently, other work has considered hypothesis testing in certain, specific optimization contexts. In [7], the authors show how hypothesis tests can be used to construct uncertainty sets for robust optimization problems, and establish a finite-sample guarantee that is similar in spirit to our own. They do not, however, consider asymptotic performance. In [3], the authors consider robust optimization problems described by phi-divergences over uncertain, discrete probability distributions with finite support and provide tractable reformulations of these constraints. The authors mention that these divergences are related to GoF tests for discrete distributions, but do not explicitly explore asymptotic convergence of their approach to the full-information optimum or the case of continuous distributions. In [2], the authors survey the use of these methods for stochastic optimization with known, finitely-many scenarios but with unknown probabilities in data-driven setting and provide a classification of the various phi-divergences in terms of their ability to "suppress" observed scenarios or "pop" unobserved scenarios. Similarly, in [28], the authors study a stochastic lot-sizing problem under discrete distributional uncertainty described by Pearson's $\chi^2$ GoF test and develop a dynamic programming approach to this particular problem. The authors establish conditions for asymptotic convergence for this problem but do not discuss finite sample guarantees.

By contrast, we provide a systematic study of GoF testing and data-driven DRO. By connecting these problems with the existing statistics literature, we provide a unified treatment of both discrete and continuous distributions, finite-sample guarantees, and asymptotic convergence. Moreover, our results apply in a general optimization context for a large variety of cost functions. We consider this viewpoint to both unify and extend these previous results.

## 1.2 Setup

In the remainder, we denote the support of $\xi$ by $\Xi$. We assume $\Xi \subseteq \mathbb{R}^d$ is closed, and denote by $\mathcal{P}(\Xi)$ the set of (Borel) probability distributions over $\Xi$. For any probability distribution $F_0 \in \mathcal{P}(\Xi)$, $F_0(A)$ denotes the probability of the event $\xi \in A$. To streamline the notation when $d = 1$, we let $F_0(t) = F_0((-\infty, t])$. When $d > 1$ we also denote by $F_{0,i}$ the univariate marginal distribution of the $i^{\text{th}}$ component, i.e., $F_{0,i}(A) = F_0(\{\xi : \xi_i \in A\})$. We assume that $X \subseteq \mathbb{R}^{d_x}$ is closed and that for any $x \in X$, $\mathbb{E}_F[c(x;\xi)] < \infty$ with respect to the true distribution, i.e., the objective function of the full-information stochastic problem (1) is well-defined.

When $\Xi$ is unbounded, (3) may not admit an optimal solution. (We will see non-pathological examples of this behavior in Section 3.2.) To be completely formal in what follows, we first establish sufficient conditions for the existence of an optimal solution. Recall the definition of equicontinuity:

**Definition 1** A set of functions $\mathcal{H} = \{h : \mathbb{R}^{m_1} \to \mathbb{R}^{m_2}\}$ is *equicontinuous* if for any given $x \in \mathbb{R}^{m_1}$ and $\epsilon > 0$ there exists $\delta > 0$ such that for all $h \in \mathcal{H}$, $\|h(x) - h(x')\| < \epsilon$ for any $x'$ with $\|x - x'\| < \delta$.

In words, equicontinuity generalizes the usual definition of continuity of a function to continuity of a set of functions. If the functions are differentiable, equicontinuity is equivalent to the existence of pointwise extrema of derivatives.

Our sufficient conditions constitute an analogue of the classical Weierstrass Theorem for deterministic optimization (see, e.g., [5], pg. 669):

**Theorem 1** *Suppose there exists $x_0 \in X$ such that $\mathcal{C}(x_0; \mathcal{F}) < \infty$ and that $c(x; \xi)$ is equicontinuous in $x$ over all $\xi \in \Xi$. If either $X$ is compact or $\lim\limits_{\|x\| \to \infty} c(x; \xi) = \infty$ for all $\xi$, then the optimal value $\overline{z}$ of (3) is finite and is achieved at some $\overline{x} \in X$.*

### 1.3 Two-sample SAA approach

As mentioned, authors such as [34] have proposed modifications of SAA to derive performance guarantees principally in the Monte Carlo setting. These same bounds can be adapted to the historical data setting. The upper bound, in particular, is based on a sample-splitting procedure that we call two-sample SAA (2-SAA). Although the authors originally conceived of this procedure for the Monte Carlo setting, it is identical to hold-out validation, a statistical technique commonly used with SAA in the historical data setting [21].

In 2-SAA, we split our data and treat it as two independent samples:

$$\{\xi^1, \ldots, \xi^{\lceil N/2 \rceil}\} \quad \text{and} \quad \{\xi^{\lceil N/2 \rceil + 1}, \ldots, \xi^N\}.$$

(We focus on the half-half split as it is most common; naturally, other splitting ratios are possible.) On the first half of the data, we train an SAA decision,

$$x_{\text{2-SAA}} \in \arg\min_{x \in X} \frac{1}{\lceil N/2 \rceil} \sum_{j=1}^{\lceil N/2 \rceil} c(x; \xi^j),$$

and, on the second half of the data, we evaluate its performance, yielding $\lfloor N/2 \rfloor$ IID observations

$$c(x_{\text{2-SAA}}; \xi^{\lceil N/2 \rceil + 1}), \cdots, c(x_{\text{2-SAA}}; \xi^N)$$

of the random variable $c(x_{\text{2-SAA}}; \xi)$, the mean of which is the true out-of-sample cost of $x_{\text{2-SAA}}$, i.e.,

$$z_{\text{2-SAA}} = \mathbb{E}_F\left[c(x_{\text{2-SAA}}; \xi)\big|\xi^1, \ldots, \xi^{\lceil N/2 \rceil}\right].$$

Hence, we may use the $\lfloor N/2 \rfloor$ IID observations of $c(x_{\text{2-SAA}}; \xi)$ to form a one-sided confidence interval for the mean in order to get a guaranteed bound on the true out-of-sample cost of $x_{\text{2-SAA}}$. Since exact, nonparametric confidence regions for the mean can be challenging to compute,[1] standard practice (cf. [34]) is to approximate a confidence region using Student's T-test, yielding the bound

$$\overline{z}_{\text{2-SAA}} = \hat{\mu}_{\lfloor N/2 \rfloor} + \hat{\sigma}_{\lfloor N/2 \rfloor} Q_{T_{\lfloor N/2 \rfloor - 1}}(\alpha) / \sqrt{\lfloor N/2 \rfloor}, \tag{5}$$

where $\hat{\mu}_{\lfloor N/2 \rfloor}$ and $\hat{\sigma}_{\lfloor N/2 \rfloor}$ are the sample mean and sample standard deviation (respectively) of the $\lfloor N/2 \rfloor$ cost observations $c(x_{\text{2-SAA}}; \xi^{\lceil N/2 \rceil + 1}), \ldots, c(x_{\text{2-SAA}}; \xi^N)$ and $Q_{T_{\lfloor N/2 \rfloor - 1}}(\alpha)$ is the $(1 - \alpha)^{th}$ quantile of Student's T-distribution with $\lfloor N/2 \rfloor - 1$ degrees of freedom (cf. [39]). Were the approximation exact, $\overline{z}_{\text{2-SAA}}$ would bound $z_{\text{2-SAA}}$ with probability $1 - \alpha$ with respect to data sampling. The accuracy of the approximation by Student's T-test depends upon the normality of the sample mean and the chi distribution of the sample standard deviation. Although these may be reasonable assumptions for "large enough" $N$, we will see in Section 7 that for some problems "large enough" is impractically large. In these settings (5) is often not a valid bound at the desired probability $1 - \alpha$.

Finally, as mentioned, 2-SAA only uses half the data to compute its solution $x_{\text{2-SAA}}$, using the other half only for "post-mortem" analysis of performance. Hence it may exacerbate any instability suffered by SAA for small or moderate $N$. We next introduce Robust SAA, which uses all of the data both to compute a solution and to assess its performance and, more importantly, seeks to proactively immunize this solution to ambiguity in the true distribution.

---

[1] Indeed, nonparametric comparison tests focus on other location parameters such as median.

## 2 Goodness-of-Fit testing and Robust SAA

In this section, we provide a brief review of GoF testing as it relates to Robust SAA. For a more complete treatment, including the wider range of testing cases possible, we refer the reader to [14, 49].

Given IID data $\xi^1, \ldots, \xi^N$ and a hypothetical distribution $F_0$ (chosen a priori, not based on the data), a GoF test considers the hypothesis

$$H_0: \text{ The data } \xi^1, \ldots, \xi^N \text{ were drawn from } F_0 \qquad (6)$$

and tests it against the alternative

$$H_1: \text{ The data } \xi^1, \ldots, \xi^N \text{ were not drawn from } F_0.$$

A GoF test rejects $H_0$ in favor of $H_1$ if there is sufficient evidence, otherwise making no particular conclusion. A test is said to be of significance level $\alpha$ if the probability of incorrectly rejecting $H_0$ is at most $\alpha$. Notice that the null hypothesis $H_0$ is fully characterized by $F_0$.[2] Consequently, we refer to $H_0$ and $F_0$ interchangeably. Thus, we may say that a test rejects *a distribution* $F_0$ if it rejects the null hypothesis that $F_0$ describes.

A typical test specifies a statistic

$$S_N = S_N(F_0, \xi^1, \ldots, \xi^N)$$

that depends on the data $\xi^1, \ldots, \xi^N$ and the hypothetical distribution $F_0$ and a threshold $Q_{S_N}(\alpha)$ that does not depend on either the data or $F_0$ (it may, however, depend on the true distribution F). The test rejects $H_0$ if $S_N > Q_{S_N}(\alpha)$.

Let $Q_{S_N}^*(\alpha)$ be the $(1-\alpha)^{\text{th}}$ quantile of the distribution of $S_N(F, \xi^1, \ldots, \xi^N)$, i.e., the statistic applied to the true distribution $F$, over data sampled from $F$. Intuitively, we seek $Q_{S_N}(\alpha)$ such that $Q_{S_N}(\alpha) \geq Q_{S_N}^*(\alpha)$, since then we have that, if $F_0 = F$ then

$$\mathbb{P}\left(S_N \leq Q_{S_N}(\alpha)\right) \geq \mathbb{P}\left(S_N \leq Q_{S_N}^*(\alpha)\right) = 1 - \alpha,$$

so that the test has significance $\alpha$. Note that $\alpha$ significance is a property that holds for finite $N$.

In this paper we will use any one of three approaches depending on the test to computing thresholds $Q_{S_N}(\alpha)$:

1. If $Q_{S_N}^*(\alpha)$ is the same regardless the unknown, true distribution $F$, we can set $Q_{S_N}(\alpha) = Q_{S_N}^*(\alpha)$, which we compute for finite $N$ either exactly in closed-form or to arbitrary precision via simulation (with any choice of $F$). Test statistics with this property are called distribution free and the Kolmgogorov-Smirnov test (see below) is a canonical example.

2. If $Q_{S_N}^*(\alpha)$ does depend on $F$, we may still be able to find an upper bound. We take this approach in Theorem 16 in the appendix. While this approach guarantees $\alpha$-significance, it can be mathematically challenging and yield loose bounds.

3. If $Q_{S_N}^*(\alpha)$ does depend on $F$, we may choose $Q_{S_N}(\alpha)$ to approximate it using the bootstrap, where $Q_{S_N}(\alpha) \approx Q_{S_N}^*(\alpha)$ is approximated as the $(1-\alpha)^{\text{th}}$ quantile of the distribution of $S_N(\hat{F}_N, \tilde{\xi}^1, \ldots, \tilde{\xi}^N)$ over data sampled from $\hat{F}_N$, i.e., over draws of $N$ samples with replacement from the data $\xi^1, \ldots, \xi^N$. Under mild assumptions, bootstrap approximations improve as $N$ grows and are asymptotically exact (see [20]). Consequently, they are routinely used throughout statistics. Implementations of bootstrap procedures for computing thresholds $Q_{S_N}(\alpha)$ are available in many popular software packages, e.g., the function *one.boot* in the [R] package *simpleboot*.

One example of a GoF test is the Kolmogorov-Smirnov (KS) test for univariate distributions. The KS test uses the statistic

$$D_N = \max_{i=1,\ldots,N} \left\{ \max\left\{ \frac{i}{N} - F_0(\xi^{(i)}), F_0(\xi^{(i)}) - \frac{i-1}{N} \right\} \right\}.$$

Tables for $Q_{D_N}(\alpha)$ are widely available and can be computed for finite $N$ (see e.g. [14, 48]). In particular, under the assumption that $F$ is continuous, the KS statistic is distribution free and $Q_{D_N}(\alpha) = Q_{D_N}^*(\alpha)$ is easily computable by simulation (see Section 3.2), and if $F$ is not continuous then $Q_{D_N}(\alpha) \geq Q_{D_N}^*(\alpha)$ [36], ensuring $\alpha$ significance universally regardless of the unknown $F$.

---

[2] Such a hypothesis is called *simple*. By contrast, a composite hypothesis is *not* defined by a single distribution, but rather a family of distributions, and asserts that the data-generating distribution $F$ is some member of this family. An example of a composite hypothesis is that $F$ is normally distributed (with some, unknown mean and variance). We do not consider composite hypotheses in this work.

Fig. 1: A visualization of an example confidence region of the Kolmogorov-Smirnov test at significance 20%. The dashed curve is the true cumulative distribution function, that of a standard normal. The solid curve is the empirical cumulative distribution function having observed 100 draws from the true distribution. The confidence region contains all distributions with cumulative distribution functions that take inside the gray region.

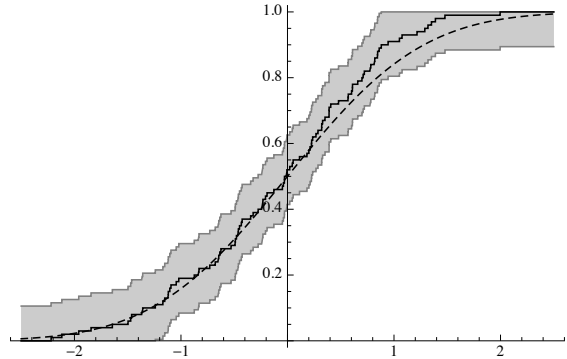

The set of all distributions $F_0$ that pass a test is called the *confidence region* of the test and is denoted by

$$\mathcal{F}^a_{S_N}(\xi^1, \ldots, \xi^N) = \left\{ F_0 \in \mathcal{P}(\Xi) : S_N(F_0, \xi^1, \ldots, \xi^N) \leq Q_{S_N}(\alpha) \right\}. \tag{7}$$

Note that in the above (7), $F_0$ is a dummy (distributional) variable. As an example, Figure 1 illustrates the confidence region of the KS test. Observe that by construction, the confidence region of a test with significance level $\alpha$ is a DUS which contains the true, unknown distribution $F$ with probability at least $1 - \alpha$ with respect to the distribution of the data drawn from $F$.

## 2.1 The Robust SAA approach

Given data $\xi^1, \ldots, \xi^N$, the Robust SAA approach involves the following steps:

---

**Robust SAA:**

1. Choose a significance level $0 < \alpha < 1$ and goodness-of-fit test at level $\alpha$ independently of the data.

2. Let $\mathcal{F} = \mathcal{F}_N(\xi^1, \ldots, \xi^N)$ be the confidence region of the test.

3. Solve

$$\overline{z} = \min_{x \in X} \sup_{F_0 \in \mathcal{F}_N(\xi^1, \ldots, \xi^N)} \mathbb{E}_{F_0}[c(x; \xi)]$$

   and let $\overline{x}$ be an optimal solution.

---

Section 5 illustrates how to solve the optimization problem in the last step for various choices of goodness-of-fit test and classes of cost functions.

## 2.2 Connections to existing methods

In Robust SAA, we use a GoF test at significance level $\alpha$ to construct a DUS that contains the true distribution $F$ with probability at least $1 - \alpha$ via its confidence region. It is possible to do the reverse as well; given a data-driven DUS $\mathcal{F}_N(\xi^1, \ldots, \xi^N)$ that contains the true distribution with probability at least $1 - \alpha$ with respect to the sampling distribution, we can construct a GoF test with significance level $\alpha$ that rejects the hypothesis (6) whenever $F_0 \notin \mathcal{F}_N(\xi^1, \ldots, \xi^N)$. This is often termed "the duality between hypothesis tests and confidence regions" (see for example §9.3 of [39]).

This reverse construction can be applied to existing data-driven DUSs in the literature such as [13, 15] to construct their corresponding hypothesis tests. In this way, hypothesis testing provides a common ground on which to understand and compare the methods.

In particular, the hypothesis tests corresponding to the DUSs of [13, 15] test only the first moments of the true distribution (cf. Section 4.3.2). By contrast, we will for the most part focus on tests (and corresponding confidence regions) that test the entire distribution, not just the first two moments. This feature is key to achieving both finite-sample and asymptotic guarantees at the same time.

## 3 Finite-sample performance guarantees

We first study the implication of a test's significance on the finite-sample performance of Robust SAA. Let us define the following random variables expressible as functions of the data $\xi^1, \ldots, \xi^N$:

The DRO solution:  $\qquad\qquad \overline{x} \in \arg\min_{x \in X} \sup_{F_0 \in \mathcal{F}_N(\xi^1, \ldots, \xi^N)} \mathbb{E}_{F_0}[c(x; \xi)].$

The DRO value:  $\qquad\qquad\quad \overline{z} = \min_{x \in X} \sup_{F_0 \in \mathcal{F}_N(\xi^1, \ldots, \xi^N)} \mathbb{E}_{F_0}[c(x; \xi)].$

The true cost of the DRO solution:  $\qquad z = \mathbb{E}_F[c(\overline{x}; \xi) | \xi^1, \ldots, \xi^N].$

The following is an immediate consequence of significance:

**Proposition 1** *If $\mathcal{F}_N(\xi^1, \ldots, \xi^N)$ is the confidence region of a valid GoF test at significance $\alpha$, then, with respect to the data sampling process,*

$$\mathbb{P}(\overline{z} \geq z) \geq 1 - \alpha.$$

*Proof* Suppose $F \in \mathcal{F}_N$. Then $\sup_{F_0 \in \mathcal{F}_N} \mathbb{E}_{F_0}[c(x; \xi)] \geq \mathbb{E}_F[c(x; \xi)]$ for any $x \in X$. Therefore, we have $\overline{z} \geq z$. In terms of probabilities, this implication yields,

$$\mathbb{P}(\overline{z} \geq z) \geq \mathbb{P}(F \in \mathcal{F}_N) \geq 1 - \alpha. \qquad\qquad \square$$

This makes explicit the connection between the statistical property of significance of a test with the objective performance of the corresponding Robust SAA decision in the full-information stochastic optimization problem.

### 3.1 Tests for distributions with known discrete support

When $\xi$ has known finite support $\Xi = \{\hat{\xi}^1, \ldots, \hat{\xi}^n\}$ there are two popular GoF tests: Pearson's $\chi^2$ test and the G-test (see [14]). Let $p(j) = F(\{\hat{\xi}^j\})$, $p_0(j) = F_0(\{\hat{\xi}^j\})$, and $\hat{p}_N(j) = \frac{1}{N} \sum_{i=1}^N \mathbb{I}\left[\xi^i = \hat{\xi}^j\right]$ be the true, hypothetical, and empirical probabilities of observing $\hat{\xi}^j$, respectively.

Pearson's $\chi^2$ test uses the statistic

$$X_N = \left( \sum_{j=1}^n \frac{(p_0(j) - \hat{p}_N(j))^2}{p_0(j)} \right)^{1/2},$$

whereas the G-test uses the statistic

$$G_N = \left( 2 \sum_{j=1}^n \hat{p}_N(j) \log\left(\frac{\hat{p}_N(j)}{p_0(j)}\right) \right)^{1/2}.$$

Thresholds for both Pearson's $\chi^2$ test and the G-test can be computed exactly via simulation (also known as Fisher's exact test); however, when $N\hat{p}_N(j) > 5 \; \forall j$, standard statistical practice is to use an asymptotic approximation equal to $1/\sqrt{N}$ times the $1 - \alpha$ quantile of the $\chi$ distribution with $n - 1$ degrees of freedom for both tests.[3]

The confidence regions of these statistics take the form of (7) for $S_N$ being either $X_N$ or $G_N$. An illustration of these ($n = 3$, $N = 50$) is given in Figure 2. Intuitively, these can be seen as generalized balls around the empirical distribution $\hat{p}_N$. The metric is given by the statistic $S_N$, and the radius diminishes as $O(N^{-1/2})$ (see [14]).

---

[3]  Note that standard definitions of the statistics have the form of $X_N^2$ and $G_N^2$ with thresholds given by the $\chi^2$ distribution. Our nonstandard but equivalent definition is so that thresholds have the rate $O(1/\sqrt{N})$ to match other tests presented.

Fig. 2: Distributional uncertainty sets (projected onto the first two components) for the discrete case with $n = 3$, $\alpha = 0.8$, $N = 50$, $p = (0.5, 0.3, 0.2)$. The dot denotes the true frequencies $p$ and the triangle the observed fractions $\hat{p}_{50}$.

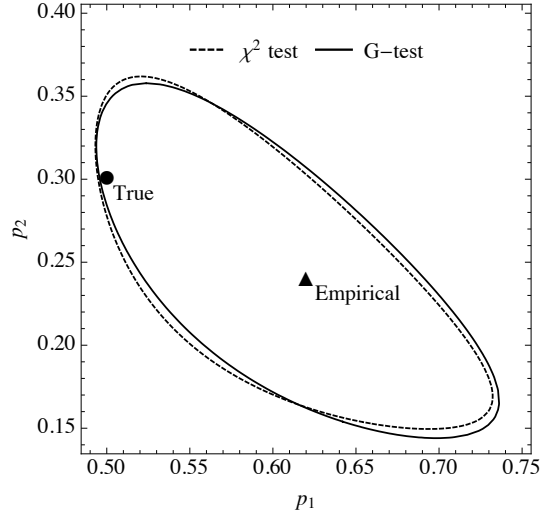

## 3.2 Tests for univariate distributions

Suppose $\xi$ is a univariate continuous random variable that is known to have lower support greater than $\underline{\xi}$ and upper support less than $\bar{\xi}$. These bounds could possibly be infinite. In this section, with $\xi$ univariate, we will use $(i)$ to denote the index of the observation that is the $i^{\text{th}}$ smallest so that $\xi^{(1)} \leq \cdots \leq \xi^{(N)}$. The most commonly used GoF tests in this setting are the Kolmogorov-Smirnov (KS) test, the Kuiper test, the Cramér-von Mises (CvM) test, the Watson test, and the Anderson-Darling (AD) test. The KS ($D_N$), the Kuiper ($V_N$), the CvM ($W_N$), the Watson ($U_N$), and the AD ($A_N$) tests use the statistics (see [14])

$$
D_N = \max_{i=1,\ldots,N} \left\{ \max \left\{ \frac{i}{N} - F_0(\xi^{(i)}), F_0(\xi^{(i)}) - \frac{i-1}{N} \right\} \right\},
$$
$$
V_N = \max_{1 \leq i \leq N} \left( F_0(\xi^{(i)}) - \frac{i-1}{N} \right) + \max_{1 \leq i \leq N} \left( \frac{i}{N} - F_0(\xi^{(i)}) \right),
$$
$$
W_N = \left( \frac{1}{12N^2} + \frac{1}{N} \sum_{i=1}^{N} \left( \frac{2i-1}{2N} - F_0(\xi^{(i)}) \right)^2 \right)^{1/2},
$$
$$
U_N = \left( W_N^2 - \left( \frac{1}{N} \sum_{i=1}^{N} F_0(\xi^{(i)}) - \frac{1}{2} \right)^2 \right)^{1/2},
$$
$$
A_N = \left( -1 - \sum_{i=1}^{N} \frac{2i-1}{N^2} \left( \log F_0(\xi^{(i)}) + \log(1 - F_0(\xi^{(N+1-i)})) \right) \right)^{1/2}.
$$

(8)

We let $S_N \in \{D_N, W_N, A_N, V_N, U_N\}$ be any one of the above statistics and $Q_{S_N}(\alpha)$ the corresponding threshold. Tables for $Q_{S_N}(\alpha)$ are widely available for finite $N$ (see [41, 48]). Alternatively, $Q_{S_N}(\alpha)$ can be computed to arbitrary precision for finite $N$ by simulation as the $(1-\alpha)^{\text{th}}$ percentile of the distribution of $S_N$ when $F_0(\xi^i)$ in (8) are replaced by IID uniform random variables on $[0, 1]$ (see [14]).

The confidence regions of these tests take the form of (7). Recall Figure 1 illustrated $\mathcal{F}_{D_N}^{\alpha}$. As in the discrete case, $\mathcal{F}_{S_N}^{\alpha}$ can also be seen as a generalized ball around the empirical distribution $\hat{F}_N$. Again, the radius diminishes as $O(N^{-1/2})$ (see [14]).

When $\underline{\xi}$ and $\bar{\xi}$ are finite, we take $\mathcal{F}_{S_N}^{\alpha}$ to be our DUS corresponding to these tests. When either $\underline{\xi}$ or $\bar{\xi}$ is infinite, however, $\bar{z}$ in (3) may also be infinite as seen in the following proposition.

**Proposition 2** *Fix $x$, $\alpha$, and $S_N \in \{D_N, W_N, A_N, V_N, U_N\}$. If $c(x; \xi)$ is continuous but unbounded on $\Xi$ then $\mathcal{C}(x; \mathcal{F}^\alpha_{S_N}) = \infty$ almost surely.*

The conditions of Proposition 2 are typical in many applications. For example, in §3 of [54], the authors briefly propose, but do not explore in detail, a data-driven DRO formulation that is equivalent to our Robust SAA formulation using the KS test with $\Xi = (-\infty, \infty)$. Over such an unbounded support, most cost functions used in practice will be unbounded over $\xi$ for any $x$, including, as an example, the newsvendor cost function. Proposition 2 implies that any such formulation will yield a trivial, infinite bound for any $x$ (and moreover there is no robust optimal $\bar{x}$), limiting the value of the approach.

Consequently, when either $\underline{\xi}$ or $\bar{\xi}$ is infinite, we will employ an alternative, non-standard, GoF test in Robust SAA. The confidence region of our proposed test will satisfy the conditions of Theorem 1, and, therefore, (3) will attain a finite, optimal solution.

Our proposed test combines one of the above GoF tests with a second test for a generalized moment of the distribution. Specifically, fix any function $\phi : \Xi \to \mathbb{R}_+$ such that $\mathbb{E}_F[\phi(\xi)] < \infty$ and $|c(x_0; \xi)| = O(\phi(\xi))$ for some $x_0 \in X$. For a fixed $\mu_0$, consider the null hypothesis

$$H'_0 : \mathbb{E}_F[\phi(\xi)] = \mu_0. \tag{9}$$

There are many possible hypothesis tests for (9). Any of these tests can be used as the second test in our proposal. For concreteness, we focus on a test that rejects (9) if

$$M_N = \left| \mu_0 - \frac{1}{N} \sum_{i=1}^N \phi(\xi^i) \right| > Q_{M_N}(\alpha). \tag{10}$$

In our numerical experiments in Section 7.1 we approximate $Q_{M_N}(\alpha)$ as $\hat{\sigma}_N Q_{T_{N-1}}(\alpha/2)/\sqrt{N}$ where $Q_{T_{N-1}}(\alpha/2)$ is the $(1 - \alpha/2)^{th}$ quantile of Student's T-distribution with $N - 1$ degrees of freedom and $\hat{\sigma}^2_N$ is the sample variance of $\phi(\xi)$. An alternative approach is to use the bootstrap.

Given $0 < \alpha_1, \alpha_2 < 1$, combining $S_N$ and (10), we propose the following GoF test:

$$\text{Reject } F_0 \text{ if either } S_N > Q_{S_N}(\alpha_1) \ \text{ or } \ \left| \mathbb{E}_{F_0}[\phi(\xi)] - \frac{1}{N} \sum_{i=1}^N \phi(\xi^i) \right| > Q_{M_N}(\alpha_2).$$

By the union bound, the probability of incorrectly rejecting $F_0$ is at most

$$\mathbb{P}(S_N > Q_{S_N}(\alpha_1)) + \mathbb{P}\left( \left| \mathbb{E}_{F_0}[\phi(\xi)] - \frac{1}{N} \sum_{i=1}^N \phi(\xi^i) \right| > Q_{M_N}(\alpha_2) \right) \le \alpha_1 + \alpha_2.$$

Thus, our proposed test has significance level $\alpha_1 + \alpha_2$.

The confidence region of the above test is given by the intersection of the confidence region of our original goodness-of-fit test and the confidence region of our test for (9):

$$\mathcal{F}^{\alpha_1, \alpha_2}_{S_N, M_N} = \mathcal{F}^{\alpha_1}_{S_N} \cap \mathcal{F}^{\alpha_2}_{M_N} = \left\{ F_0 \in \mathcal{P}(\Xi) : S_N \le Q_{S_N}(\alpha_1), \ \left| \mathbb{E}_{F_0}[\phi(\xi)] - \frac{1}{N} \sum_{i=1}^N \phi(\xi^i) \right| \le Q_{M_N}(\alpha_2) \right\}. \tag{11}$$

Observe that since $|c(x_0; \xi)| = O(\phi(\xi))$, i.e., $\exists \nu, \eta$ such that $|c(x_0; \xi)| \le \nu + \eta\phi(\xi)$, we have

$$\mathcal{C}(x_0; \mathcal{F}^{\alpha_1, \alpha_2}_{S_N, M_N}) = \sup_{F_0 \in \mathcal{F}^{\alpha_1, \alpha_2}_{S_N, M_N}} \mathbb{E}_{F_0}[c(x_0; \xi)] \le \nu + \eta \sup_{F_0 \in \mathcal{F}^{\alpha_1, \alpha_2}_{S_N, M_N}} \mathbb{E}_{F_0}[\phi(\xi)] \le \nu + \frac{\eta}{N} \sum_{i=1}^N \phi(\xi^i) + \eta Q_{M_N}(\alpha_2) < \infty,$$

so that unlike $\mathcal{F}^\alpha_{S_N}$, our new confidence region $\mathcal{F}^{\alpha_1, \alpha_2}_{S_N, M_N}$ does indeed satisfy the conditions of Theorem 1, even if $\underline{\xi}$ or $\bar{\xi}$ are infinite. This is critical to achieving useful solutions with finite guarantees and is an important distinction with the approach of [54].

### 3.3 Tests for multivariate distributions

In this section, we propose two different tests for the case $d \ge 2$. The first is a standard test based on testing marginal distributions. The second is a new test we propose that tests the full joint distribution.

### 3.3.1 Testing marginal distributions

Let $\alpha_1, \ldots, \alpha_d > 0$ be given such that $\alpha = \alpha_1 + \cdots + \alpha_d < 1$. Consider the test for the hypothesis $F = F_0$ that proceeds by testing the hypotheses $F_i = F_{0,i}$ for each $i = 1, \ldots, d$ by applying any test from Sections 3.1 and 3.2 at significance level $\alpha_i$ to the sample $\xi_i^1, \ldots, \xi_i^N$ and rejecting $F_0$ if any of these fail. The corresponding confidence region is

$$\mathcal{F}_{\text{marginals}}^\alpha = \left\{ F_0 \in \mathcal{P}(\Xi) \,:\, F_{0,i} \in \mathcal{F}_i^{\alpha_i}\left(\xi_i^1, \ldots, \xi_i^N\right) \; \forall i = 1, \ldots, d \right\},$$

where $\mathcal{F}_i^{\alpha_i}$ denotes the confidence region corresponding to the test applied on the $i^{\text{th}}$ component. By the union bound we have

$$\mathbb{P}\left(F \notin \mathcal{F}_{\text{marginals}}^\alpha\right) \leq \sum_{i=1}^d \mathbb{P}\left(F_i \notin \mathcal{F}_i^{\alpha_i}\right) \leq \sum_{i=1}^d \alpha_i = \alpha,$$

so the test has significance $\alpha$.

This test has the limitation that each $\alpha_i$ may need to be small to achieve a particular total significance $\alpha$, especially when $d$ is large. Moreover, two multivariate distributions may be distinct while having identical marginals, limiting the ability of this test to distinguish general distributions, which has adverse consequences for convergence as we will see in Section 4. The test we develop next does not suffer from either of these limitations.

### 3.3.2 Testing linear-convex ordering

In this section, we first provide some background on the linear-convex ordering (LCX) of random vectors first proposed in [43], and then use LCX to motivate a new GoF test for multivariate distributions. To the best of our knowledge, we are the first to propose GoF tests based on LCX.

Given two multivariate distributions $G$ and $G'$, we define

$$G \preceq_{\text{LCX}} G' \iff \mathbb{E}_G[\phi(a^T \xi)] \leq \mathbb{E}_{G'}[\phi(a^T \xi)] \;\; \forall \, a \in \mathbb{R}^d \text{ and convex functions } \phi : \mathbb{R} \to \mathbb{R} \tag{12}$$

$$\iff \mathbb{E}_G[\max\{a^T \xi - b, 0\}] \leq \mathbb{E}_{G'}[\max\{a^T \xi - b, 0\}] \;\; \forall \, |a_1| + \ldots + |a_d| + |b| \leq 1, \tag{13}$$

where the second equivalence follows from Theorem. 3.A.1 of [44].

Our interest in LCX stems from the following result. Theorem 1 of [43] stats that, assuming $\mathbb{E}_G[\|\xi\|_2^2] < \infty$,

$$\mathbb{E}_G[\xi_i^2] = \mathbb{E}_{G'}[\xi_i^2] \;\forall i = 1, \ldots, d \text{ and } G \preceq_{\text{LCX}} G' \implies G = G'.$$

Since the converse is clearly true and since $G \preceq_{\text{LCX}} G'$ implies that $\mathbb{E}_G[\xi_i^2] \leq \mathbb{E}_{G'}[\xi_i^2] \;\forall i = 1, \ldots, d$, we can rewrite this equivalently as

$$\mathbb{E}_G[\|\xi\|_2^2] \geq \mathbb{E}_{G'}[\|\xi\|_2^2] \text{ and } G \preceq_{\text{LCX}} G' \iff G = G', \tag{14}$$

where $\|\xi\|_2^2 = \sum_{i=1}^d \xi_i^2$. Equation (14) motivates our GoF test. Intuitively, the key idea of our test is that if $F \neq F_0$, i.e., we should reject $F_0$, then by (14) either $\mathbb{E}_{F_0}[\|\xi\|_2^2] < \mathbb{E}_F[\|\xi\|_2^2]$ or $F_0 \not\preceq_{LCX} F$. Thus, we can create a GoF test by testing for each of these cases separately.

More precisely, for a fixed $\mu_0$, first consider the hypothesis

$$H_0' : \mathbb{E}_F[\|\xi\|_2^2] = \mu_0. \tag{15}$$

As in Section 3.2, there are many possible tests for (15). For concreteness, we focus on a one-sided test (i.e., against the alternative $H_1' : \mathbb{E}_F[\|\xi\|_2^2] > \mu_0$) which rejects (15) if $R_N = \frac{1}{N} \sum_{i=1}^N \|\xi^i\|_2^2 - \mu_0 > Q_{R_N}(\alpha)$. where $Q_{R_N}(\alpha)$ is a threshold which can be computed by bootstrapping.

Next, define the statistic

$$C_N(F_0) = \sup_{|a_1| + \cdots + |a_d| + |b| \leq 1} \left( \mathbb{E}_{F_0}[\max\{a^T \xi - b, 0\}] - \frac{1}{N} \sum_{i=1}^N [\max\{a^T \xi^i - b, 0\}] \right).$$

From (13), $C_N(F_0) \leq 0 \iff F_0 \preceq_{LCX} \hat{F}_N$. (Recall that $\hat{F}_N$ denotes the empirical distribution.)

Finally, combining these pieces and given $0 < \alpha_1, \alpha_2 < 1$, our LCX-based GoF test is

$$\text{Reject } F_0 \text{ if either } C_N(F_0) > Q_{C_N}(\alpha_1) \text{ or } \mathbb{E}_{F_0}[\|\xi\|_2^2] < \frac{1}{N} \sum_{i=1}^N \|\xi^i\|_2^2 - Q_{R_N}(\alpha_2). \tag{16}$$

In Section 10.3 in the appendix, we provide an exact closed-form formula for a valid threshold $Q_{C_N}(\alpha_1)$. Approximate thresholds can also be computed via the bootstrap. We use the bootstrap in our numerical experiments in Section 7.5 because, in this case, it provides much tighter thresholds for the sizes $N$ studied. The particular bootstrap procedure is provided in Section 10.3.

From a union bound we have that the LCX-based GoF test (16) has significance level $\alpha_1 + \alpha_2$. The confidence region of the LCX-based GoF test is

$$\mathcal{F}_{C_N, R_N}^{\alpha_1, \alpha_2} = \left\{ F_0 \in \mathcal{P}(\Xi) : C_N(F_0) \leq Q_{C_N}(\alpha_1), \; \mathbb{E}_{F_0}[\|\xi\|_2^2] \geq \frac{1}{N} \sum_{i=1}^{N} \|\xi^i\|_2^2 - Q_{R_N}(\alpha_2) \right\}. \tag{17}$$

## 4 Convergence

In this section, we study the relationship between the GoF test underlying an application of Robust SAA and convergence properties of the Robust SAA optimal values $\overline{z}$ and solutions $\overline{x}$. Recall from Section 2.2 that since many existing DRO formulations can be recast as confidence regions of hypothesis tests, our analysis will simultaneously also allow us to study the convergence properties of these methods as well.

The convergence conditions we seek are

Convergence of objective function:    $\mathcal{C}(x; \mathcal{F}_N) \to \mathbb{E}_F[c(x; \xi)]$    (18)
                                       uniformly over any compact subset of $X$,

Convergence of optimal values:        $\min_{x \in X} \mathcal{C}(x; \mathcal{F}_N) \to \min_{x \in X} \mathbb{E}_F[c(x; \xi)],$    (19)

Convergence of optimal solutions:     Every sequence $x_N \in \arg\min_{x \in X} \mathcal{C}(x; \mathcal{F}_N)$ has at least one limit    (20)
                                       point, and all of its limit points are in $\arg\min_{x \in X} \mathbb{E}_F[c(x; \xi)],$

all holding almost surely (a.s.). The key to these will be a new, restricted form of statistical consistency of GoF tests that we introduce and term uniform consistency.

### 4.1 Uniform consistency and convergence of optimal solutions

In statistics, *consistency of a test* (see Def. 2 below) is a well-studied property that a GoF test may exhibit. In this section, we define a new property of GoF tests that we call *uniform consistency*. Uniform consistency is a strictly stronger property than consistency, in the sense that every uniformly consistent test is consistent, but some consistent tests are not uniformly consistent. More importantly, we will prove that uniform consistency of the underlying GoF test tightly characterizes when conditions (18)-(20) hold. In particular, we show that when $X$ and $\Xi$ are bounded, uniform consistency of the underlying test implies conditions (18)-(20) for any cost function $c(x; \xi)$ which is equicontinuous in $x$, and if the test is not uniformly consistent, then there exist cost functions (equicontinuous in $x$) for which conditions (18)-(20) do not hold. When $X$ or $\Xi$ are unbounded, the same conclusions hold for all cost functions which are equicontinuous in $x$ and satisfy an additional, mild regularity condition. (See Theorem 2 for a precise statement.) In other words, we can characterize the convergence of Robust SAA and other data-driven, DRO formulations by studying if their underlying GoF test is uniformly consistent. In our opinion, these results highlight a new, fundamental connection between statistics and data-driven optimization. We will use this result to assess the strength of various DRO formulations for certain applications in what follows.

First, we recall the definition of consistency of a GoF test (cf. entry for *consistent test* in [17]):

**Definition 2** A GoF test is *consistent* if, for any distribution $F$ generating the data and any $F_0 \neq F$, the probability of rejecting any $F_0$ approaches 1 as $N \to \infty$. That is, $\lim_{N \to \infty} \mathbb{P}(F_0 \notin \mathcal{F}_N) = 1$ whenever $F \neq F_0$.

Observe

**Proposition 3** *If a test is consistent, then any $F_0 \neq F$ is a.s. rejected infinitely often (i.o.) as $N \to \infty$.*

*Proof*          $\mathbb{P}(F_0 \text{ rejected i.o.}) = \mathbb{P}\left(\limsup_{N \to \infty} \{F_0 \notin \mathcal{F}_N\}\right) \geq \limsup_{N \to \infty} \mathbb{P}(F_0 \notin \mathcal{F}_N) = 1,$

where the first inequality follows from Fatou's Lemma, and the second equality holds since the test is consistent.                                                                                                      □

Consistency describes the test's behavior with respect to a single, fixed null hypothesis $F_0$. In particular, the conclusion of Proposition 3 holds only when we consider the same, fixed distribution $F_0$ for each $N$. We would like to extend consistency to describe the test's behavior with respect to many possible $F_0$ simultaneously. We define uniform consistency by requiring that a condition similar to the conclusion of Proposition 3 hold for a sequence of distributions sufficiently different from $F$:

**Definition 3** A GoF test is *uniformly consistent* if, for any distribution $F$ generating the data, every sequence $F_N$ that does not converge weakly to $F$ is rejected i.o, a.s., that is,

$$\mathbb{P}\left(F_N \not\to F \implies F_N \notin \mathcal{F}_N \text{ i.o.}\right) = 1.$$

The requirement that $F_N$ does not converge weakly to $F$ parallels the requirement that $F_0 \neq F$.

Uniform consistency is a strictly stronger requirement than consistency.

**Proposition 4** *If a test is uniformly consistent, then it is consistent. Moreover, there exist tests which are consistent, but not uniformly consistent.*

Uniform consistency is the key property for the convergence of Robust SAA. Besides uniform consistency, convergence will be contingent on three assumptions.

*Assumption 1* $c(x; \xi)$ is equicontinuous in $x$ over all $\xi \in \Xi$.

*Assumption 2* $X$ is closed and either

a. $X$ is bounded or

b. $\lim_{||x|| \to \infty} c(x; \xi) = \infty$ uniformly over $\xi$ in some $D \subseteq \Xi$ with $F(D) > 0$ and $\liminf_{||x|| \to \infty} \inf_{\xi \notin D} c(x; \xi) > -\infty$.

*Assumption 3* Either

a. $\Xi$ is bounded or

b. $\exists \phi : \Xi \to \mathbb{R}_+$ such that $\sup_{F_0 \in \mathcal{F}_N} \left| \mathbb{E}_{F_0} \phi(\xi) - \frac{1}{N} \sum_{i=1}^{N} \phi(\xi^i) \right| \to 0$ almost surely and $c(x; \xi) = O(\phi(\xi))$ for each $x \in X$.

Assumptions 1 and 2 are only slightly stronger than those required for the existence of an optimal solution in Theorem 1. Assumption 1 is the same. Assumption 2 is the same when $X$ is bounded, otherwise, requiring uniformity in coerciveness over some small set with nonzero measure. The second portion of Assumption 2b is trivially satisfied by cost functions which are bounded from below, such as nonnegative. Finally, observe that in the case that $\Xi$ is unbounded, our proposed DUS in (11) satisfies Assumption 3b by construction.

Under these assumptions, the following theorem provides a tight characterization of convergence.

**Theorem 2** $\mathcal{F}_N$ *is the confidence region of a uniformly consistent test if and only if Assumptions 1, 2, and 3 imply that conditions (18)-(20) hold a.s.*

Table 1: Summary of convergence results.
* denotes the result is tight in the sense that there are examples in this class that do not converge.
† denotes that a transformation of the data may be necessary; see Section 4.3.1.

| GoF test | Support | Consistent | Uniformly consistent | (18)-(20) hold a.s. for any $c(x; \xi)$ that is | | |
|---|---|---|---|---|---|---|
| | | | | equicontinuous in $x$ | separable as in (21) | factorable as in (26) |
| $\chi^2$ and G-test | Finite | Yes | Yes | Yes | Yes | Yes |
| KS, Kuiper, CvM, Watson, and AD tests | Univariate | Yes | Yes | Yes | Yes | Yes |
| Test of marginals using the above tests | Multivariate | No | No | No* | Yes | Yes† |
| LCX-based test | Multivariate, bounded | Yes | Yes | Yes | Yes | Yes |
| LCX-based test | Multivariate, unbounded | Yes | ? | ? | ? | ? |
| Tests implied by DUSs of [13, 15] | Multivariate | No | No | No* | No* | Yes |

Thus, in one direction, we can guarantee convergence (i.e., conditions (18)-(20) hold a.s.) if Assumptions 1, 2, and 3 are satisfied and we use a uniformly consistent test in applying Robust SAA. In the other direction, if we use a test that is not uniformly consistent, there will exist instances satisfying Assumptions 1, 2, and 3 for which convergence fails.

Some of the GoF tests in Section 3 are not consistent, and therefore, cannot be uniformly consistent. By Theorem 2, DROs built from these tests cannot exhibit asymptotic convergence for all cost functions. One might argue, then, that that these DRO formulations should be avoided in modeling and applications in favor of DROs based on uniformly consistent tests.

In most applications, however, we are not concerned with asymptotic convergence *for all* cost functions, but rather only for the *given* cost function $c(x; \xi)$. It may happen a DRO may exhibit asymptotic convergence for this particular cost function, even when its DUS is given by the confidence region of an inconsistent test. (We will see an example of this behavior with the multi-item newsvendor problem in Section 7.4.)

To better understand when this convergence may occur despite the fact that the test is not consistent, we introduce a more relaxed form of uniform consistency.

**Definition 4** Given $c(x; \xi)$, we say that $F_N$ *c-converges* to $F$ if $\mathbb{E}_{F_N}[c(x; \xi)] \to \mathbb{E}_F[c(x; \xi)]$ for all $x \in X$.

**Definition 5** Given $c(x; \xi)$, a GoF test is *c-consistent* if, for any distribution $F$ generating the data, every sequence $F_N$ that does not *c-converge* to $F$ is rejected i.o., a.s., that is,

$$\mathbb{P}\left(F_N \text{ does not } c\text{-converge to } F \implies F_N \notin \mathcal{F}_N \text{ i.o.}\right) = 1.$$

This notion may potentially be weaker than consistency, but is sufficient for convergence for a given instance as shown below.

**Theorem 3** *Suppose Assumptions 1 and 3 hold and that $\mathcal{F}_N$ always contains the empirical distribution. If $\mathcal{F}_N$ is the confidence region of a c-consistent test, then conditions (18)-(20) hold a.s.*

In the next sections we will explore the consistency of the various tests introduced in Section 3. We summarize our results in Table 1.

## 4.2 Tests for distributions with discrete or univariate support

All of the classical tests we considered in Section 3 are uniformly consistent.

**Theorem 4** *Suppose $\Xi$ has known discrete support. Then, the $\chi^2$ and G-tests are uniformly consistent.*

**Theorem 5** *Suppose $\Xi$ is univariate. Then, the KS, Kuiper, CvM, Watson, and AD tests are uniformly consistent.*

## 4.3 Tests for multivariate distributions

### 4.3.1 Testing marginal distributions

We first claim that the test of marginals is not consistent. Indeed, consider a multivariate distribution $F_0 \neq F$ which has the same marginal distributions, but a different joint distribution. By construction, the probability of rejecting $F_0$ is at most $\alpha$ for all $N$, and hence does not converge to 1. Since the test of marginals is not consistent, it cannot be uniformly consistent.

We next show that the test is, however, *c-consistent* whenever the cost is separable over the components of $\xi$.

**Proposition 5** *Suppose $c(x; \xi)$ is separable over the components of $\xi$, that is, can be written as*

$$c(x; \xi) = \sum_{i=1}^{d} c_i(x; \xi_i), \tag{21}$$

*and Assumptions 1, 2, and 3 hold for each $c_i(x; \xi_i)$. Then, the test of marginals is c-consistent if each univariate test is uniformly consistent.*

That is to say, if the cost can be separated as in (21), applying the tests from Section 3.2 to the marginals is sufficient to guarantee convergence.

It is important to note that some cost functions may only be separable after a transformation of the data, potentially into a space of different dimension. If that is the case, we may transform $\xi$ and apply the tests to the transformed components in order to achieve convergence.

### 4.3.2 Tests implied by DUSs of [13, 15]

The DUS of [15] has the form

$$
\mathcal{F}_{\mathrm{DY},N}^{\alpha} = \left\{ F_0 \in \mathcal{P}(\Xi) \ : \ \begin{matrix} \left(\mathbb{E}_{F_0}[\xi] - \hat{\mu}_N\right)^T \hat{\Sigma}_N^{-1} \left(\mathbb{E}_{F_0}[\xi] - \hat{\mu}_N\right) \le \gamma_{1,N}(\alpha), \\ \gamma_{3,N}(\alpha)\hat{\Sigma}_N \preceq \mathbb{E}_{F_0}[(\xi - \hat{\mu}_N)\,(\xi - \hat{\mu}_N)^T] \preceq \gamma_{2,N}(\alpha)\hat{\Sigma}_N \end{matrix} \right\}
\tag{22}
$$

where $\hat{\mu}_N = \dfrac{1}{N}\sum_{i=1}^{N}\xi^i, \quad \hat{\Sigma}_N = \dfrac{1}{N}\sum_{i=1}^{N}\left(\xi^i - \hat{\mu}_N\right)\left(\xi^i - \hat{\mu}_N\right)^T.$  (23)

The thresholds $\gamma_{1,N}(\alpha)$, $\gamma_{2,N}(\alpha)$ $\gamma_{3,N}(\alpha)$ are developed therein (for $\Xi$ bounded) so as to guarantee a significance of $\alpha$ (in our GoF interpretation) and, in particular, have the property that

$$
\gamma_{1,N}(\alpha) \downarrow 0,\ \gamma_{2,N}(\alpha) \downarrow 1,\ \gamma_{3,N}(\alpha) \uparrow 1.
\tag{24}
$$

Seen from the perspective of GoF testing, valid thresholds can also be approximated by the bootstrap, as discussed in Section 2. The resulting threshold will be significantly smaller (see Section 7) but exhibit the same asymptotics.

The DUS of [13] has the form

$$
\mathcal{F}_{\mathrm{CEG},N}^{\alpha} = \left\{ F_0 \in \mathcal{P}(\Xi) \ : \ \begin{matrix} ||\mathbb{E}_{F_0}[\xi] - \hat{\mu}_N||_2 \le \gamma_{1,N}(\alpha), \\ \left|\left|\mathbb{E}_{F_0}\left[(\xi - \mathbb{E}_{F_0}[\xi])\,(\xi - \mathbb{E}_{F_0}[\xi])^T\right] - \hat{\Sigma}_N\right|\right|_{\mathrm{Frobenius}} \le \gamma_{2,N}(\alpha) \end{matrix} \right\}.
$$

The thresholds $\gamma_{1,N}(\alpha)$, $\gamma_{2,N}(\alpha)$ are developed in [47] (for $\Xi$ bounded) so as to guarantee a significance of $\alpha$ and with the property that

$$
\gamma_{1,N}(\alpha) \downarrow 0,\ \gamma_{2,N}(\alpha) \downarrow 0.
\tag{25}
$$

Again, seen from the perspective of GoF testing, valid thresholds can also be approximated by the bootstrap, having similar asymptotics.

The GoF tests implied by these DUSs consider only the first two moments of a distribution (mean and covariance). Therefore, the probability of rejecting a multivariate distribution different from the true one but with the same mean and covariance is by construction never more than $\alpha$, instead of converging to 1. That is, these tests are not consistent and therefore they are not uniformly consistent. We next provide conditions on the cost function that guarantee that the tests are nonetheless $c$-consistent.

**Proposition 6** *Suppose $c(x;\xi)$ can be written as*

$$
c(x;\xi) = c_0(x) + \sum_{i=1}^{d} c_i(x)\xi_i + \sum_{i=1}^{d}\sum_{j=1}^{i} c_{ij}(x)\xi_i\xi_j
\tag{26}
$$

*and that $\mathbb{E}_F[\xi_i\xi_j]$ exists. Then, the tests with confidence regions given by $\mathcal{F}_{DY,N}^{\alpha}$ or $\mathcal{F}_{CEG,N}^{\alpha}$ are c-consistent.*

Note that because we may transform the data to include components for each pairwise multiplication, the conditions on the cost function in Proposition 6 are stronger than those in Proposition 5. In particular, in one dimension $d = 1$, every cost function is separable but not every cost function satisfies the decomposition (26).

### 4.3.3 Testing linear-convex ordering

The previous two multivariate GoF tests were neither consistent, nor uniformly consistent. By contrast,

**Theorem 6** *The LCX-based test is consistent.*

**Theorem 7** *Suppose $\Xi$ is bounded. Then the LCX-based test is uniformly consistent.*

It is an open question whether the LCX-based test is uniformly consistent for unbounded $\Xi$ (Theorem 6 proves it is consistent). We conjecture that it is. Moreover, in our numerical experiments involving the LCX test, we have observed convergence of the Robust SAA solutions to the full-information optimum even when $\Xi$ is unbounded. (See Section 7.5 for an example.)

## 5 Tractability

In this section, we characterize conditions under which problem (3) is theoretically tractable, i.e., can be solved with a polynomial-time algorithm. Additionally, we are interested in cases where (3) is practically tractable, i.e., can be solved using off-the-shelf linear or second-order cone optimization solvers. In some specific cases, we show that Robust SAA using the KS test admits a closed-form solution.

### 5.1 Tests for distributions with known discrete support

In the case of known discrete support, we considered two GoF tests: Pearson's $\chi^2$ test, corresponding to the DUS $\mathcal{F}^\alpha_{X_N}$, and the G-test, corresponding to the DUS $\mathcal{F}^\alpha_{G_N}$. In this section we present a reformulation of (3) for these DUSs as a single-level convex optimization problem, from which tractability results will follow.

The DUSs of these GoF tests are a special case of those considered in [3]. As corollaries of the results therein we have the following:

**Theorem 8** *Under the assumptions of Theorem 1, we have*

$$
\mathcal{C}\left(x; \mathcal{F}^\alpha_{X_N}\right) = \min_{r,s,t,y,c} \quad r + \left(Q_{X_N}(\alpha)\right)^2 s - \sum_{j=1}^n \hat{p}_N(j)t_j
$$
$$
s.t. \quad r \in \mathbb{R},\ s \in \mathbb{R}_+,\ t \in \mathbb{R}^n,\ y \in \mathbb{R}^n_+,\ c \in \mathbb{R}^n
$$
$$
s + r \geq c_j \qquad\qquad\qquad \forall j = 1, \ldots, n
$$
$$
2s + t_j \leq y_j \qquad\qquad\qquad \forall j = 1, \ldots, n
$$
$$
(2s - c_j + r,\ r - c_j,\ y_j) \in C^3_{SOC} \qquad \forall j = 1, \ldots, n
$$
$$
c_j \geq c\left(x; \hat{\xi}^j\right) \qquad\qquad\qquad \forall j = 1, \ldots, n
$$

$$
\mathcal{C}\left(x; \mathcal{F}^\alpha_{G_N}\right) = \min_{r,s,t,c} \quad r + \frac{1}{2}\left(Q_{G_N}(\alpha)\right)^2 s - \sum_{j=1}^n \hat{p}_N(j)t_j
$$
$$
s.t. \quad r \in \mathbb{R},\ s \in \mathbb{R}_+,\ t \in \mathbb{R}^n,\ c \in \mathbb{R}^n
$$
$$
(t_j,\ s,\ s + r - c_j) \in C_{XC} \qquad\qquad \forall j = 1, \ldots, n \qquad\qquad (27)
$$
$$
c_j \geq c\left(x; \hat{\xi}^j\right) \qquad\qquad\qquad \forall j = 1, \ldots, n
$$

*where* $C^3_{SOC} = \left\{(x,y,z) \in \mathbb{R}^3 : x \geq \sqrt{y^2 + z^2}\right\}$ *is the three-dimensional second-order cone and* $C_{XC} =$ *closure* $\left(\{(x,y,z) : ye^{x/y} \leq z,\ y > 0\}\right)$ *is the exponential cone.*

The DRO problem (3) is $\min_{x \in X} \mathcal{C}(x; \mathcal{F})$. Therefore, for $\mathcal{F}^\alpha_{X_N}$ and $\mathcal{F}^\alpha_{G_N}$, (3) can be formulated as a single-level optimization problem by augmenting the corresponding minimization problem above with the decision variable $x \in X$. Note that apart from the constraints $x \in X$ and

$$
c_j \geq c\left(x; \hat{\xi}^j\right), \qquad\qquad\qquad\qquad (28)
$$

the rest of the constraints, as seen in the problems in Theorem 8, are convex. The following result characterizes in general when solving these problems is tractable in a theoretical sense.

**Theorem 9** *Suppose that* $X \subseteq \mathbb{R}^{d_x}$ *is a closed convex set for which a weak separation oracle is given and that*

$$
c\left(x; \hat{\xi}^j\right) = \max_{k=1,\ldots,K_j} c_{jk}(x)
$$

*where each* $c_{jk}(x)$ *is a convex function in* $x$ *for which evaluation and subgradient oracles are given. Then, under the assumptions of Theorem 1, we can find an* $\epsilon$*-optimal solution to (3) in the discrete case for* $S_N = X_N,\ G_N$ *in time and oracle calls polynomial in* $n, d_x, K_1, \ldots, K_n, \log(1/\epsilon)$.

For some problems the constraints $x \in X$ and (28) can also be conically formulated as the Example 1 below shows. In such a case, the DRO can be solved directly as a conic optimization problem. Optimization over the exponential cone – a non-symmetric cone — although theoretically tractable, is numerically challenging. Fortunately, the particular exponential cone constraints (27) can be recast as second-order cone constraints, albeit with constraint complexity growing in both $n$ and $N$ (see [33]).

*Example 1 Two-stage problem with linear recourse and a non-increasing, piece-wise-linear convex disutility.* Consider the following problem

$$c(x; \hat{\xi}^j) = \max_{k=1,\dots,K} \left( \gamma_k (c^T x + R_j(x)) + \beta_k \right), \quad \gamma_k \leq 0 \tag{29}$$

$$\text{where} \quad R_j(x) = \min_{y \in \mathbb{R}_+^{d_y}} \ f_j^T y$$

$$\text{s.t.} \ A_j x + B_j y = b_j$$

$$X = \{x \geq 0 : Hx = h\}.$$

This problem was studied in a non-data-driven DRO settings in [6, 19, 56]. To formulate (3), we may introduce variables $y \in \mathbb{R}_+^{n \times d_y}$ and replace (28) with

$$c_j \geq \gamma_k \left( c^T x + f_j^T y_j \right) + \beta_k \qquad\qquad \forall j = 1, \dots, n, \ \forall k = 1, \dots, K,$$

$$A_j x + B_j y_j = b_j \qquad\qquad\qquad\quad \forall j = 1, \dots, n.$$

The resulting problem is then a second-order cone optimization problem for $\mathcal{F}_{X_N}^\alpha$ and $\mathcal{F}_{G_N}^\alpha$.

## 5.2 Tests for univariate distributions

We now consider the case where $\xi$ is a general univariate random variable. We proceed by reformulating (3) as a single-level optimization problem by leveraging semi-infinite conic duality. This leads to corresponding tractability results. In the following we will use the notation $\xi^{(0)} = \underline{\xi}$ and $\xi^{(N+1)} = \overline{\xi}$. Recall that these may be infinite.

The first observation is that the constraint $S_N(\zeta_1, \dots, \zeta_N) \leq Q_{S_N}(\alpha)$ is convex in $\zeta_i = F_0(\xi^{(i)})$ and representable using *canonical cones*. By a canonical cone, we mean any Cartesian product of the cones $\mathbb{R}^k$, $\{0\}$, $\mathbb{R}_+^k$ (positive orthant), $C_{\text{SOC}}^k$ (second-order cone), and semidefinite cone. Optimization over canonical cones is tractable both theoretically and practically using state-of-the-art interior point algorithms [4].

**Theorem 10** *For each of $S_N \in \{D_N, V_N, W_N, U_N, A_N\}$*

$$S_N(\zeta_1, \dots, \zeta_N) \leq Q_{S_N}(\alpha) \quad \Longleftrightarrow \quad A_{S_N} \zeta - b_{S_N, \alpha} \in K_{S_N}$$

*for convex cones $K_{S_N}$, matrices $A_{S_N}$, and vectors $b_{S_N, \alpha}$ as follows:*

$$K_{D_N} = \mathbb{R}_+^{2N}, \quad b_{D_N, \alpha} = \begin{pmatrix} \frac{1}{N} - Q_{D_N}(\alpha) \\ \vdots \\ \frac{N}{N} - Q_{D_N}(\alpha) \\ -\frac{0}{N} - Q_{D_N}(\alpha) \\ \vdots \\ -\frac{N-1}{N} - Q_{D_N}(\alpha) \end{pmatrix}, \quad A_{D_N} = \begin{pmatrix} [I_N] \\ [-I_N] \end{pmatrix},$$

$$K_{V_N} = \left\{ (x, y) \in \mathbb{R}^{2N} : \min_i x_i + \min_i y_i \geq 0 \right\}, \quad b_{V_N, \alpha} = \begin{pmatrix} \frac{1}{N} - Q_{V_N}(\alpha)/2 \\ \vdots \\ \frac{N}{N} - Q_{V_N}(\alpha)/2 \\ -\frac{0}{N} - Q_{V_N}(\alpha)/2 \\ \vdots \\ -\frac{N-1}{N} - Q_{V_N}(\alpha)/2 \end{pmatrix}, \quad A_{V_N} = \begin{pmatrix} [I_N] \\ [-I_N] \end{pmatrix},$$

$$K_{W_N} = C_{\text{SOC}}^{N+1}, \quad b_{W_N, \alpha} = \begin{pmatrix} \sqrt{N \left( Q_{W_N}(\alpha) \right)^2 - \frac{1}{2N}} \\ \frac{1}{2N} \\ \frac{3}{2N} \\ \vdots \\ \frac{2N-1}{2N} \end{pmatrix}, \quad A_{W_N} = \begin{pmatrix} 0 \cdots 0 \\ [I_N] \end{pmatrix},$$

$$K_{U_N} = C_{\mathrm{SOC}}^{N+2}, \quad b_{U_N,\alpha} = \begin{pmatrix} \frac{-1}{2} + \left( \frac{N}{24} - \frac{N}{2} \left( Q_{U_N}(\alpha) \right)^2 \right) \\ \frac{-1}{2} - \left( \frac{N}{24} - \frac{N}{2} \left( Q_{U_N}(\alpha) \right)^2 \right) \\ 0 \\ \vdots \\ 0 \end{pmatrix}, \quad A_{U_N} = \begin{pmatrix} \frac{1-N}{2N} & \frac{3-N}{2N} & \cdots & \frac{N-1}{2N} \\ \frac{N-1}{2N} & \frac{N-3}{2N} & \cdots & \frac{1-N}{2N} \\ & & & \\ [I_N - \frac{1}{N} E_N] & & & \end{pmatrix},$$

$$K_{A_N} = \left\{ (z,x,y) \in \mathbb{R} \times \mathbb{R}_+^{2N} \; : \; |z| \leq \prod_{i=1}^{N} (x_i y_i)^{\frac{2i-1}{2N^2}} \right\}, \quad b_{A_N,\alpha} = \begin{pmatrix} e^{-\left( Q_{A_N}(\alpha) \right)^2 - 1} \\ 0 \\ \vdots \\ 0 \\ -1 \\ \vdots \\ -1 \end{pmatrix}, \quad A_{A_N} = \begin{pmatrix} 0 \cdots 0 \\ [I_N] \\ [-\tilde{I}_N] \end{pmatrix},$$

where $I_N$ is the $N \times N$ identity matrix, $\tilde{I}_N$ is the skew identity matrix ($[\tilde{I}_N]_{ij} \equiv \mathbb{I}\,[i = N - j]$), and $E_N$ is the $N \times N$ matrix of all ones.

Note that the cones $K_{D_N}, K_{W_N}, K_{U_N}$ are canonical cones. The other cones can be expressed using canonical cones. The cone $K_{V_N}$ is an orthogonal projection of an affine slice of $\mathbb{R}^{2N+2} \times \mathbb{R}_+^3$. The cone $K_{A_N}$ is an orthogonal projection of an affine slice of the product of $2^{\lceil \log_2(2N^2) \rceil + 1} - 2 = O(N^2)$ three-dimensional second-order cones (see [33]). Therefore, the constraint $A_{S_N}\zeta - b_{S_N,\alpha} \in K_{S_N}$ can be expressed using polynomially-sized canonical cones in each case.

Problem (3) is a two-level optimization problem. To formulate it as a single-level problem, we dualize the inner problem, $\mathcal{C}(x; \mathcal{F})$. For a cone $K \subseteq \mathbb{R}^k$, we use the notation $K^*$ to denote the dual cone $K^* = \{y \in \mathbb{R}^k : y^T z \geq 0 \; \forall z \in K\}$. The following is a direct consequence of Theorem 10 and semi-infinite conic duality, namely Proposition 3.4 of [45].[4]

**Theorem 11** *Let $S_N \in \{D_N, V_N, W_N, U_N, A_N\}$. Under the assumptions of Theorem 1,*

$$\mathcal{C}\left(x; \mathcal{F}_{S_N}^\alpha\right) = \min_{r,c} \quad b_{S_N,\alpha}^T r + c_{N+1}$$
$$\text{s.t} \quad -r \in K_{S_N}^*, \; c \in \mathbb{R}^{N+1}$$
$$\left( A_{S_N}^T r \right)_i = c_i - c_{i+1} \qquad\qquad \forall i = 1, \ldots, N$$
$$c_i \geq \sup_{\xi \in (\xi^{(i-1)}, \xi^{(i)})} c(x; \xi) \qquad\qquad \forall i = 1, \ldots, N+1 \qquad (30)$$

**Theorem 12** *Let $S_N \in \{D_N, V_N, W_N, U_N, A_N\}$. Under the assumptions of Theorem 1,*

$$\mathcal{C}\left(x; \mathcal{F}_{S_N, M_N}^{\alpha_1, \alpha_2}\right) = \min_{r,t,s,c} \quad b_{S_N,\alpha_1}^T r + c_{N+1} + \left( \hat{\mu} + Q_{M_N}^{\alpha_2} \right) t - \left( \hat{\mu} - Q_{M_N}^{\alpha_2} \right) s$$
$$\text{s.t} \quad -r \in K_{S_N}^*, \; t \geq 0, \; s \geq 0, \; c \in \mathbb{R}^{N+1}$$
$$\left( A_{S_N}^T r \right)_i = c_i - c_{i+1} \qquad\qquad \forall i = 1, \ldots, N$$
$$c_i \geq \sup_{\xi \in (\xi^{(i-1)}, \xi^{(i)})} \left( c(x; \xi) - (t - s)\phi(\xi) \right) \qquad\qquad \forall i = 1, \ldots, N+1. \qquad (31)$$

Recall that when Proposition 2 applies, i.e., when either $\underline{\xi}$ or $\overline{\xi}$ are infinite and the cost may be unbounded, then an additional test to bound $\mathbb{E}_F[\phi(\xi)]$ is necessary in order to satisfy the assumptions of Theorem 1 and get a finite solution. It is in this case, that Theorem 12 is relevant.

Note that the cones $K_{D_N}, K_{W_N}, K_{U_N}$ are self-dual ($K^* = K$) and therefore the dual cones remain canonical cones. For $K_{V_N}$ and $K_{A_N}$, the dual cones are

$$K_{V_N}^* = \left\{ (x,y) \in \mathbb{R}_+^{2N} \; : \; \sum_{i=1}^{N} x_i = \sum_{i=1}^{N} y_i \right\}$$

$$K_{A_N}^* = \left\{ (z,x,y) \; : \; (z/\gamma, x, y) \in K_{A_N} \right\} \text{ where } \gamma = \prod_{i=1}^{d} \left( \frac{2i-1}{2N^2} \right)^{\frac{2i-1}{N^2}},$$

---

[4]  The only nuance is that Proposition 3.4 of [45] requires a generalized Slater point. We use the empirical distribution function, $\hat{F}_N$, as the generalized Slater point in the space of distributions.

and therefore they remain expressible using canonical cones.

Note that in the case of $\mathcal{F}_{S_N}^\alpha$, the worst-case distribution has discrete support on no more than $N + 1$ points. This is because shifting probability mass inside the interval $(\xi^{(i-1)}, \xi^{(i)}]$ does not change any $F_0(\xi^{(i)})$. In the worst-case, all mass in the interval (if any) will be placed on the point in the interval with the largest cost (including the left endpoint in the limit).

The DRO problem (3) is $\min_{x \in X} \mathcal{C}(x; \mathcal{F})$. Therefore, for $\mathcal{F}_{S_N}^\alpha$ and $\mathcal{F}_{S_N, M_N}^{\alpha_1, \alpha_2}$, (3) can formulated as a single-level optimization problem by augmenting the corresponding minimization problem above with the decision variable $x \in X$. We next give general conditions that ensure the theoretical tractability of the problem.

**Theorem 13** *Suppose that $X \subseteq \mathbb{R}^{d_x}$ is a closed convex set for which a weak separation oracle is given and that*

$$c(x; \xi) = \max_{k=1,\ldots,K} c_k(x; \xi) \tag{32}$$

*where each $c_k(x; \xi)$ is convex in $x$ for each $\xi$ and continuous in $\xi$ for each $x$ and for which an oracle is given for the subgradient in $x$. If $\mathcal{F} = \mathcal{F}_{S_N}^\alpha$, suppose also that an oracle is given for maximizing $c_k(x; \xi)$ over $\xi$ in any closed (possibly infinite) interval for fixed $x$. If $\mathcal{F} = \mathcal{F}_{S_N, M_N}^{\alpha_1, \alpha_2}$, suppose also that an oracle is given for maximizing $c_k(x; \xi) + \eta \phi(\xi)$ over $\xi$ in a closed interval for fixed $x$ and $\eta \in \mathbb{R}$. Then, under the assumptions of Theorem 1, we can find an $\epsilon$-optimal solution to (3) in time and oracle calls polynomial in $N, d_x, K, \log(1/\epsilon)$ for $\mathcal{F} = \mathcal{F}_{S_N}^\alpha$ or $\mathcal{F} = \mathcal{F}_{S_N, M_N}^{\alpha_1, \alpha_2}$.*

As in the discrete case, when the constraints $x \in X$ and (30) (or, (31)) can be conically formulated, Theorem 11 (or, Theorem 12, respectively) provides an explicit single-level conic optimization formulation of the problem (3). In Examples 2, 3, and 4 below, we consider specific problems for which this is the case and study this formulation.

*Example 2 The newsvendor problem.* In the newsvendor problem, one orders in advance $x \geq 0$ units of a product to satisfy an unknown future demand for $\xi \geq 0$ units. Unmet demand is penalized by $b > 0$, representing either backlogging costs or lost profit. Left over units are penalized by $h > 0$, representing either holding costs or recycling costs. The cost function is therefore $c(x; \xi) = \max\{b(\xi - x), h(x - \xi)\}$, the lower support of $\xi$ is $\underline{\xi} \geq 0$, and the space of decisions is $X = \mathbb{R}_+$. For the bounded-support case $\overline{\xi} < \infty$, the constraints (30) become

$$c_i \geq b(\xi^{(i)} - x), \quad c_i \geq h(x - \xi^{(i-1)}) \qquad\qquad \forall i = 1, \ldots, N + 1$$

and $x \in X$ becomes $x \in \mathbb{R}_+$. In the unbounded case $\overline{\xi} = \infty$, we may use $\phi(\xi) = |\xi|$ in the construction of (11). Because $\underline{\xi} \geq 0$, we have $\phi(\xi) = \xi$. The constraints (31) then become

$$c_i \geq b(\xi^{(i)} - x) - (t - s)\xi^{(i)}, \quad c_i \geq h(x - \xi^{(i-1)}) - (t - s)\xi^{(i-1)} \qquad \forall i = 1, \ldots, N + 1$$

where the $(N + 1)^{\text{st}}$ left constraint is equivalent to $b \leq t - s$ because $\xi^{(N+1)} = \infty$. Substituting these constraints in this way the DRO (3) becomes a conic optimization problem.

In the specific case of bounded support and $\mathcal{F} = \mathcal{F}_{D_N}^\alpha$ this reformulation yields a linear optimization problem, which admits a closed-form solution given next.

**Proposition 7** *Suppose that $\Xi = [\underline{\xi}, \overline{\xi}]$ is compact, and $N$ is large enough so that $Q_{D_N}(\alpha) < \frac{\min\{b, h\}}{b + h}$. Then, the DRO (3) for the newsvendor problem with $\mathcal{F} = \mathcal{F}_{D_N}^\alpha$ admits the closed-form solution:*

$$\overline{x} = (1 - \theta)\xi^{(i_{lo})} + \theta\xi^{(i_{hi})}$$

$$\overline{z} = \frac{1}{N} \sum_{1 \leq i \leq i_{lo} \lor i_{hi} \leq i \leq N} c\left(\overline{x}; \xi^{(i)}\right) + Q_{D_N}(\alpha) c\left(\overline{x}; \underline{\xi}\right) + Q_{D_N}(\alpha) c\left(\overline{x}; \overline{\xi}\right)$$

$$\quad - \left(\frac{\lceil N(\theta - Q_{D_N}(\alpha)) \rceil}{N} - (\theta - Q_{D_N}(\alpha))\right) c\left(\overline{x}; \xi^{(i_{lo})}\right) - \left((\theta + Q_{D_N}(\alpha)) - \frac{\lfloor N(\theta + Q_{D_N}(\alpha)) \rfloor}{N}\right) c\left(\overline{x}; \xi^{(i_{hi})}\right)$$

*where $\theta = b/(b + h)$, $i_{lo} = \lceil N(\theta - Q_{D_N}(\alpha)) \rceil$, and $i_{hi} = \lfloor N(\theta + Q_{D_N}(\alpha)) \rfloor + 1 \rfloor$.*

Importantly, this means that solving the Robust SAA newsvendor problem is no more difficult than solving the SAA newsvendor problem.

*Example 3 Max of bilinear functions.* More generally, we may consider cost functions of the form (32) with bilinear parts $c_k(x; \xi) = p_{k0} + p_{k1}^T x + p_{k2}\xi + \xi p_{k3}^T x$. In this case, (30) is equivalent to

$$c_i \geq p_{k0} + p_{k1}^T x + p_{k2}\xi^{(i-1)} + \xi^{(i-1)} p_{k3}^T x, \qquad \forall i = 1, \ldots, N, \quad \forall k = 1, \ldots, K \qquad (33)$$

$$c_i \geq p_{k0} + p_{k1}^T x + p_{k2}\xi^{(i)} + \xi^{(i)} p_{k3}^T x, \qquad \forall i = 1, \ldots, N, \quad \forall k = 1, \ldots, K. \qquad (34)$$

If the cost is fully linear, $p_{k3} = 0$ (as in the case of the newsvendor example), then (30) can be written in one linear inequality:

$$c_i \geq p_{k0} + p_{k1}^T x + \max\left\{ p_{k2}\xi^{(i-1)}, p_{k2}\xi^{(i)} \right\} \qquad \forall i = 1, \ldots, N, \quad \forall k = 1, \ldots, K. \qquad (35)$$

For $\mathcal{F} = \mathcal{F}_{S_N, M_N}^{\alpha_1, \alpha_2}$ we may use $\phi(\xi) = |\xi|$ and simply add $\left|\xi^{(i-1)}\right|$ and $\left|\xi^{(i)}\right|$ to the left-hand sides of (33) and (34), respectively, or to the corresponding branches of the max in (35).

*Example 4 Two-stage problem.* Consider a two-stage problem similar to the one studied in Example 1:

$$c(x; \xi) = \max_{k=1, \ldots, K} \left( \gamma_k(c^T x + R(x; \xi)) + \beta_k \right), \quad \gamma_k \leq 0 \qquad (36)$$

$$\text{where } R(x; \xi) = \min_{y \in \mathbb{R}_+^{d_y}} (f + g\xi)^T y$$

$$\text{s.t. } Ax + By = b + p\xi$$

$$X = \{x \geq 0 : Hx = h\}.$$

When only the right-hand-side vector is uncertain ($g = 0$), the recourse $R(x; \xi)$ is convex in $\xi$ so that the supremum in (30) is taken at one of the endpoints and we may use a similar construction as in Example 3.

When only the cost vector is uncertain ($p = 0$), the recourse $R(x; \xi)$ is concave in $\xi$. By linear optimization duality we may reformulate (30) by introducing variables $R \in \mathbb{R}^{N+1}$, $y \in \mathbb{R}_+^{d_y \times (N+1)}$, $\eta \in \mathbb{R}_+^{N+1}$, $\theta \in \mathbb{R}_+^{N+1}$ and constraints

$$c_i \geq \gamma_k(c^T x + R_i) + \beta_k \qquad \forall i = 1, \ldots, N+1, \quad \forall i = k, \ldots, K$$

$$\eta_i - \theta_i = f^T y_i, \ Ax + By_i = b \qquad \forall i = 1, \ldots, N+1$$

$$R_i \geq g^T y_i + \xi^{(i)}\eta_i - \xi^{(i-1)}\theta_i \qquad \forall i = 1, \ldots, N+1.$$

## 5.3 Tests for multivariate distributions

### 5.3.1 Testing marginal distributions

Recall that when $c(x; \xi)$ is separable over the components of $\xi$, i.e.,

$$c(x; \xi) = \sum_{i=1}^{d} c_i(x; \xi_i),$$

Robust SAA using the test of marginals converges to the full-information optimum (cf. Section 4.3.1). We next show that Robust SAA is also tractable in this case. When $\mathcal{F} = \mathcal{F}_{\text{marginals}}^{\alpha}$ and costs are separable, (3) can be written as

$$\min_{x \in X} \sup_{F_0 \in \mathcal{F}} \mathbb{E}_{F_0}[c(x; \xi)] = \min_{x \in X} \sum_{i=1}^{d} \sup_{F_{0,i} \in \mathcal{F}_i^{\alpha_i}} \mathbb{E}_{F_{0,i}}[c_i(x; \xi_i)].$$

Applying Theorems 8, 11, and 12 separately to each of these $d$ subproblems yields a single-level optimization problem (the theorem applied in each corresponds to the DUS $\mathcal{F}_i^{\alpha_i}$ used and whether it has discrete or general support). This problem is theoretically tractable when each subproblem satisfies the corresponding conditions in Theorems 9 and 13. Similarly, when each subproblem is of one of the forms treated in Examples 1, 2, 3, and 4, (3) can be formulated as a linear or second-order cone optimization problem.

*5.3.2 Testing linear-convex ordering*

Next, we consider the case of the test based on LCX. We assume that we can represent the cost function as

$$c(x;\xi) = \max_{k=1,\ldots,K} c_k(x;\xi). \tag{37}$$

where each $c_k(x;\xi)$ is concave in $\xi$.

**Theorem 14** *Suppose that we can express $c(x;\xi)$ as in (37) with each $c_k(x;\xi)$ closed concave in $\xi$. Let $c_{k*}(x;\zeta) = \inf_{\xi \in \Xi}\left(\zeta^T \xi - c_k(x;\xi)\right)$ be $c_k$'s concave conjugate and $\mathcal{G} = \{0,1\}^K \backslash \{(0,\ldots,0)\}$. Let*

$$
\begin{aligned}
\mathcal{C}'(x;\nu,\eta) = \min_{\substack{r,s,z,y,y' \\ w,w',f,g,h}} \quad & \nu f + \eta r + \sum_{\gamma \in \mathcal{G}} \sum_{i=1}^N z_{\gamma,i} \\
\text{s.t} \quad & r \in \mathbb{R},\, z \in \mathbb{R}_+^{|\mathcal{G}| \times N},\, y \in \mathbb{R}_+^{|\mathcal{G}| \times d},\, y' \in \mathbb{R}_+^{|\mathcal{G}|},\, w \in \mathbb{R}_+^{|\mathcal{G}| \times d},\, w' \in \mathbb{R}_+^{|\mathcal{G}|} \\
& f \in \mathbb{R}_+,\, g \in \mathbb{R}^K,\, h \in \mathbb{R}^{K \times d} \\
& g_k - \sum_{\gamma \in \mathcal{G}} \gamma_k w'_\gamma + r \geq 0 && \forall k = 1,\ldots,K \\
& \frac{1}{N}\left(w_\gamma^T \xi^i - w'_\gamma\right) \leq z_{\gamma,i} && \forall \gamma \in \mathcal{G},\, i = 1,\ldots,N \\
& w_\gamma \leq y_\gamma,\, -w_\gamma \leq y_\gamma,\, w'_\gamma \leq y'_\gamma,\, -w'_\gamma \leq y'_\gamma && \forall \gamma \in \mathcal{G} \\
& \sum_{\gamma \in \mathcal{G}}\left(\sum_{j=1}^d y_{\gamma,j} + y'_\gamma\right) \leq f \\
& h_k = \sum_{\gamma \in \mathcal{G}} \gamma_k w_\gamma && \forall k = 1,\ldots,K \\
& g_k \leq c_{k*}(x;h_k) && \forall k = 1,\ldots,K. \tag{39}
\end{aligned}
$$

*Then, for any $\xi' \in \Xi$, under the assumptions of Theorem 1, we have*

$$\mathcal{C}'(x;Q_{C_N}(\alpha_1),1) \geq \mathcal{C}\left(x;\mathcal{F}_{C_N,R_N}^{\alpha_1,\alpha_2}\right) \geq \mathcal{C}'\left(x;Q_{C_N}(\alpha_1) - \frac{||\xi'||_2 + 1}{||\xi'||_2^2}Q_{R_N}(\alpha_2), 1 - \frac{1}{||\xi'||_2^2}Q_{R_N}(\alpha_2)\right) - \frac{c(x;\xi')}{||\xi'||_2^2}.$$

*Moreover, if for a given $x$ there exists a sequence $\xi'_i \in \Xi$ such that*

$$\lim_{i \to \infty} ||\xi'_i|| = \infty, \quad \limsup_{i \to \infty} \frac{c(x;\xi'_i)}{||\xi'_i||_2^2} \geq 0, \tag{40}$$

*i.e., $\Xi$ is unbounded and the negative part of the cost function is sub-quadratic in $\xi$, then, under the assumptions of Theorem 1, we have*

$$\mathcal{C}\left(x;\mathcal{F}_{C_N,R_N}^{\alpha_1,\alpha_2}\right) = \mathcal{C}'(x;Q_{C_N}(\alpha_1),1).$$

The first part of the Theorem 14 gives a finite linear optimization problem to approximate $\mathcal{C}\left(x;\mathcal{F}_{C_N,R_N}^{\alpha_1,\alpha_2}\right)$, bounding it from above, and a range of similar problems – one for each $\xi' \in \Xi$ – that bound it from below and quantify the quality of approximation. The upper bound and lower bound problems differ only in their objective coefficients and a constant term. The objective coefficients differ by a quantity that shrinks with $||\xi'||_2$ suggesting that the approximation becomes better the larger $\Xi$ is. The second part of Theorem 14 says that if $\Xi$ is in fact unboundedly large, i.e., there exist sequences of $\xi'$ with magnitude approaching infinity, and if there is one such sequence that does not grow the constant term, then the reformulation is in fact exact. Note that any cost function that is bounded below (e.g., nonnegative) is trivially sub-quadratic in $\xi$ and always satisfies this second condition with any sequence. Similarly, linear or bi-linear cost functions are sub-quadratic in both the negative and positive directions and satisfy this condition.

Notice that in the case where (40) holds, the reformulation of $\mathcal{C}\left(x;\mathcal{F}_{C_N,R_N}^{\alpha_1,\alpha_2}\right)$ (i.e., $\mathcal{C}'(x;Q_{C_N}(\alpha_1),1)$) does not explicitly involve $\alpha_2$ – the significance of the test for $\mathbb{E}\left[||\xi||_2^2\right]$. This a consequence of the structure of the sub-quadratic cost function and the unbounded support, which implies that the lower bound on $\mathbb{E}\left[||\xi||_2^2\right]$ is not active for the worst-case distribution. This is similar to an often observed phenomenon where a worst-case distribution over a DUS restricting only mean and variance always attains maximal (and not minimal) variance for many optimization problems (see e.g. [15, 42]). The implication is that we may let $\alpha_2 \to 0$, increasing the probability of the finite-sample guarantee without affecting the solution $\overline{x}$ or the bound $\overline{z}$. This is the approach we take in the empirical study in Section 7.5 where we effectively set $\alpha_2 = 0$ and ignore the test for a lower bound on $\mathbb{E}\left[||\xi||_2^2\right]$.

The DRO problem (3) is $\min_{x \in X} \mathcal{C}(x; \mathcal{F})$. Therefore, for $\mathcal{F}_{C_N, R_N}^{\alpha_1, \alpha_2}$, (3) can be formulated as a single-level optimization problem by augmenting the minimization problem given by $\mathcal{C}'(x; Q_{C_N}(\alpha_1), 1)$ with the decision variable $x \in X$ when (40) holds. When (40) does not hold, this still yields a bound – i.e., the optimizer $\tilde{x}$ satisfies $\mathcal{C}'\left(\tilde{x}; Q_{C_N}(\alpha_1), 1\right) \geq \mathcal{C}\left(\tilde{x}; \mathcal{F}_{C_N, R_N}^{\alpha_1, \alpha_2}\right) \geq \mathcal{C}\left(\overline{x}; \mathcal{F}_{C_N, R_N}^{\alpha_1, \alpha_2}\right)$, which means it inherits $\overline{x}$'s performance guarantee. Apart from the constraint $x \in X$ and the constraint (39), the optimization problem that results is linear. The following provides general conditions that ensure theoretical tractability.

**Theorem 15** *Suppose that $X \subseteq \mathbb{R}^{d_x}$ is a closed convex set for which a weak separation oracle is given, that we can express $c(x; \xi)$ as in (37) with each $c_k(x; \xi)$ closed concave in $\xi$ for each $x$ and convex in $x$ for each $\xi$, and that oracles for subgradient in $x$ and concave conjugate in $\xi$ (minimizing $\zeta^T \xi - c(x; \xi)$ for any $\zeta$) are given. Then we can find an $\epsilon$-optimal solution to*

$$\min_{x \in X} \mathcal{C}'(x; \nu, \eta)$$

*in time and oracle calls polynomial in $N, d_x, 2^K, \log(1/\epsilon)$ for $\mathcal{F} = \mathcal{F}_{C_N, R_N}^{\alpha_1, \alpha_2}$.*

Note that reformulation size, and accordingly the complexity above, grows exponentially with $K$, i.e., as $2^K - 1$. At the same time, however, each piece $c_k(x; \xi)$ may be a general convex-concave function. Consequently, for many examples of interest, as detailed below, $K$ may be 1 or 2. For these, and other examples, the reformulation becomes a single linear optimization problem, suitable for off-the-shelf solvers.

*Example 5 Bilinear cost pieces.* One example of (37) with convex-concave pieces is the max of bilinear pieces:

$$c_k(x; \xi) = p_{k0} + p_{k1}^T x + p_{k2}^T \xi + x^T P_k \xi. \tag{41}$$

If $\Xi = \{\xi : B\xi \geq b\}$ is polyhedral the corresponding concave conjugate is

$$c_{k_*}(x; \zeta) = -p_{k0} - p_{k1}^T x + \max_{\rho \geq 0, \, B^T \rho = \zeta - p_{k2} - P_k^T x} b^T \rho,$$

and if $\Xi = \mathbb{R}^d$ this simplifies to

$$c_{k_*}(x; \zeta) = \begin{cases} -p_{k0} - p_{k1}^T x & \zeta = p_{k2} + P_k^T x \\ -\infty & \text{otherwise} \end{cases}.$$

In either case, we can represent (39) with a small number of linear inequalities and variables. If $X$ is also polyhedral, the DRO problem (3) becomes a single-level linear optimization problem.

*Example 6 portfolio allocation.* An important special case of the above is that of portfolio allocation. Consider $d$ securities with unknown future returns $\xi_i$. We must divide our budget into fractions $x_i \geq 0$ invested in security $i$ with $\sum_i x_i = 1$. The return on a unit budget is $x^T \xi$.

A popular metric for the quality of a portfolio is its conditional value at risk (CVaR). The CVaR at level $\epsilon$ of a portfolio [40] is defined as

$$\text{CVaR}_\epsilon(x) = \inf_{\beta \in \mathbb{R}} \mathbb{E}\left[\beta + \frac{1}{\epsilon} \max\{-x^T \xi - \beta, 0\}\right].$$

Under continuity assumptions, $\text{CVaR}_\epsilon(x)$ is the conditional expectation of negative returns given that returns are below the $\epsilon$-quantile. We can formulate the min CVaR problem using (41) with only $K = 2$ pieces. We augment the decision vector as $(\beta, x)$ and write the cost function as

$$c((\beta, x); \xi) = \max\left\{\beta, \, (1 - 1/\epsilon)\beta - \xi^T x/\epsilon\right\}. \tag{42}$$

As in Example 5, if $X$ is also polyhedral, the DRO problem (3) becomes a linear optimization problem. Since $K = 2$, the size of the problem only grows with $N$ and the dimension of $X$.

*Example 7 Two-stage problem.* Consider a two-stage linear optimization problem similar to the one studied in Example 4:

$$c(x;\xi) = c^T x + \min_{y \in \mathbb{R}_+^{d_y}} \ (f + G\xi)^T y$$

$$\text{s.t. } A_1 x + A_2 y = a,$$

where $\Xi = \{B\xi \geq b\}$, $b \in \mathbb{R}^m$. Then $c(x;\xi)$ is convex-concave and therefore representable as in (37) with $K = 1$ convex-concave pieces. By linear optimization duality, its concave conjugate is

$$c_*(x;\zeta) = -c^T x + \max_{y \in \mathbb{R}_+^{d_y},\, z \in \mathbb{R}^m} \ b^T z - f^T y$$

$$\text{s.t. } A_1 x + A_2 y = a,$$
$$B^T z + G^T y = \zeta.$$

Hence, if $X$ is also polyhedral, the DRO problem (3) becomes a linear optimization problem. Since $K = 1$, the size of the problem only grows with $N$ and the dimension of $X$.

## 6 Estimating the price of data

Our framework allows one to compute the price one would be willing to pay for further data gathering. Given the present dataset, we define the price of data (PoD) as follows:

$$\text{PoD} = \overline{z}\left(\xi^1, \ldots, \xi^N\right) - \mathbb{E}\left[\overline{z}\left(\xi^1, \ldots, \xi^N, \xi^{N+1}\right) \Big| \xi^1, \ldots, \xi^N\right]. \tag{43}$$

PoD is equal to the expected marginal benefit of one additional data point in reducing our bound on costs.

One way to estimate the above quantity is via resampling:

$$\text{PoD} \approx \overline{z}\left(\xi^1, \ldots, \xi^N\right) - \frac{1}{N}\sum_{i=1}^N \overline{z}\left(\xi^1, \ldots, \xi^N, \xi^i\right). \tag{44}$$

The resampled average can also be, in turn, estimated by an average over a smaller random subsample from the data. This approach is illustrated numerically in Section 7.5.

In the case of the newsvendor problem using the KS test, the closed form solution yields a simpler approximation. Observe that in Proposition 7, small changes to the data change $\overline{x}$ very little and the costs for $\xi$ near $\overline{x}$ (in particular, between $i_{\text{lo}}$ and $i_{\text{hi}}$) are small compared to costs far away from $\overline{x}$. Thus, we suggest the approximation

$$\text{PoD} \approx \left(Q_{D_N}(\alpha) - Q_{D_{N+1}}(\alpha)\right)\left(c\left(\overline{x};\underline{\xi}\right) + c\left(\overline{x};\overline{\xi}\right)\right). \tag{45}$$

This approximation is illustrated numerically in Section 7.1.

We can write a more explicit approximation using the asymptotic approximation of $Q_{D_N}(\alpha)$ (see [49]) and $1/\sqrt{N} - 1/\sqrt{N+1} \approx 1/(2N^{3/2})$ for large $N$:

$$\text{PoD} \approx \frac{q_\alpha}{2N^{3/2}}\left(c\left(\overline{x};\underline{\xi}\right) + c\left(\overline{x};\overline{\xi}\right)\right) \quad \text{where} \quad q_\alpha = \begin{cases} 1.36, & \alpha = 0.05, \\ 1.22, & \alpha = 0.1, \\ 1.07, & \alpha = 0.2. \end{cases}$$

## 7 Empirical study

We now turn to an empirical study of Robust SAA as applied to specific problems in inventory and portfolio management. The cost functions are all specific cases of the examples studied in Section 5. Recall that in these examples the resulting formulations were all linear and second-order cone optimization problems.

Before proceeding to the details of our experiments, we summarize our main findings.

– Robust SAA yields stable, low-variance solutions, in contrast to SAA and 2-SAA. Indeed, for small to moderate $N$, Robust SAA performs comparably to both methods in terms of expected out-of-sample performance, but significantly outperforms both in terms of out-of-sample variability.

Fig. 3: The PDFs of the demand distributions considered for the newsvendor problem.

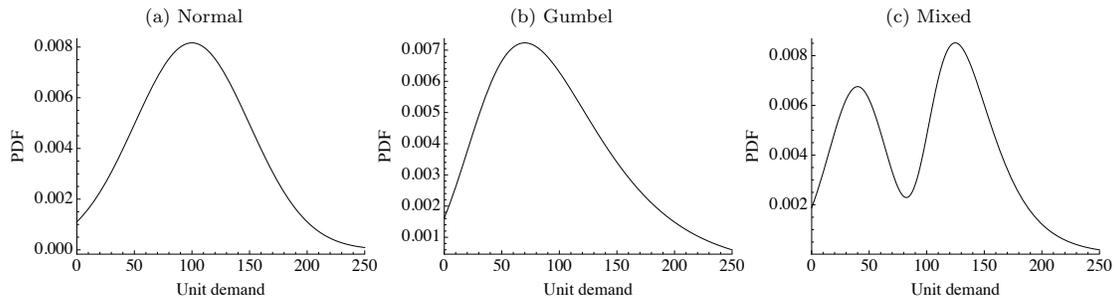

– In agreement with our theoretical findings, Robust SAA's bounds hold with probability $1 - \alpha$, even for small values of $N$. By contrast, for some problems, 2-SAA may yield invalid bounds that do *not* hold with probability $1 - \alpha$, even for very large $N$.
– Robust SAA typically produces significantly tighter bounds than the method of [15], especially for small to moderate $N$.
– The above observations regarding Robust SAA seem to hold generally for many choices of GoF test, so long as this test is $c$−consistent for the relevant problem.

### 7.1 Single-item newsvendor and the KS test

We begin with an application to the classic newsvendor problem with a bounded, continuous demand distribution, as studied in Example 2, with the KS GoF test. As noted in Proposition 7, Robust SAA admits a closed-form solution for $N$ sufficiently large, and can be solved as a small linear optimization problem otherwise. We compare its solution to SAA, 2-SAA, and the method of [15] (which in this case can also be solved closed-form). We use a significance level of $\alpha = 20\%$ (i.e., 80% confidence) for all methods. All optimization problems are solved using GUROBI 5.5 [24].

We consider a 95% service-level requirement ($b = 19$, $h = 1$) and each of the following distributions, all truncated to lie between 0 and 250:

1. Normal distribution with mean 100 and standard deviation 50.
2. Right-skewed Gumbel distribution with location 70 and scale $30/\gamma$ (the Euler constant).
3. Mixture model of 40% normal with mean 40 and standard deviation 25 and 60% right-skewed Gumbel with location 125 and scale $15/\gamma$.

For reference, we plot their PDFs in Figure 3.

In Figure 4 we focus on the results for the bounded normal distribution. (The results for the other distributions are shown in Figure 6b.) In Panel 4a, we present the full-information optimum (1), SAA estimates (2), 2-SAA upper-bound (5), the bound of [15] (labeled "Delage & Ye '10 DRO Bound" in the figure), and Robust SAA upper bounds $\overline{z}$ as a function of $N$. For reference, we also present the non-data-driven bound of [42] that requires knowledge of the true mean and true standard deviation (labeled "Scarf '58 DRO Bound" in the figure). The center lines represent the mean over several draws of the data and error bars represent the $20^{\text{th}}$ and $80^{\text{th}}$ percentiles over the draws of the data. The number of draws for each $N$ was chosen so that the mean has a standard error of no more than 0.1 standard deviations – this ensures that estimation errors are imperceptible in the figure.

We first observe that both the SAA estimates and 2-SAA bounds converge to the full-information optimum $z_{\text{stoch}}$ as expected. However, the SAA exhibits a strong bias; it frequently provides estimates that are below the full information optimum (65% of the time for $N = 100$) which is an unattainable performance. This phenomenon is not unique to this problem. The SAA estimate for any fixed $N$ will always have a

Fig. 4: Comparison of SAA estimates, 2-SAA guarantees, data-driven DRO of [15], bound of [42] and Robust SAA guarantees. All data-driven methods use a significance of $\alpha = 20\%$. Note the log scales.

(a) Estimates and bounds: center lines are means over replicates and error bars are the $20^{\text{th}}$ and $80^{\text{th}}$ percentiles.

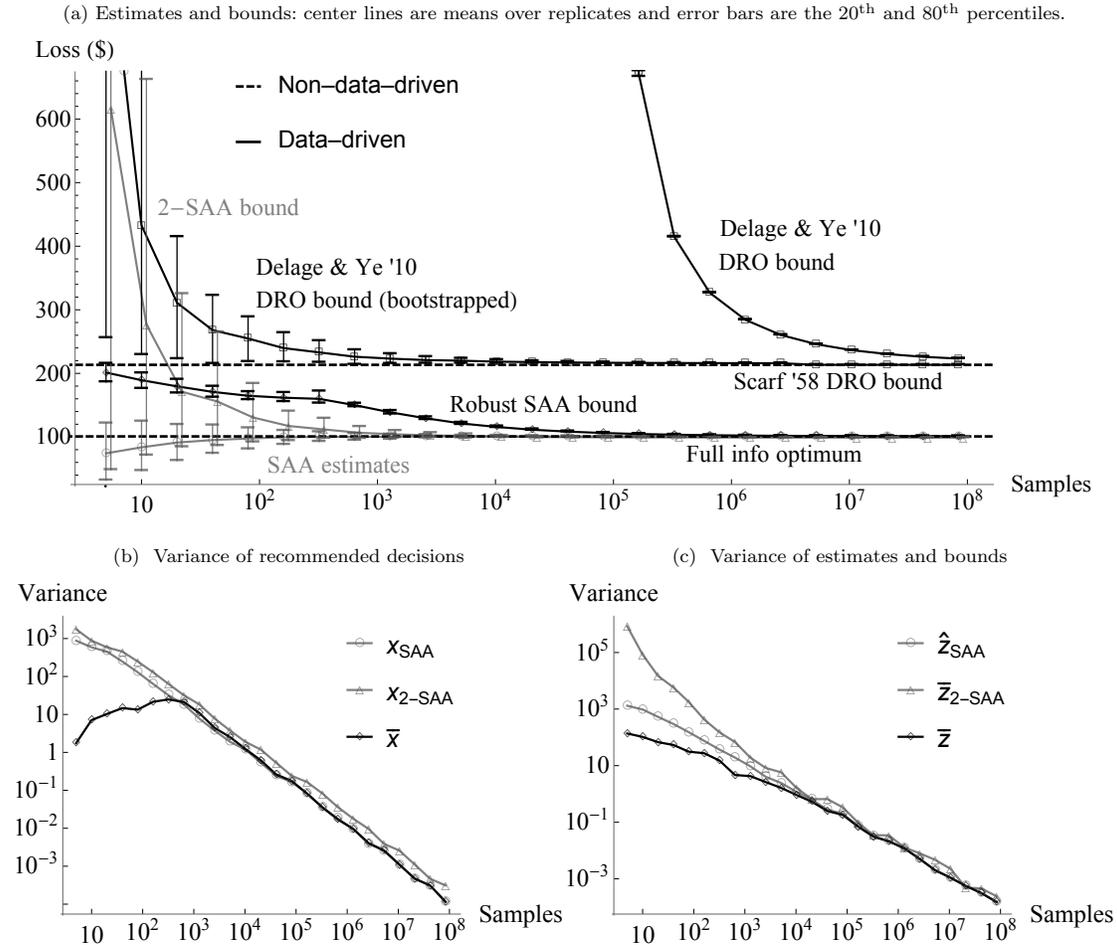

(b) Variance of recommended decisions     (c) Variance of estimates and bounds

downward bias in estimating the full-information optimum $z_{\text{stoch}}$.[5] Fortunately, 2-SAA seemingly corrects for this bias, yielding valid bound in this example.

Nonetheless, both SAA and 2-SAA exhibit large variability (large error bars), especially for small samples. This is perhaps more easily seen in Panels 4b and 4c where we plot the variance of the recommended order quantity and the estimate or bound. Indeed, it can be seen in these plots that by sacrificing half the data when computing an ordering decision, 2-SAA indeed exacerbates SAA's instability, generally exhibiting about twice the variance.

By contrast, the data-driven guarantees of [15] do not converge to the full-information optimum, but rather to the bound of [42]. This is intuitive since as $N \to \infty$, the asymptotics in (24) suggest the DUS of [15] converges to the DUS of [42]. As discussed, by interpreting the DUS of [15] as a hypothesis test, we can improve these thresholds using bootstrapping, yielding the "Delage & Ye '10 (bootstrapped)" bound in the figure. While bootstrapping improves the bound significantly, the asymptotic performance remains the same. On the other hand, this bound is valid at level $1 - \alpha$ and is significantly less variable than both SAA and 2-SAA.

---

[5] This is most easily seen by noting:

$$\mathbb{E}\left[\hat{z}_{\text{SAA}}\right] = \mathbb{E}\left[\min_{x \in X} \frac{1}{N} \sum_{j=1}^{N} c(x; \xi^j)\right] \leq \min_{x \in X} \frac{1}{N} \sum_{j=1}^{N} \mathbb{E}\left[c(x; \xi^j)\right] = z_{\text{stoch}} \leq \mathbb{E}\left[c(x_{\text{SAA}}; \xi)\right]$$

Fig. 5: The price of data in the newsvendor problem: average of true PoD in solid black and of its distribution-agnostic approximation (45) in dashed gray. Note the log scales.

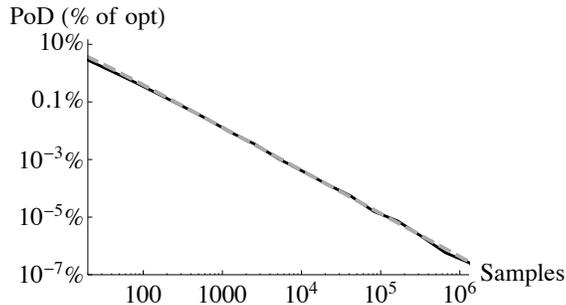

Finally, Robust SAA combines the strengths of all three approaches. It correctly converges to the full-information optimum, it is a valid $1 - \alpha$ bound on the true optimum (cf. Panel 4a), and it exhibits much smaller variance than both SAA and 2-SAA. In particular, in Panel 4b we see that the Robust SAA order quantity has significantly smaller variance than both the SAA and 2-SAA order quantities for $N \le 10^3$, and, for large $N$, the variance of Robust SAA and SAA are indistinguishable on a log-scale. (Notice though, that 2-SAA does have noticeably larger variance.) Similar effects are seen in Figure 4c.

### 7.2 Price of Data for a Newsvendor

Panel 4a demonstrates another typical property of SAA. Because of their downward bias, SAA estimates initially increase as we obtain more data, suggesting (incorrectly) that the price of data is negative. Neither 2-SAA nor the method of [15] share this defect, but even by roughly "eye-balling" the derivatives of these curves, it seems clear that they cannot provide a good estimate of the value of an additional data point.

By contrast, the PoD as estimated using Robust SAA is extremely accurate for this example. In Figure 5 we show the actual PoD and our estimate from Section 6 using Robust SAA. They are virtually indistinguishable on a log-scale.

### 7.3 Single-item newsvendor with unbounded distributions and other GoF tests

We next extend our previous test to consider unbounded distributions and our other univariate GoF tests. We consider unbounded variants of each of our previous three distributions, using the DUS $\mathcal{F}_{S_N, M_N}^{\alpha_1, \alpha_2}$ with $\phi(\xi) = \xi$ as in Example 2, and Student's T-test to compute $Q_{M_N(\alpha_2)}$ as outlined in Section 3.2. We set $\alpha_1 = 15\%$ and $\alpha_2 = 5\%$, yielding a total significance of 20%, comparable to previous example. We use IPOPT 3.11 [53] to solve the optimization problem for the AD test, and Gurobi otherwise. We compare to the same methods as the previous section.

Fig. 6: Probabilistic guarantees of Robust SAA (i.e., $\overline{z}$) for the singled-item newsvendor problem with significance 20%.

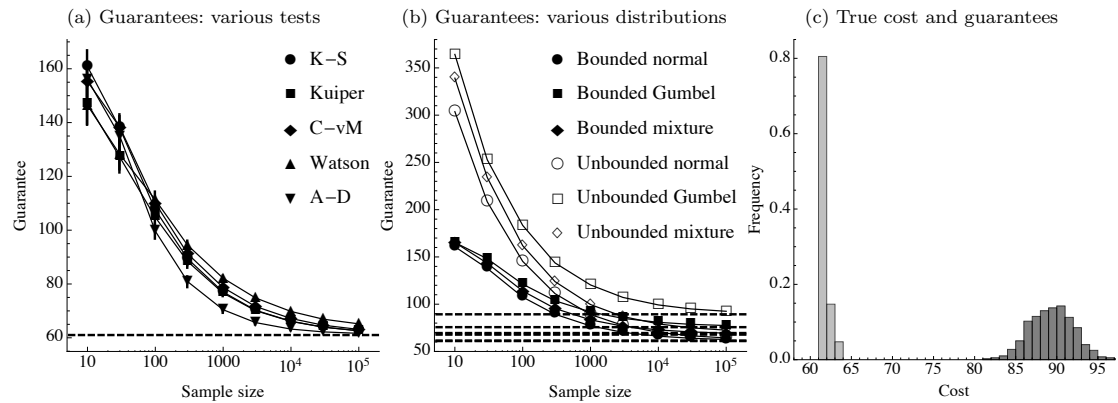

Fig. 7: The probabilistic guarantees of Robust SAA for the multi-item newsvendor problem. The vertical lines denote the span from the $20^{\text{th}}$ to the $80^{\text{th}}$ percentile. The dashed line denotes the full-information optimum.

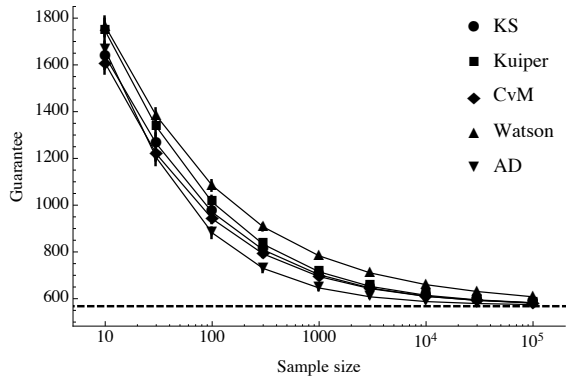

Figure 6 summarizes our findings. In Panel 6a, we consider the same bounded normal distribution from the previous section and see similar convergence and variability for the various tests. This observation holds for the other distributions as well. In Panel 6b, we see the performance of Robust SAA with the KS test for the different distributions. Notice, we observe asymptotic convergence (at a similar rate) for each choice. In both panels, dashed lines indicate the full-information optimum $z_{\text{stoch}}$. Panel 6c displays the distribution of true costs $z = \mathbb{E}_F[c(\overline{x}; \xi)]$ (light gray) and guarantees $\overline{z}$ (dark gray) for the Kolmogorov-Smirnov test with $N = 300$ samples from the bounded normal distribution.

We believe these numerical results confirm that the guarantees and asymptotic convergence of Robust SAA hold irrespective of the unknown distribution $F$ and which uniformly consistent test is used. Different tests also seem to yield mostly comparable results, with the AD test providing slightly better results when $N$ is at least 100. With small $N$, the Kuiper and Watson tests seem to perform the best. These observations should not, however, be taken as general conclusions about the relative performance of these tests for general problems. The conservatism of the guarantees depends both on the structure of the cost function as well as the true, unknown distribution and how we test it. For practical purposes, if the convergence rates are comparable as they are here, we recommend to choose the test that yields the simplest optimization problem, which in this case is the KS test.

### 7.4 Multi-item newsvendor

We now consider the multi-item newsvendor problem, which is a special case of a separable cost function as considered in Section 5.3.1. In the multi-item newsvendor we have $X = \{x \in \mathbb{R}_+^d : \sum_{i=1}^d x_i \leq \kappa\}$ for some capacity $\kappa$ and

$$c(x; \xi) = \sum_{i=1}^d c_i(x_i; \xi_i),$$

where each $c_i$ takes the form of a newsvendor cost function with its own parameters $b_i$, $h_i$.

We consider the case of three items, each having demand distributed as one of the three bounded distributions considered in the single-item case, with the parameters $\kappa = 250$, $b_1 = 9$, $b_2 = 6$, $b_3 = 3$, $h_1 = 3$, $h_2 = 2$, $h_3 = 1$. In our application of Robust SAA we employ the test based on marginals where, for different choices of univariate test, we use the same GoF test for each marginal, each at significance of $\alpha_i = 6.67\%$ (total significance 20%).

We present the results in Figure 7. As before, the center lines represent the mean over several draws of the data and error bars represent the $20^{\text{th}}$ and $80^{\text{th}}$ percentiles, and the number of draws for each $N$ was chosen so that the mean has a standard error of no more than 0.1 standard deviations. The main observation is that the prediction of the $c$-consistency theory holds: we observe convergence of guarantees even though testing marginals is not generally a uniformly consistent test but rather a $c$-consistent test for this particular problem.

### 7.5 Portfolio Allocation

We now consider the minimum-CVaR portfolio allocation problem as studied in Example 6. Portfolio allocation is well-known to be a challenging problem for SAA and, consequently, has been studied in the data-driven setting by a number of authors [15, 16, 26, 32]. We minimize the 10%-level CVaR of negative returns of a portfolio of $d = 10$ securities. The random returns are supported on the unbounded domain $\mathbb{R}^{10}$ and given by the factor model

$$\xi_i = \frac{i}{11}\tau + \frac{11-i}{11}\zeta_i, \quad i = 1, \ldots, 10, \tag{46}$$

where $\tau$ is a common market factor following a normal distribution with mean 2.5% and standard deviation 3% and $\zeta_i$'s are independent idiosyncratic contributions all following a negative Pareto distribution with upper support 3.7%, mean 2.5%, and standard deviation 3.8% (i.e. $\zeta_i \sim 0.05 - \text{Pareto}(0.013, 2.05)$). All securities have the same average return. Lower indexed securities are more volatile but are also more independent of the market factor. In our opinion, these features – a low-dimensional factor model and long tails – are common in financial data. We plot the PDFs of the returns of a few of the securities in Figure 8.

For samples drawn from this distribution, we consider data-driven solutions by the SAA, 2-SAA, the DRO of [15], and Robust SAA using the test for LCX. We use a total significance of $\alpha = 20\%$ for each method. We use the bootstrap to compute $Q_{C_N}(\alpha)$ (see Section 10.3). The constants $\gamma_{1,N}(\alpha)$, $\gamma_{2,N}(\alpha)$, $\gamma_{3,N}(\alpha)$ (see (22)) for the DRO of [15] are only developed therein for the case of bounded support, so in order to offer a fair comparison, we bootstrap these thresholds. We use GUROBI 5.5 [24] to solve all optimization problems.

We summarize our results in Figure 9. In Panel 9a we present the various estimates and bounds as a function of $N$. As before, the center lines represent the mean over several draws of the data and error bars represent the 20th and 80th percentiles, and the number of draws for each $N$ was chosen so that the mean has a standard error of no more than 0.1 standard deviations.

Notice that similar to our single-item newsvendor example, both SAA and 2-SAA converge to the full information optimum. Unlike our previous example, however, we observe that SAA has a very low variance. Unfortunately, it is significantly biased downwards. Indeed, it (incorrectly) appears to outperform the full-information optimum for $N \leq 10^3$.

By contrast, 2-SAA does *seem* to provide an upper bound, albeit a highly variable one. We emphasize "seem" because the full-information optimum falls well-within the error bars for 2-SAA in Panel 9a. In other words, the 2-SAA bound does not hold with probability close to the desired $1 - \alpha = 80\%$. This fact is more clearly seen in Panel 9b where we plot the frequency with which the 2-SAA guarantee actually upper bounds its out-of-sample performance. One can see that the guarantee only holds about 50-60% of the time, well below the desired 80%. This failure of 2-SAA can be attributed to the high-degree of non-normality in the cost estimates on the validation sample, which is in turn caused by the long-tails in the factor model (46). In other words, the Student T approximation in (5) is very inaccurate. Numerical simulation confirms that $N$ must be very large before this approximation is acceptably accurate ($N > 10^7$). In our experiments, we have also used other approximations in place of the Student T approximation in (5), e.g., the bootstrap, and observed similar behavior. The fundamental issue is that it is difficult to create a small, valid confidence interval for the mean for a non-parametric, long-tailed distribution without large amounts of data. This difficulty poses a significant obstacle for 2-SAA and similar methods. We emphasize that that this example is not pathological. Long tails are common in many financial applications.

Fig. 8: The PDFs of the distributions of security returns considered for the portfolio allocation problem.

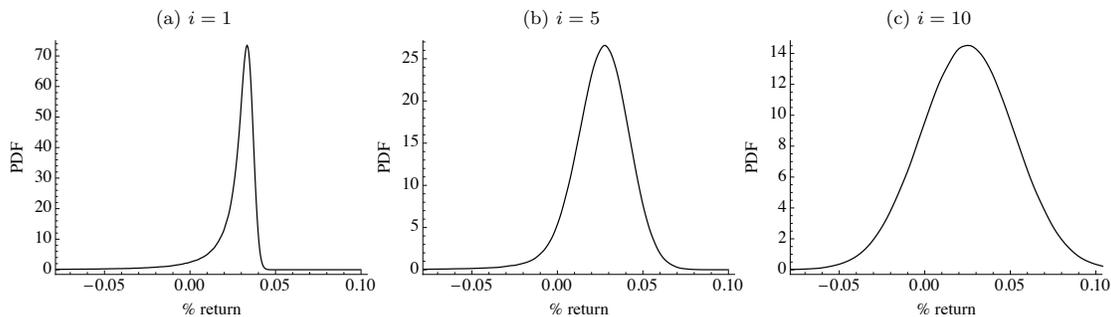

Fig. 9: Comparison of various data-driven approaches to the portfolio allocation problem (42). Where applicable, center lines are means over replicates and error bars are the $20^{\text{th}}$ and $80^{\text{th}}$ percentiles.

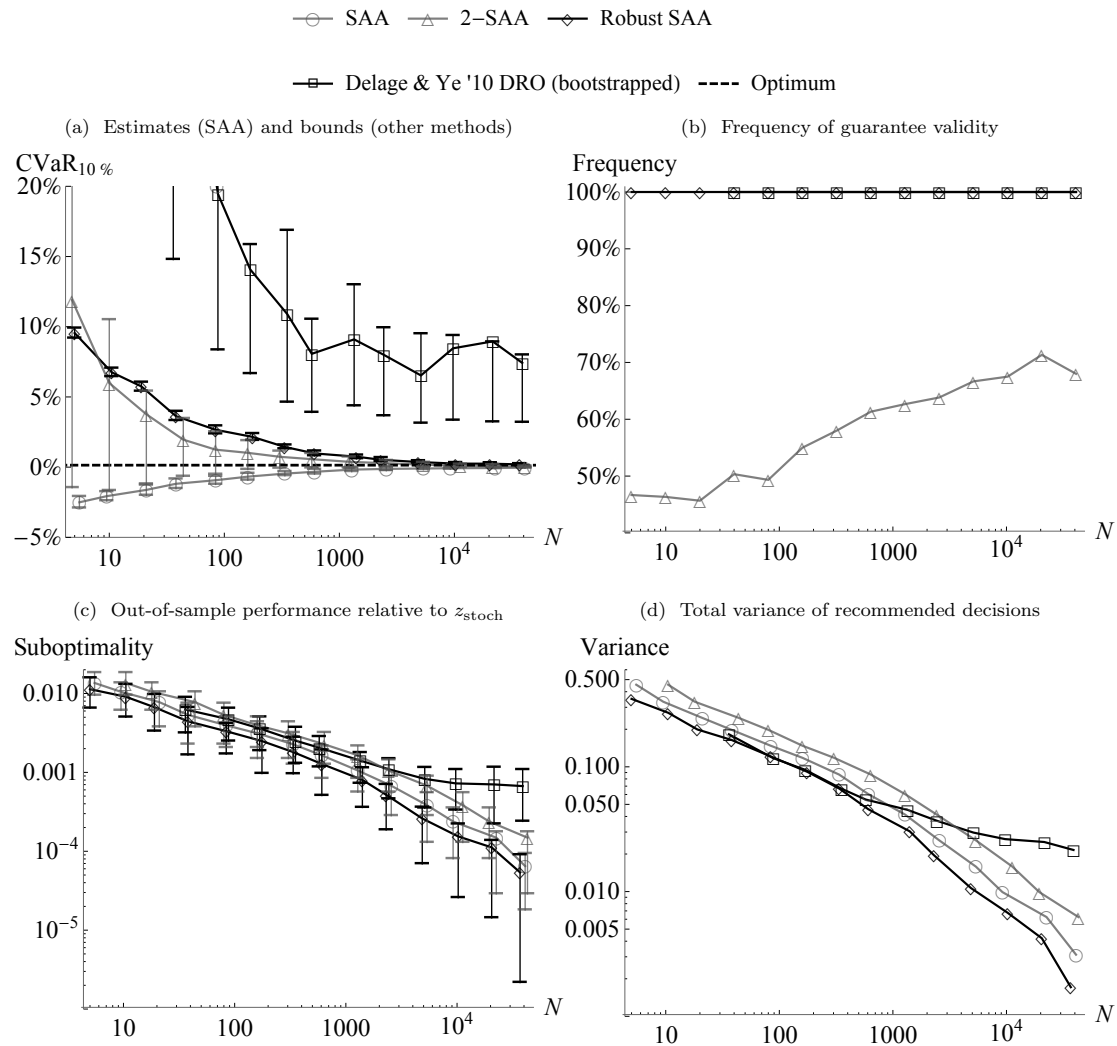

Finally, the method of [15] does provide a valid bound that holds with the desired significance, but does not converge to the full-information optimum (cf. panel 9a). Moreover, it is extremely variable, even for large $N$.

By contrast, Robust SAA again seems to draw on the strengths of the three approaches. It converges correctly to the full-information optimum, and it exhibits low-variability similar to SAA (cf. Panel 9a). Moreover, it is a valid bound with probability at least 80% on the out-of-sample performance. In fact, as seen in Panel 9b the bound is actually valid over all of our samples, i.e., with probability close to 1 for all $N$.

While Panel 9a describes the bounds and (in-sample) estimates, Panel 9c illustrates the out-of-sample performance of each method in terms of suboptimality to the full-information optimum. One can see that Robust SAA outperforms all of the previous methods on average for all values of $N$.

Finally, Panel 9d plots the total variance (sum of the variance of each component) of the recommended portfolio for each method. Roughly speaking, portfolios with large variance may incur large transaction costs from the frequent need to rebalance. One can see clearly that again, Robust SAA yields the lowest variance solutions of the compared approaches. Admittedly, in a real-life application of portfolio allocation we would likely introduce additional constraints over time to explicitly control the transaction costs, but we still feel that studying the variance of the unconstrained portfolio is insightful.

Fig. 10: The price of data in portfolio optimization: average of true PoD in solid black and of its resampling-based approximation (44) in dashed gray. Note the vertical log scale.

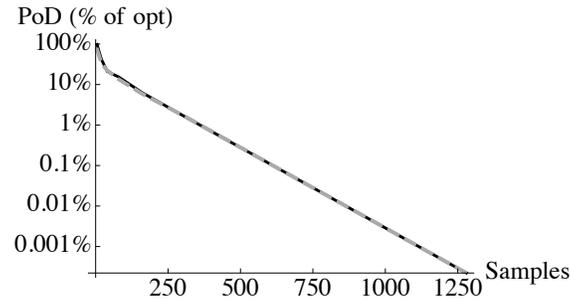

To conclude the experiment, we can again consider the price of data. In Figure 10 we compare the true price of data (43) for the LCX-based DRO bound and the resampling based approximation of it (44). On a log-scale, the two values are again indistinguishable.

## 8 Conclusion

In this paper, we proposed a novel, tractable approach to data-driven optimization called robust sample average approximation (Robust SAA). Robust SAA enjoys the tractability and finite-sample performance guarantees of many existing data-driven methods, but, unlike those methods, additionally exhibits asymptotic behavior similar to traditional sample average approximation (SAA). The key to the approach is a novel connection between SAA, DRO, and statistical hypothesis testing.

In particular, we were able to link properties of a data-driven *optimization* problem, i.e., its finite sample and asymptotic performance, to *statistical* properties of an associated goodness-of-fit hypothesis test, i.e., its significance and consistency. As a theoretical consequence, this connection allow us to describe the finite sample and asymptotic performance of both Robust SAA and other data-driven DRO formulations. As a practical consequence, our hypothesis testing perspective first, sheds light on which data-driven DRO formulations are likely to perform well in particular applications and second, enables us to use powerful, numerical tools like bootstrapping to improve their performance. Numerical experiments in inventory management and portfolio allocation confirm that our new method Robust SAA is tractable and can outperform existing data-driven methods in these applications.

## 9 Acknowledgements

The authors would like to thank the anonymous reviewers and associate editor for their extremely helpful suggestions and very thorough review of the paper.

# 10 Appendix

## 10.1 Proof of Theorem 1

*Proof* We will show that $\mathcal{C}(x;\mathcal{F})$ is continuous in $x$ and that when $c(x;\xi)$ satisfies the coerciveness conditions, $\mathcal{C}(x;\mathcal{F})$ is also coercive. The result will then follow from the usual Weierstrass extreme value theorem for deterministic optimization [5].

Let $S = \{x \in X : \mathcal{C}(x;\mathcal{F}) < \infty\}$. By assumption $x_0 \in S$, so $S \neq \varnothing$. Consequently, we can restrict to minimizing over $S$ instead of over $X$.

Fix any $x \in X$. Let $\epsilon > 0$ be given. By equicontinuity of the cost at $x$ there is a $\delta > 0$ such that any $y \in X$ with $||x - y|| \leq \delta$ has $|c(x;\xi) - c(y;\xi)| \leq \epsilon$ for all $\xi \in \Xi$. Fix any such $y$. Then

$$\mathcal{C}(y;\mathcal{F}) = \sup_{F_0 \in \mathcal{F}} \mathbb{E}_{F_0}[c(y;\xi)] \leq \sup_{F_0 \in \mathcal{F}} \mathbb{E}_{F_0}[c(x;\xi)] + \epsilon = \mathcal{C}(x;\mathcal{F}) + \epsilon, \tag{47}$$

$$\mathcal{C}(x;\mathcal{F}) = \sup_{F_0 \in \mathcal{F}} \mathbb{E}_{F_0}[c(x;\xi)] \leq \sup_{F_0 \in \mathcal{F}} \mathbb{E}_{F_0}[c(y;\xi)] + \epsilon = \mathcal{C}(y;\mathcal{F}) + \epsilon. \tag{48}$$

Note that (48) implies that $S$ is closed relative to $X$, which is itself closed. Hence, $S$ is closed. Furthermore, (47) and (48) imply that $\mathcal{C}(x;\mathcal{F})$ is continuous in $x$ on $S$.

If $S$ is compact, the usual Weierstrass extreme value theorem provides that the continuous $\mathcal{C}(x;\mathcal{F})$ attains its minimal (finite) value at an $x \in S \subseteq X$.

Suppose $S$ is not compact. Since $S \subseteq X$ is closed, this must mean $S$ is unbounded and then $X$ is unbounded and hence not compact. Then, by assumption, $c(x;\xi)$ satisfies the coerciveness assumption. Because $S$ is unbounded, there exists a sequence $x_i \in S$ such that $\lim_{i\to\infty} ||x_0 - x_i|| = \infty$. Then by the coerciveness assumption, $c_i(\xi) = c(x_i;\xi)$ diverges pointwise to infinity. Fix any $F_0 \in \mathcal{F}$. Let $c_i'(\xi) = \inf_{j \geq i} c_j(\xi)$, which is then pointwise monotone nondecreasing and pointwise divergent to infinity. Then, by Lebesgue's monotone convergence theorem, $\lim_{i\to\infty} \mathbb{E}_{F_0}[c_i'(\xi)] = \infty$. Since $c_i' \leq c_j$ pointwise for any $j \geq i$, we have $\mathbb{E}_{F_0}[c_i'(\xi)] \leq \inf_{j \geq i} \mathbb{E}_{F_0}[c_j(\xi)]$ and therefore

$$\infty = \lim_{i\to\infty} \mathbb{E}_{F_0}[c_i'(\xi)] \leq \lim_{i\to\infty} \inf_{j \geq i} \mathbb{E}_{F_0}[c_i(\xi)] = \liminf_{i\to\infty} \mathbb{E}_{F_0}[c_i(\xi)].$$

Thus $\mathcal{C}(x;\mathcal{F}) \geq \mathbb{E}_{F_0}[c(x;\xi)]$ is also coercive in $x$ over $S$. Then, the usual Weierstrass extreme value theorem provides that the continuous $\mathcal{C}(x;\mathcal{F})$ attains its minimal (finite) value at an $x \in S \subseteq X$.   □

## 10.2 Proof of Proposition 2

*Proof* Suppose that $c(x;\xi) \to \infty$ as $\xi \to \infty$. The case of unboundedness in the negative direction is similar. Choose $\rho > 0$ small so that $\xi^{(i)} - \xi^{(i-1)} > 2\rho$ for all $i$. For $\delta > 0$ and $\xi' \geq \xi^{(N)} + \rho$, let $F_{\delta,\xi'}$ be the measure with density function

$$f(\xi;\delta,\xi') = \begin{cases} 1/(2N\rho) & \xi^{(i)} - \rho \leq \xi \leq \xi^{(i)} + \rho \text{ for } 1 \leq i \leq N-1, \\ 1/(2N\rho) & \xi^{(N)} - \rho + \delta\rho \leq \xi \leq \xi^{(N)} + \rho - \delta\rho, \\ 1/(2N\rho) & \xi' \leq \xi \leq \xi' + 2\delta\rho, \\ 0 & \text{otherwise.} \end{cases}$$

In words, this density equally distributes mass on a $\rho$-neighborhood of every point $\xi^{(i)}$, except $\xi^{(N)}$. For $\xi^{(N)}$, we "steal" $\delta$ of the mass to place around $\xi'$. Notice that for any $\xi'$, $F_{0,\xi'}$ (i.e., take $\delta = 0$) minimizes $S_N(F_0)$ over distributions $F_0$. Since $\alpha > 0$, $Q_{S_N}(\alpha)$ is strictly greater than this minimum. Since $S_N(F_{\delta,\xi'})$ increases continuously with $\delta$ independently of $\xi'$, there must exist $\delta > 0$ small enough so that $F_{\delta,\xi'} \in \mathcal{F}_{S_N}^\alpha$ for any $\xi' > \xi^{(N)} + \rho$.

Let any $M > 0$ be given. By infinite limit of the cost function, there exists $\xi' > \xi^{(N)} + \rho$ sufficiently large such that $c(x;\xi) \geq MN/\delta$ for all $\xi \geq \xi'$. Then, we have $\mathcal{C}(x;\mathcal{F}_{S_N}^\alpha) \geq \mathbb{E}_{F_{\delta,\xi'}}[c(x;\xi)] \geq \mathbb{P}(\xi \geq \xi') MN/\delta = M$.

Since we have shown this for every $M > 0$, we have $\mathcal{C}(x;\mathcal{F}_{S_N}^\alpha) = \infty$.   □

10.3 Computing a threshold $Q_{C_N}(\alpha)$

We provide two ways to compute $Q_{C_N}(\alpha)$ for use with the LCX-based GoF test. One is an exact, closed form formula, but which may be loose for moderate $N$. Another uses the bootstrap to compute a tighter, but approximate threshold.

The theorem below employs a bound on $\mathbb{E}_F\left[\|\xi\|_2^2\right]$ to provide a valid threshold. This bound could either stem from known support bounds or from changing (15) to a two-sided hypothesis with two-sided confidence interval, using the lower bound as in (17) and the upper bound in (49) given below.

**Theorem 16** *Let $N \geq 2$. Suppose that with probability at least $1 - \alpha_2$, $\mathbb{E}_F\left[\|\xi\|_2^2\right] \leq \overline{Q_{R_N}}(\alpha_2)$. Let $\alpha_1 \in (0,1)$ be given and suppose $F_0 \preceq_{LCX} F$. Then, with probability at least $1 - \alpha_1 - \alpha_2$,*

$$\mathbb{E}_F\left[\|\xi\|_2^2\right] \leq \overline{Q_{R_N}}(\alpha_2) \quad and$$

$$C_N(F_0) \leq \left(1 + \overline{Q_{R_N}}(\alpha_2)\right)\left(1 + \frac{p}{2-p}\right)\frac{2^{\frac{1}{2}+\frac{1}{p}}}{N^{1-\frac{1}{p}}}\sqrt{d+1+(d+1)\log\left(\frac{N}{d+1}\right)+\log\left(\frac{4}{\alpha_1}\right)}, \quad (49)$$

*where*

$$p = \frac{1}{2}\left(\sqrt{\log(256)+8\log(N)+(\log(2N))^2}-\log(2N)\right) \in (1,2). \quad (50)$$

*Hence, defining $Q_{C_N}(\alpha_1)$ equal to the right-hand side of (49), we get a valid threshold for $C_N$ in testing $F_0 \preceq_{LCX} F$ at level $\alpha_1$.*

*Proof* Observe

$$
\begin{aligned}
C_N(F_0) &\leq \sup_{\|a\|_1+|b|\leq 1}\left(\mathbb{E}_{F_0}[\max\{a^T\xi-b,0\}]-\mathbb{E}_F[\max\{a^T\xi-b,0\}]\right)\\
&\quad + \sup_{\|a\|_1+|b|\leq 1}\left(\mathbb{E}_F[\max\{a^T\xi-b,0\}]-\frac{1}{N}\sum_{i=1}^N\max\{a^T\xi^i-b,0\}\right)\\
&\leq \sup_{\|a\|_1+|b|\leq 1}\left(\mathbb{E}_F[\max\{a^T\xi-b,0\}]-\frac{1}{N}\sum_{i=1}^N\max\{a^T\xi^i-b,0\}\right), \quad (51)
\end{aligned}
$$

where the first inequality follows by distributing the sup and the second inequality follows because $F_0 \preceq_{LCX} F$. Next, we provide a probabilistic bound on the last sup, which is the maximal difference over $\|a\|_1 + |b| \leq 1$ between the true expectation of the hinge function $\max\{a^T\xi-b,0\}$ and its empirical expectation.

The class of level sets of such functions, i.e., $\mathcal{S} = \left\{\{\xi \in \Xi : \max\{a^T\xi-b,0\} \leq t\} : \|a\|_1 + |b| \leq 1, \, t \in \mathbb{R}\right\}$, is contained in the class of the empty set and all halfspaces. Therefore, it has Vapnik-Chervonenkis dimension at most $d+1$ (cf. [52]). Therefore, Theorem 5.2 and equation (5.12) of [52] provide that for any $\epsilon > 0$ and $p \in (1,2]$,

$$
\begin{aligned}
\mathbb{P}\Bigg(\sup_{\|a\|_1+|b|\leq 1}\frac{1}{D_p(a,b)}&\left(\mathbb{E}_F[\max\{a^T\xi-b,0\}]-\frac{1}{N}\sum_{i=1}^N\max\{a^T\xi^i-b,0\}\right) > \epsilon\Bigg) \quad (52)\\
&< 4\exp\left(\left(\frac{(d+1)(\log\left(\frac{N}{d+1}\right)+1)}{N^{2-2/p}}-\frac{\epsilon^2}{2^{1+2/p}}\right)N^{2-2/p}\right),
\end{aligned}
$$

where $D_p(a,b) = \int_0^\infty\left(\mathbb{P}_F\left(\max\{a^T\xi-b,0\}>t\right)\right)^{1/p}dt$.

Notice that for any $\|a\|_1 + |b| \leq 1$, we have $0 \leq \max\{a^T\xi-b,0\} \leq \max\{1,\|\xi\|_\infty\} \leq \max\{1,\|\xi\|_2\}$ and hence $\mathbb{E}_F\left[\max\{a^T\xi-b,0\}^2\right] \leq \max\left\{1^2,\mathbb{E}_F\left[\|\xi\|_2^2\right]\right\} \leq 1 + \mathbb{E}_F\left[\|\xi\|_2^2\right]$. These observations combined with Markov's inequality yields that, for any $\|a\|_1 + |b| \leq 1$ and $p \in (1,2)$, we have

$$
\begin{aligned}
D_p(a,b) = \int_0^\infty\left(\mathbb{P}_F\left(\max\{a^T\xi-b,0\}>t\right)\right)^{1/p}dt &\leq 1 + \int_1^\infty\frac{\left(\mathbb{E}_F\left[\max\{a^T\xi-b,0\}^2\right]\right)^{1/p}}{t^{2/p}}dt\\
&\leq \left(1+\mathbb{E}_F\left[\|\xi\|_2^2\right]\right)^{1/p}\left(1+\frac{p}{2-p}\right) \leq \left(1+\mathbb{E}_F\left[\|\xi\|_2^2\right]\right)\left(1+\frac{p}{2-p}\right).
\end{aligned}
$$

This yields a bound on $D_p(a, b)$ that is independent of $(a, b)$. Using this bound takes $D_p(a, b)$ out of the sup in (52) and by bounding $\mathbb{E}_F\left[\|\xi\|_2^2\right] \leq \overline{Q_{R_N}}(\alpha_2)$ (which holds with probability $1 - \alpha_2$), we conclude from (52) that (49) holds with probability $1 - \alpha_1 - \alpha_2$ for any $p \in (1, 2)$. The $p$ given in (50) optimizes the bound for $N \geq 2$. □

Next we show how to bootstrap an approximate threshold $Q_{C_N}(\alpha)$. Recall that we seek a threshold $Q_{C_N}(\alpha)$ such that $\mathbb{P}\left(C_N(F_0) > Q_{C_N}(\alpha)\right) \leq \alpha$ whenever $F_0 \preceq_{\mathrm{LCX}} F$. Employing (51), we see that a sufficient threshold is the $(1 - \alpha)^{\mathrm{th}}$ quantile of

$$\sup_{\|a\|_1 + |b| \leq 1} \left( \mathbb{E}_F[\max\{a^T\xi - b, 0\}] - \frac{1}{N}\sum_{i=1}^{N} \max\{a^T\xi^i - b, 0\} \right),$$

where $\xi^i$ are drawn IID from $F$. The bootstrap [20] approximates this by replacing $F$ with the empirical distribution $\hat{F}_N$. In particular, given an iteration count $B$, for $t = 1, \ldots, B$ it sets

$$Q^t = \sup_{\|a\|_1 + |b| \leq 1} \left( \frac{1}{N}\sum_{i=1}^{N} \max\{a^T\xi^i - b, 0\} - \frac{1}{N}\sum_{i=1}^{N} \max\{a^T\tilde{\xi}^{t,i} - b, 0\} \right) \tag{53}$$

where $\tilde{\xi}^{t,i}$ are drawn IID from $\hat{F}_N$, i.e., IID random choices from $\{\xi^1, \ldots, \xi^N\}$. Then the bootstrap approximates $Q_{C_N}(\alpha)$ by the $(1-\alpha)^{\mathrm{th}}$ quantile of $\{Q^1, \ldots, Q^B\}$. However, it may be difficult to compute (53) as the problem is non-convex. Fortunately (53) can be solved with a standard MILP formulation or by discretizing the space and enumerating (the objective is Lipschitz).

In particular, our bootstrap algorithm for computing $Q_{C_N}(\alpha)$ is as follows:

---

Input: $\xi^1, \ldots, \xi^N$ drawn from $F$, significance $0 < \alpha < 1$, precision $\delta > 0$, iteration count $B$
Output: Threshold $Q_{C_N}(\alpha)$ such that $\mathbb{P}\left(C_N(F_0) > Q_{C_N}(\alpha)\right) \lesssim \alpha$ whenever $F_0 \preceq_{\mathrm{LCX}} F$.

---

For $t = 1, \ldots, B$:
1. Draw $\tilde{\xi}^{t,1}, \ldots, \tilde{\xi}^{t,N}$ IID from $\hat{F}_N$.
2. Solve $Q^t = \sup_{\|a\|_1 + |b| \leq 1} \left( \frac{1}{N}\sum_{i=1}^{N} \max\{a^T\xi^i - b, 0\} - \frac{1}{N}\sum_{i=1}^{N} \max\{a^T\tilde{\xi}^{t,i} - b, 0\} \right)$ to precision $\delta$.
Sort $Q^{(1)} \leq \cdots \leq Q^{(B)}$ and return $Q^{(\lceil(1-\alpha)B\rceil)} + \delta$.

---

### 10.4 Proof of Proposition 4

*Proof* We first prove that a uniformly consistent test is consistent. Let $G_0 \neq F$ be given. Denote by $d$ the Lévy-Prokhorov metric, which metrizes weak convergence (cf. [9]), and observe that $d(G_0, F) > 0$.

Define $R_N = \sup_{F_0 \in \mathcal{F}_N} d(F_0, F)$. We claim that if the test is uniformly consistent, then $\mathbb{P}(R_N \to 0) = 1$. Suppose that for some sample path, $R_N \not\to 0$. By the definition of the supremum, there must exist $\delta > 0$ and a sequence $F_N \in \mathcal{F}_N$ such that $d(F_N, F) \geq \delta$ i.o. Since $d$ metrizes weak convergence, this means that $F_N$ does not converge to $F$. However, $F_N \in \mathcal{F}_N$ for all $N$, i.e. it is never rejected. Therefore, by uniform consistency of the test, the event that $R_N \not\to 0$ must have probability 0. I.e., $R_N \to 0$ a.s.

Since a.s. convergence implies convergence in probability, we have that $\mathbb{P}(R_N < \epsilon) \to 1$ for every $\epsilon > 0$, and, in particular, for $\epsilon = d(G_0, F) > 0$. Hence we have,

$$\mathbb{P}(G_0 \in \mathcal{F}_N) \leq \mathbb{P}(d(G_0, F) \leq R_N) \to 0.$$

This proves the first part of the proposition.

For the second part, we describe a test which is consistent but not uniformly consistent. Consider testing a continuous distribution $F$ with the following univariate GoF test:

Given data $\xi^1, \ldots, \xi^N$ drawn from $F$ and a hypothetical continuous distribution $F_0$:

Let $\ell = \lfloor \log_2 N \rfloor$, $k = N - 2^\ell$.

If $\frac{k}{2^\ell} \leq F_0(\xi^1) \leq \frac{k+1}{2^\ell}$ then $F_0$ is not rejected.

Otherwise, reject $F_0$ if it is rejected by the KS test at level $\frac{\alpha}{1 - 2^{-\ell}}$ applied to the data $\xi^2, \ldots, \xi^N$.

Notice that under the null-hypothesis, the probability of rejection is

$$\mathbb{P}(F_0 \text{ rejected }) = \mathbb{P}\left(F_0(\xi^1) \notin \left[\frac{k}{2^\ell}, \frac{k+1}{2^\ell}\right]\right) \mathbb{P}(F_0 \text{ is rejected by the KS test }) = (1 - 2^{-\ell})\frac{\alpha}{1 - 2^{-\ell}} = \alpha,$$

where we've used that $\xi^1$ is independent of the rest of the sample, and $F_0(\xi^1)$ is uniformly distributed for $F_0$ continuous. Consequently, the test is a valid GoF test and it has significance $\alpha$.

We claim this test is also consistent. Specifically, consider any $F_0 \neq F$. By continuity of $F_0$ and consistency of the KS test,

$$\mathbb{P}\left(F_0 \text{ is rejected}\right) = \mathbb{P}\left(F_0(\xi^1) \notin \left[\frac{k}{2^\ell}, \frac{k+1}{2^\ell}\right]\right) \mathbb{P}\left(F_0 \text{ is rejected by the KS test}\right) \longrightarrow 1.$$

However, the test is not uniformly consistent. Fix any continuous $F_0 \neq F$ and let

$$F_N = \begin{cases} F_0 & \text{if } \frac{k}{2^\ell} \leq F_0(\xi^1) \leq \frac{k+1}{2^\ell}, \\ \hat{F}_N & \text{otherwise.} \end{cases}$$

Observe that $0 \leq F_0(\xi^1) \leq 1$ and $[0,1] = \bigcup_{k=0}^{2^\ell - 1}\left[\frac{k}{2^\ell}, \frac{k+1}{2^\ell}\right]$. That is, for every $\ell$, $F_N = F_0$ at least once for $N \in \{2^\ell, \dots, 2^{\ell+1} - 1\}$. Therefore $F_N = F_0$ i.o., so it does not converge weakly to $F$. However, as constructed, $F_N$ is never rejected by the above test. This is done for every sample path so the test cannot be uniformly consistent.                                                                                                                  □

### 10.5 Proofs of Theorems 2 and 3

We first establish two useful results.

**Proposition 8** *Suppose $\mathcal{F}_N$ is the confidence region of a uniformly consistent test and that Assumptions 1 and 3 hold. Then, almost surely, $\mathbb{E}_{F_N}[c(x;\xi)] \to E_F[c(x;\xi)]$ for any $x \in X$ and sequences $F_N \in \mathcal{F}_N$.*

*Proof* Restrict to the a.s. event that $(F_N \not\to F \implies F_N \notin \mathcal{F}_N$ i.o.). Fix $F_N \in \mathcal{F}_N$. Then the contrapositive gives $F_N \to F$. Fix $x$. If $\Xi$ is bounded (Assumption 3a) then the result follows from the portmanteau lemma (see for example Theorem 2.1 of [9]). Suppose otherwise (Assumption 3b). Then $\mathbb{E}_{F_N}[\phi(\xi)] \to \mathbb{E}_F[\phi(\xi)]$. By Theorem 3.6 of [9], $\phi(\xi)$ is uniformly integrable over $\{F_1, F_2, \dots\}$. Since $c(x;\xi) = O(\phi(\xi))$, it is also uniformly integrable over $\{F_1, F_2, \dots\}$. Then the result follows by Theorem 3.5 of [9].                                                  □

The following is a restatement of the equivalence between local uniform convergence of continuous functions and convergence of evaluations along a convergent path. We include the proof for completeness.

**Proposition 9** *Suppose Assumption 1 holds and $\mathcal{C}(x_N; \mathcal{F}_N) \to \mathbb{E}_F[c(x;\xi)]$ for any convergent sequence $x_N \to x$. Then (18) holds.*

*Proof* Let $E \subseteq X$ compact be given and suppose to the contrary that $\sup_{x \in E} |\mathcal{C}(x; \mathcal{F}_N) - \mathbb{E}_F[c(x;\xi)]| \not\to 0$. Then $\exists \epsilon > 0$ and $x_N \in E$ such that $|\mathcal{C}(x_N; \mathcal{F}_N) - \mathbb{E}_F[c(x_N;\xi)]| \geq \epsilon$ i.o. This, combined with compactness, means that there exists a subsequence $N_1 < N_2 < \cdots < N_k \to \infty$ such that $x_{N_k} \to x \in E$ and $|\mathcal{C}(x_{N_k}; \mathcal{F}_{N_k}) - \mathbb{E}_F[c(x_{N_k};\xi)]| \geq \epsilon$ $\forall k$. Then,

$$0 < \epsilon \leq |\mathcal{C}(x_{N_k}; \mathcal{F}_{N_k}) - \mathbb{E}_F[c(x_{N_k};\xi)]| \leq |\mathcal{C}(x_{N_k}; \mathcal{F}_{N_k}) - \mathbb{E}_F[c(x;\xi)]| + |\mathbb{E}_F[c(x;\xi)] - \mathbb{E}_F[c(x_{N_k};\xi)]|.$$

By assumption, $\exists k_1$ such that $|\mathcal{C}(x_{N_k}; \mathcal{F}_{N_k}) - \mathbb{E}_F[c(x;\xi)]| \leq \epsilon/4$ $\forall k \geq k_1$. By equicontinuity and $x_{N_k} \to x$, $\exists k_2$ such that $|c(x;\xi) - c(x_{N_k};\xi)| \leq \epsilon/4$ $\forall \xi, k \geq k_2$. Then,

$$|\mathbb{E}_F[c(x;\xi)] - \mathbb{E}_F[c(x_{N_k};\xi)]| \leq \mathbb{E}_F[|c(x;\xi) - c(x_{N_k};\xi)|] \leq \epsilon/4 \quad \forall \xi, k \geq k_2.$$

Combining and considering $k = \max\{k_1, k_2\}$, we get the contradiction $\epsilon \leq \epsilon/2$ for strictly positive $\epsilon$.                                                  □

We prove the "if" and "only if" sides of Theorem 2 separately.

*Proof (Proofs of Theorem 3 and the "only if" side of Theorem 2)* For either theorem restrict to the a.s. event that

$$\mathbb{E}_{F_N}[c(x;\xi)] \to \mathbb{E}_F[c(x;\xi)] \text{ for every } x \in X \text{ and sequences } F_N \in \mathcal{F}_N \tag{54}$$

(using Proposition 8 for Theorem 2 or by assumption of $c$-consistency for Theorem 3).

Let any convergent sequence $x_N \to x$ and $\epsilon > 0$ be given. By equicontinuity and $x_N \to x$, $\exists N_1$ such that $|c(x_N;\xi) - c(x;\xi)| \leq \epsilon/2 \; \forall \xi, \; N \geq N_1$. Then, $\forall N \geq N_1$,

$$|\mathcal{C}(x_N;\mathcal{F}_N) - \mathcal{C}(x;\mathcal{F}_N)| \leq \sup_{F_0 \in \mathcal{F}_N} |\mathbb{E}_{F_0}[c(x_N;\xi) - c(x;\xi)]| \leq \sup_{F_0 \in \mathcal{F}_N} \mathbb{E}_{F_0}[|c(x_N;\xi) - c(x;\xi)|] \leq \epsilon/2.$$

By definition of supremum, $\exists F_N \in \mathcal{F}_N$ such that $\mathcal{C}(x;\mathcal{F}_N) \leq \mathbb{E}_{F_N}[c(x;\xi)] + \epsilon/4$. By (54), $\mathbb{E}_{F_N}[c(x;\xi)] \to \mathbb{E}_F[c(x;\xi)]$. Hence, $\exists N_2$ such that $|\mathbb{E}_{F_N}[c(x;\xi)] - \mathbb{E}_F[c(x;\xi)]| \leq \epsilon/4 \; \forall N \geq N_2$. Combining these with the triangle inequality

$$|\mathcal{C}(x_N;\mathcal{F}_N) - \mathbb{E}_F[c(x;\xi)]| \leq |\mathcal{C}(x_N;\mathcal{F}_N) - \mathcal{C}(x;\mathcal{F}_N)| + |\mathcal{C}(x;\mathcal{F}_N) - \mathbb{E}_F[c(x;\xi)]|,$$

we get

$$|\mathcal{C}(x_N;\mathcal{F}_N) - \mathbb{E}_F[c(x;\xi)]| \leq \epsilon \quad \forall N \geq \max\{N_1, N_2\}.$$

Thus, by Proposition 9, we get that (18) holds.

Let $A_N = \arg\min_{x \in X} \mathcal{C}(x;\mathcal{F}_N)$. We now show that $\bigcup_N A_N$ is bounded. If $X$ is compact (Assumption 2a) then this is trivial. Suppose $X$ is not compact (Assumption 2b). Using the same arguments as in the proof of Theorem 1, we have in particular that $\lim_{||x|| \to \infty} \mathbb{E}_F[c(x;\xi)] = \infty$, $z_{\text{stoch}} = \min_{x \in X} \mathbb{E}_F[c(x;\xi)] < \infty$, that $A = \arg\min_{x \in X} \mathbb{E}_F[c(x;\xi)]$ is compact, and each $A_N$ is compact. Let $x^* \in A$. Fix $\epsilon > 0$. By definition of supremum $\exists F_N \in \mathcal{F}_N$ such that $\mathcal{C}(x^*;\mathcal{F}_N) \leq \mathbb{E}_{F_N}[c(x^*;\xi)] + \epsilon$. By (54), $\mathbb{E}_{F_N}[c(x^*;\xi)] \to \mathbb{E}_F[c(x^*;\xi)] = z_{\text{stoch}}$. As this is true for any $\epsilon$ and since $\min_{x \in X} \mathcal{C}(x;\mathcal{F}_N) \leq \mathcal{C}(x^*;\mathcal{F}_N)$, we have $\limsup_{N \to \infty} \min_{x \in X} \mathcal{C}(x;\mathcal{F}_N) \leq z_{\text{stoch}}$. Now, suppose for contradiction that $\bigcup_N A_N$ is unbounded, i.e. there is a subsequence $N_1 < N_2 < \cdots < N_k \to \infty$ and $x_{N_k} \in A_{N_k}$ such that $||x_{N_k}|| \to \infty$. Let $\delta' = \limsup_{k \to \infty} \inf_{\xi \notin D} c(x_{N_k};\xi) \geq \liminf_{N \to \infty} \inf_{\xi \notin D} c(x_N;\xi) > -\infty$ and $\delta = \min\{0, \delta'\}$. By $D$-uniform coerciveness, $\exists k_0$ such that $c(x_{N_k};\xi) \geq (z_{\text{stoch}} + 1 - \delta)/F(D) \; \forall \xi \in D, \; k \geq k_0$. In the case of Theorem 2, let $F_N$ be any $F_N \in \mathcal{F}_N$. In the case of Theorem 3, let $F_N$ be the empirical distribution $F_N = \hat{F}_N \in \mathcal{F}_N$. In either case, we get $F_N \to F$ weakly. In particular, $F_N(D) \to F(D)$. Then $\mathbb{E}_{F_N}[c(x_{N_k};\xi)] \geq F_N(D) \times (z_{\text{stoch}} + 1 - \delta)/F(D) + \min\{0, \inf_{\xi \notin D} c(x_{N_k};\xi)\} \; \forall k \geq k_0$. Thus $\limsup_{N \to \infty} \min_{x \in X} \mathcal{C}(x;\mathcal{F}_N) \geq \limsup_{k \to \infty} \min_{x \in X} \mathcal{C}(x;\mathcal{F}_{N_k}) \geq z_{\text{stoch}} + 1 - \delta + \delta = z_{\text{stoch}} + 1$, yielding the contradiction $z_{\text{stoch}} + 1 \leq z_{\text{stoch}}$.

Thus $\exists A_\infty$ compact such that $A \subseteq A_\infty$, $A_N \subseteq A_\infty$. Then, by (18),

$$\delta_N = \left| \min_{x \in X} \mathcal{C}(x;\mathcal{F}_N) - \min_{x \in X} \mathbb{E}_F[c(x;\xi)] \right| = \left| \min_{x \in A_\infty} \mathcal{C}(x;\mathcal{F}_N) - \min_{x \in A_\infty} \mathbb{E}_F[c(x;\xi)] \right|$$
$$\leq \sup_{x \in A_\infty} |\mathcal{C}(x;\mathcal{F}_N) - \mathbb{E}_F[c(x;\xi)]| \to 0,$$

yielding (19). Let $x_N \in A_N$. Since $A_\infty$ is compact, $x_N$ has at least one convergent subsequence. Let $x_{N_k} \to x$ be any convergent subsequence. Suppose for contradiction $x \notin A$, i.e., $\epsilon = \mathbb{E}_F[c(x;\xi)] - z_{\text{stoch}} > 0$. Since $x_{N_k} \to x$ and by equicontinuity, $\exists k_1$ such that $|c(x_{N_k};\xi) - c(x;\xi)| \leq \epsilon/4 \; \forall \xi, \; k \geq k_1$. Then, $|\mathbb{E}_F[c(x_{N_k};\xi)] - \mathbb{E}_F[c(x;\xi)]| \leq \mathbb{E}_F[|c(x_{N_k};\xi) - c(x;\xi)|] \leq \epsilon/4 \; \forall k \geq k_1$. Also $\exists k_2$ such that $\delta_{N_k} \leq \epsilon/4 \; \forall k \geq k_2$. Then, for $k \geq \max\{k_1, k_2\}$,

$$\min_{x \in X} \mathcal{C}(x;\mathcal{F}_{N_k}) = \mathcal{C}(x_{N_k};\mathcal{F}_{N_k}) \geq \mathbb{E}_F[c(x_{N_k};\xi)] - \delta_{N_k} \geq \mathbb{E}_F[c(x;\xi)] - \epsilon/2 \geq z_{\text{stoch}} + \epsilon/2.$$

Taking limits, we contradict (19). □

*Proof (Proof of the "if" side of Theorem 2)* Consider any $\Xi$ bounded ($R = \sup_{\xi \in \Xi} ||\xi|| < \infty$). Let $X = \mathbb{R}^d$, and

$$c_1(x;\xi) = ||x|| \left(2 + \cos\left(x^T\xi\right)\right), \quad c_2(x;\xi) = ||x|| \left(2 - \cos\left(x^T\xi\right)\right),$$
$$c_3(x;\xi) = ||x|| \left(2 + \sin\left(x^T\xi\right)\right), \quad c_4(x;\xi) = ||x|| \left(2 - \sin\left(x^T\xi\right)\right).$$

Since $|c_i(x;\xi)| \leq 3 ||x||$, expectations exist. The gradient of each $c_i$ at $x$ has magnitude bounded by $R ||x|| + 3$ uniformly over $\xi$, so Assumption 1 is satisfied. Also, $\lim_{||x|| \to \infty} c_i(x;\xi) \geq \lim_{||x|| \to \infty} ||x|| = \infty$ uniformly over all $\xi \in \Xi$ and $c_i(x;\xi) \geq 0$, so Assumption 2 is satisfied.

Restrict to the a.s. event that (18) applies simultaneously for $c_1, c_2, c_3, c_4$. Let any sequence $F_N \not\to F$ be given. Let $I$ denote the imaginary unit. Then, by the Lévy continuity theorem and Cramér-Wold device, there exists $x$ such that $\mathbb{E}_{F_N}\left[e^{Ix^T\xi}\right] \not\to \mathbb{E}_F\left[e^{Ix^T\xi}\right]$. On the other hand, by (18),

$$2\,||x|| + ||x|| \sup_{F_0 \in \mathcal{F}_N} \mathbb{E}_{F_0}\left[\cos\left(x^T\xi\right)\right] \longrightarrow 2\,||x|| + ||x||\,\mathbb{E}_F\left[\cos\left(x^T\xi\right)\right],$$

$$2\,||x|| - ||x|| \inf_{F_0 \in \mathcal{F}_N} \mathbb{E}_{F_0}\left[\cos\left(x^T\xi\right)\right] \longrightarrow 2\,||x|| - ||x||\,\mathbb{E}_F\left[\cos\left(x^T\xi\right)\right],$$

$$2\,||x|| + ||x|| \sup_{F_0 \in \mathcal{F}_N} \mathbb{E}_{F_0}\left[\sin\left(x^T\xi\right)\right] \longrightarrow 2\,||x|| + ||x||\,\mathbb{E}_F\left[\sin\left(x^T\xi\right)\right],$$

$$2\,||x|| - ||x|| \inf_{F_0 \in \mathcal{F}_N} \mathbb{E}_{F_0}\left[\sin\left(x^T\xi\right)\right] \longrightarrow 2\,||x|| - ||x||\,\mathbb{E}_F\left[\sin\left(x^T\xi\right)\right].$$

The first two limits imply that $\sup_{F_0 \in \mathcal{F}_N} \left|\mathbb{E}_{F_0}\left[\cos\left(x^T\xi\right)\right] - \mathbb{E}_F\left[\cos\left(x^T\xi\right)\right]\right| \to 0$ and the second two imply that $\sup_{F_0 \in \mathcal{F}_N} \left|\mathbb{E}_{F_0}\left[\sin\left(x^T\xi\right)\right] - \mathbb{E}_F\left[\sin\left(x^T\xi\right)\right]\right| \to 0$. Together, recalling that $e^{It} = \cos(t) + I\sin(t)$, this implies that

$$\sup_{F_0 \in \mathcal{F}_N} \left|\mathbb{E}_{F_0}\left[e^{Ix^T\xi}\right] - \mathbb{E}_F\left[e^{Ix^T\xi}\right]\right| \to 0.$$

However, since, $\mathbb{E}_{F_N}\left[e^{Ix^T\xi}\right] \not\to \mathbb{E}_F\left[e^{Ix^T\xi}\right]$, it must be that $F_N \notin \mathcal{F}_N$ i.o. □

## 10.6 Proof of Theorem 4

*Proof* In the case of finite support $\Xi = \{\hat{\xi}^1, \dots, \hat{\xi}^n\}$, total variation metrizes weak convergence:

$$d_{\mathrm{TV}}(q, q') = \frac{1}{2}\sum_{j=1}^n |q(j) - q'(j)|.$$

Restrict to the almost sure event $d_{\mathrm{TV}}(\hat{p}_N, p) \to 0$ (see Theorem 11.4.1 of [18]). We need only show that now $\sup_{p_0 \in \mathcal{F}_N} d_{\mathrm{TV}}(\hat{p}_N, p_0) \to 0$, yielding the contrapositive of the uniform consistency condition.

By an application of the Cauchy-Schwartz inequality (cf. [22]),

$$d_{\mathrm{TV}}(\hat{p}_N, p_0) = \frac{1}{2}\sum_{j=1}^n \frac{|\hat{p}_N(j) - p_0(j)|}{\sqrt{p_0(j)}} \times \sqrt{p_0(j)} \le \frac{1}{2}\left(\sum_{j=1}^n \frac{(\hat{p}_N(j) - p_0(j))^2}{p_0(j)}\right)^{1/2} = \frac{X_N(p_0)}{2}.$$

By Pinsker's inequality (cf. [30]),

$$d_{\mathrm{TV}}(\hat{p}_N, p_0) \le \frac{1}{\sqrt{2}}\left(\sum_{j=1}^n \sum_{j=1}^n \hat{p}_N(j)\log\left(\hat{p}_N(j)/p_0(j)\right)\right)^{1/2} = \frac{G_N(p_0)}{2}.$$

Since both the $\chi^2$ and G-tests use a rejection threshold equal to $\sqrt{Q/N}$ where $Q$ is the $(1-\alpha)^{\mathrm{th}}$ quantile of a $\chi^2$ distribution with $n-1$ degrees of freedom ($Q$ is independent of $N$), we have that $d_{\mathrm{TV}}(\hat{p}_N, p_0)$ is uniformly bounded over $p_0 \in \mathcal{F}_N$ by a quantity diminishing with $N$. □

## 10.7 Proof of Theorem 5

*Proof* In the case of univariate support, the Lévy metric metrizes weak convergence:

$$d_{\mathrm{Lévy}}(G, G') = \inf\{\epsilon > 0 : G(\xi - \epsilon) - \epsilon \le G'(\xi) \le G(\xi + \epsilon) + \epsilon\ \forall \xi \in \mathbb{R}\}.$$

Restrict to the almost sure event $d_{\mathrm{Lévy}}(\hat{F}_n, F) \to 0$ (see Theorem 11.4.1 of [18]). We need only show that now $\sup_{F_0 \in \mathcal{F}_N} d_{\mathrm{Lévy}}(\hat{F}_N, F_0) \to 0$, yielding the contrapositive of the uniform consistency condition.

Fix $F_0$ and let $0 \le \epsilon < d_{\mathrm{Lévy}}(\hat{F}_N, F_0)$. Then $\exists \xi_0$ such that either (1) $\hat{F}_N(\xi_0 - \epsilon) - \epsilon > F_0(\xi_0)$ or (2) $\hat{F}_N(\xi_0 + \epsilon) + \epsilon < F_0(\xi_0)$. Since $F_0$ is monotonically non-decreasing, (1) implies $D_N(F_0) \ge \hat{F}_N(\xi_0 - \epsilon) - F_0(\xi_0 - \epsilon) > \epsilon$ and (2) implies $D_N(F_0) \ge F_0(\xi_0 + \epsilon) - \hat{F}_N(\xi_0 + \epsilon) > \epsilon$. Hence $d_{\mathrm{Lévy}}(\hat{F}_N, F_0) \le D_N(F_0)$. Moreover, $D_N \le V_N$ by definition. Since $\sup_{F_0 \in \mathcal{F}^s_{S_N}} S_N(F_0) = Q_{S_N}(\alpha) = O(N^{-1/2})$ for either statistic, both the KS and Kuiper tests are uniformly consistent.

Consider $D'_N(F_0) = \max_{i=1,\,...,\,N}\left|F_0(\xi^{(i)}) - \frac{2i-1}{2N}\right| = \sigma\left(F_0(\xi^{(j)}) - \frac{2j-1}{2N}\right)$, where $j$ and $\sigma$ are the maximizing index and sign, respectively. Suppose $D'_N(F_0) \geq 1/\sqrt{N} + 1/N$. If $\sigma = +1$, this necessarily means that $1 - \frac{2j-1}{2N} \geq 1/\sqrt{N} + 1/N$ and therefore $N - j \geq \lceil\sqrt{N}\rceil + 1$. By monotonicity of $F_0$ we have for $0 \leq k \leq \lceil\sqrt{N}\rceil$ that $j + k \leq N$ and

$$F_0(\xi^{(j+k)}) - \frac{2(j+k)-1}{2N} \geq F_0(\xi^{(j)}) - \frac{2j-1}{2N} - \frac{k}{N} = D'_N(F_0) - \frac{k}{N} \geq 0.$$

If instead $\sigma = -1$, this necessarily means that $\frac{2j-1}{2N} \geq 1/\sqrt{N} + 1/N$ and therefore $j \geq \lceil\sqrt{N}\rceil + 1$. By monotonicity of $F_0$ we have for $0 \leq k \leq \lceil\sqrt{N}\rceil$ that $j - k \geq 1$ and

$$\frac{2(j-k)-1}{2N} - F_0(\xi^{(j-k)}) \geq \frac{2j-1}{2N} - F_0(\xi^{(j)}) - \frac{k}{N} = D'_N(F_0) - \frac{k}{N} \geq 0.$$

In either case we have that

$$W_N^2 = \frac{1}{12N^2} + \frac{1}{N}\sum_{i=1}^{N}\left(F_0(\xi^{(i)}) - \frac{2i-1}{2N}\right)^2 \geq \frac{1}{12N^2} + \frac{1}{N}\sum_{k=0}^{\lceil\sqrt{N}\rceil}\left(D'_N - \frac{k}{N}\right)^2 \geq \frac{D_N^2}{\sqrt{N}} - \frac{2}{N}$$

using $D'_N(F_0) \geq 1/\sqrt{N} + 1/N$ and $|D'_N(F_0) - D_N(F_0)| \leq 1/(2N)$ in the last inequality. Therefore,

$$D_N^2(F_0) \leq \max\left\{\frac{1}{\sqrt{N}} + \frac{3}{2N}, \sqrt{N}W_N^2(F_0) + \frac{2}{\sqrt{N}}\right\}.$$

Since $F_0(\xi)(1 - F_0(\xi)) \leq 1$ and by using the integral formulation of CvM and AD (see [49]) the same is true replacing $W_N^2$ by $A_N^2$. For either of the CvM or AD statistic $Q_{S_N}(\alpha) = O(N^{-1/2})$ and hence $\sup_{F_0 \in \mathcal{F}_{S_N}^\alpha} S_N^2(F_0) = O(N^{-1})$. Therefore, both the CvM and AD tests are uniformly consistent.

$$W_N^2 - U_N^2 = \left(\frac{1}{N}\sum_{i=1}^{N}F_0(\xi^{(i)}) - \frac{1}{2}\right)^2 \leq \max\left\{\left(\frac{1}{N}\sum_{i=1}^{N}\min\left\{1, \frac{2i-1}{2N} + D'_N(F_0)\right\} - \frac{1}{2}\right)^2,\right.$$
$$\left.\left(\frac{1}{N}\sum_{i=1}^{N}\max\left\{0, \frac{2i-1}{2N} - D'_N(F_0)\right\} - \frac{1}{2}\right)^2\right\}.$$

Letting $M = \lfloor\frac{1}{2} + N(1 - D'_N(F_0))\rfloor$ we have $\sum_{i=1}^{N}\min\left\{1, \frac{2i-1}{2N} + D'_N(F_0)\right\} = \frac{M^2}{2N} + MD'_N(F_0) + N - M$ so that in the case of $D'_N(F_0) \geq 1/\sqrt{N} + 1/N$, $\left(\frac{1}{N}\sum_{i=1}^{N}\min\left\{1, \frac{2i-1}{2N} + D'_N(F_0)\right\} - \frac{1}{2}\right)^2 = O(1/N)$. Thus, the Watson test is also uniformly consistent. □

## 10.8 Proof of Proposition 5

*Proof* Apply Theorem 2 to each $i$ and restrict to the almost sure event that (18) holds for all $i$. Fix $F_N$ such that $F_N \in \mathcal{F}_N$ eventually. Then, (18) yields $\mathbb{E}_{F_N}[c_i(x;\xi_i)] \to \mathbb{E}_F[c_i(x;\xi_i)]$ for every $x \in X$. Summing over $i$ yields the contrapositive of the $c$-consistency condition. □

## 10.9 Proof of Proposition 6

*Proof* Restrict to a sample path in the almost sure event $\mathbb{E}_{\hat{F}_N}[\xi_i] \to \mathbb{E}_{\hat{F}}[\xi_i]$, $\mathbb{E}_{\hat{F}_N}[\xi_i\xi_j] \to \mathbb{E}_{\hat{F}}[\xi_i\xi_j]$ for all $i, j$. Consider any $F_N$ such that $F_N \in \mathcal{F}_{\text{CEG},N}^\alpha$ eventually. Then clearly $\mathbb{E}_{F_N}[\xi_i] \to \hat{\mathbb{E}}_F[\xi_i]$, $\mathbb{E}_{F_N}[\xi_i\xi_j] \to \mathbb{E}_F[\xi_i\xi_j]$.

Consider any $F_N$ such that $F_N \in \mathcal{F}_{\text{DY},N}^\alpha$ eventually. Because covariances exist, we may restrict to $N$ large enough so that $\left\|\left|\hat{\Sigma}_N\right|\right\|_2 \leq M$ (operator norm) and $F_N \in \mathcal{F}_{\text{DY},N}^\alpha$. Then we get

$$\|\mathbb{E}_{F_N}[\xi] - \hat{\mu}_N\| \leq M\gamma_{1,N}(\alpha) \to 0$$

and

$$(\gamma_{3,N}(\alpha) - 1)\,\hat{\Sigma}_N \preceq \mathbb{E}_{F_N}[(\xi - \hat{\mu}_N)(\xi - \hat{\mu}_N)^T] - \hat{\Sigma}_N \preceq (\gamma_{2,N}(\alpha) - 1)\,\hat{\Sigma}_N,$$

which gives $\left\|\left|\mathbb{E}_{F_0}[(\xi - \hat{\mu}_N)(\xi - \hat{\mu}_N)^T] - \hat{\Sigma}_N\right|\right\|_2 \leq M\max\{\gamma_{2,N}(\alpha) - 1, 1 - \gamma_{3,N}(\alpha)\} \to 0$. Then again, we have $\mathbb{E}_{F_N}[\xi_i] \to \mathbb{E}_F[\xi_i]$, $\mathbb{E}_{F_N}[\xi_i\xi_j] \to \mathbb{E}_F[\xi_i\xi_j]$.

In either case we get $\mathbb{E}_{F_N}[c(x;\xi)] \to \mathbb{E}_{F_N}[c(x;\xi)]$ for any $x$ due to factorability as in (26). This yields the contrapositive of the $c$-consistency condition. □

### 10.10 Proof of Theorem 6

*Proof* If $F_0 \neq F$ then Theorem 1 of [43] yields that either $F_0 \npreceq_{\text{LCX}} F$ or there is some $j = 1, \ldots, d$ such that $\mathbb{E}_{F_0}[\xi_j^2] \neq \mathbb{E}_F[\xi_j^2]$. If $F_0 \npreceq_{\text{LCX}} F$ then probability of rejection approaches one since $C_N > 0$ but $Q_{C_N}(\alpha_1) \to 0$. Otherwise, $F_0 \preceq_{\text{LCX}} F$ yields $\mathbb{E}_{F_0}[\xi_i^2] \leq \mathbb{E}_F[\xi_i^2]$ for all $i$ via (12) using $a = e_i$ and $\phi(\zeta) = \zeta^2$. Then $\mathbb{E}_{F_0}[\xi_j^2] \neq \mathbb{E}_F[\xi_j^2]$ must mean that $\mathbb{E}_{F_0}\left[||\xi||_2^2\right] < \mathbb{E}_F\left[||\xi||_2^2\right]$ and probability of rejection goes to one.   □

### 10.11 Proof of Theorem 7

*Proof* Let $R = \sup_{\xi \in \Xi} ||\xi||_2 < \infty$. Restrict to the almost sure event that $\hat{F}_N \to F$. Consider $F_N$ such that $F_N \in \mathcal{F}_N$ eventually. Let $N$ be large enough so that it is so. Fix $||a||_2 = 1$. Let $a_1 = a$ and complete an orthonormal basis for $\mathbb{R}^d$: $a_1, a_2, \ldots, a_d$. On the one hand we have $Q_{R_N}(\alpha_2) \geq \mathbb{E}_{\hat{F}_N}\left[\sum_{i=1}^d (a_i^T \xi)^2\right] - \mathbb{E}_{F_N}\left[\sum_{i=1}^d (a_i^T \xi)^2\right]$. On the other hand, for each $i$,

$$
\begin{aligned}
\mathbb{E}_{\hat{F}_N}\left[(a_i^T \xi)^2\right] - \mathbb{E}_{F_N}\left[(a_i^T \xi)^2\right] &= 2\int_{b=-R}^0 \left(\mathbb{E}_{\hat{F}_N}[\max\{b - a_i^T \xi, 0\}] - \mathbb{E}_{F_N}[\max\{b - a_i^T \xi, 0\}]\right) db \\
&\quad + 2\int_{b=0}^R \left(\mathbb{E}_{\hat{F}_N}[\max\{a_i^T \xi - b, 0\}] - \mathbb{E}_{F_N}[\max\{a_i^T \xi - b, 0\}]\right) db \\
&\geq 4\int_{b=0}^R (||a||_1 + |b|)Q_{C_N}(\alpha_1) db \geq 4\left(\sqrt{d} + R^2/2\right)Q_{C_N}(\alpha_1) = p_N.
\end{aligned}
$$

Therefore, $q_N = Q_{R_N}(\alpha_2) + (d-1)p_N \geq \mathbb{E}_{\hat{F}_N}\left[(a^T \xi)^2\right] - \mathbb{E}_{F_N}\left[(a^T \xi)^2\right]$ and $Q_{R_N}(\alpha_2), Q_{C_N}(\alpha_1), p_N, q_N \to 0$. Let $G_N(t) = F_N(\{\xi : a^T \xi \leq t\}) \in [0,1]$ and $\hat{G}_N(t) = \hat{F}_N(\{\xi : a^T \xi \leq t\}) \in [0,1]$ be the CDFs of $a^T \xi$ under $F_N$ and $\hat{F}_N$, respectively. Then,

$$
\begin{aligned}
q_N \geq \mathbb{E}_{\hat{F}_N}\left[(a^T \xi)^2\right] - \mathbb{E}_{F_N}\left[(a^T \xi)^2\right] &= 2\int_{b=-R}^0 \left(\mathbb{E}_{\hat{F}_N}[\max\{b - a^T \xi, 0\}] - \mathbb{E}_{F_N}[\max\{b - a^T \xi, 0\}]\right) db \\
&\quad + 2\int_{b=0}^R \left(\mathbb{E}_{\hat{F}_N}[\max\{a^T \xi - b, 0\}] - \mathbb{E}_{F_N}[\max\{a^T \xi - b, 0\}]\right) db \\
&= 2\int_{b=-R}^0 \int_{t=-R}^b \left(\hat{G}_N(t) - G_N(t)\right) dt\, db \\
&\quad + 2\int_{b=0}^R \int_{t=b}^R \left(G_N(t) - \hat{G}_N(t)\right) dt\, db \geq p_N, \\
\int_{t=-R}^b \left(\hat{G}_N(t) - G_N(t)\right) dt &\geq -(\sqrt{d} + R)Q_{C_N}(\alpha) \quad \forall b \in [-R, 0], \\
\int_{t=b}^R \left(G_N(t) - \hat{G}_N(t)\right) dt &\geq -(\sqrt{d} + R)Q_{C_N}(\alpha) \quad \forall b \in [0, R],
\end{aligned}
$$

Because $\hat{F}_N \to F$, we get $\hat{G}_N(t) \to F(\{\xi : a^T \xi \leq t\})$ and therefore $G_N(t) \to F(\{\xi : a^T \xi \leq t\})$ at every continuity point $t$. Because this is true for every $a$, the Cramer-Wold device yields $F_N \to F$. This is the contrapositive of the uniform consistency condition.   □

### 10.12 Proof of Theorem 8

*Proof* The proof amounts to dualizing a finite-support phi-divergence as done in Theorem 1 of [3] and is included here only for the sake of self-containment.

Let $\phi_X(t) = (t-1)^2/t$ and $\phi_G(t) = -\log(t) + t - 1$ (corresponding to "$\chi^2$-distance" and "Burg entropy" in Table 2 of [3], respectively). Then we can write

$$
\mathcal{F}_{X_N}^\alpha = \left\{ p_0 \geq 0 : \sum_{j=1}^n p_0(j) = 1, \sum_{j=1}^n \tilde{p}(j)\phi_X(p_0(j)/\tilde{p}_N(j)) \leq (Q_{X_N}(\alpha))^2 \right\},
$$

$$\mathcal{F}^{\alpha}_{G_N} = \left\{ p_0 \geq 0 : \sum_{j=1}^{n} p_0(j) = 1, \ \sum_{j=1}^{n} \hat{p}(j)\phi_G(p_0(j)/\hat{p}_N(j)) \leq \frac{1}{2}(Q_{G_N}(\alpha))^2 \right\}.$$

Therefore, letting $c_j = c(x; \hat{\xi}^j)$ and fixing $\phi$ and $Q$ appropriately, the inner problem (4) is equal to

$$\max_{p_0 \in \mathbb{R}^n_+} \quad c^T p_0$$

$$\text{s.t.} \quad \sum_{j=1}^{n} p_0(j) = 1$$

$$\sum_{j=1}^{n} \hat{p}_N(j)\phi(p_0(j)/\hat{p}_N(j)) \leq Q^2.$$

By Fenchel duality, the above is equal to

$$\min_{r \in \mathbb{R}, \ s \in \mathbb{R}_+} \max_{p_0 \in \mathbb{R}^n_+} \quad c^T p_0 + r(1 - e^T p_0) + s(Q^2 - \sum_{j=1}^{n} \hat{p}_N(j)\phi(p_0(j)/\hat{p}_N(j)))$$

$$= \min_{r \in \mathbb{R}, \ s \in \mathbb{R}_+} \quad r + Q^2 s + \sum_{j=1}^{n} \hat{p}_N(j) s \max_{\rho \in \mathbb{R}_+} \left( \frac{c_j - r}{s}\rho - \phi(\rho) \right)$$

$$= \min_{r \in \mathbb{R}, \ s \in \mathbb{R}_+} \quad r + Q^2 s + \sum_{j=1}^{n} \hat{p}_N(j) s \phi^* \left( \frac{c_j - r}{s} \right),$$

where $\phi^*(\tau) = \sup_{\rho \geq 0}(\tau\rho - \phi(\rho))$ is the convex conjugate. Therefore, the inner problem (4) is equal to

$$\min_{r,s,t} \quad r + Q^2 s - \sum_{j=1}^{n} \hat{p}_N(j) t_j$$

$$\text{s.t.} \quad r \in \mathbb{R}, \ s \in \mathbb{R}_+, \ t \in \mathbb{R}^n, \ c \in \mathbb{R}^n$$

$$t_j \leq -s\phi^* \left( \frac{c_j - r}{s} \right) \qquad\qquad \forall j = 1, \ldots, N$$

$$c_j \geq c(x; \hat{\xi}^j) \qquad\qquad\qquad \forall j = 1, \ldots, N.$$

It is easy to verify that

$$\phi^*_X(\tau) = \begin{cases} 2 - 2\sqrt{1 - \tau} & \tau \leq 1 \\ \infty & \text{otherwise} \end{cases}, \quad \text{and} \quad \phi^*_G(\tau) = \begin{cases} -\log(1 - \tau) & \tau \leq 1 \\ \infty & \text{otherwise} \end{cases}$$

(see e.g. Table 4 of [3]). In the case of $X_N$, since $s \geq 0$, we get

$$t_j \leq -s\phi^*_X \left( \frac{c_j - r}{s} \right) \iff \begin{aligned} & c_j - r \leq s, \\ & 2s + t_j \leq 2\sqrt{s(s - c_j + r)} \end{aligned}$$

$$\iff \exists y_j : \begin{aligned} & c_j - r \leq s, \ y_j \geq 0, \\ & 2s + t_j \leq y_j, \\ & y_j^2 \leq (2s)(2s - 2c_j + 2r) \end{aligned}$$

$$\iff \exists y_j : \begin{aligned} & c_j - r \leq s, \ y_j \geq 0, \\ & 2s + t_j \leq y_j, \\ & y_j^2 + (r - c_j)^2 \leq (2s - c_j + r)^2. \end{aligned}$$

In the case of $G_N$, since $s \geq 0$, we get

$$t_j \leq -s\phi^*_G \left( \frac{c_j - r}{s} \right) \iff \begin{aligned} & c_j - r \leq s, \\ & se^{t_j/s} \leq s + r - c_j \end{aligned}$$

$$\iff (t_j, \ s, \ s + r - c_j) \in C_{\mathrm{XC}}.$$

$\square$

### 10.13 Proof of Theorem 9

*Proof* Problem (3) is equal to the optimization problems of Theorem 8 augmented with the variable $x \in X$ and weak optimization is polynomially reducible to weak separation (see [23]). Tractable weak separation for all constraints except $x \in X$ and (28) is given by the tractable weak optimization over these standard conic-affine constraints. A weak separation oracle is assumed given for $x \in X$. We polynomially reduce separation over $c_j \geq \max_k c_{jk}(x)$ for fixed $c'_j, x'$ to the oracles. We first call the evaluation oracle for each $k$ to check violation and if there is a violation and $k^* \in \arg\max_k c_{jk}(x')$ then we call the subgradient oracle to get $s \in \partial c_{jk^*}(x')$ with $||s||_\infty \leq 1$ and produce the separating hyperplane $0 \geq c_{jk^*}(x') - c_j + s^T(x - x')$. □

### 10.14 Proof of Theorem 10

*Proof* Substituting the given formulas for $K_{S_N}$, $A_{S_N}$, $b_{S_N,\alpha}$ for each $S_N \in \{D_N, V_N, W_N, U, N, A_N\}$ in $A_{S_N}\zeta - b_{S_N,\alpha} \in K_{S_N}$ we obtain exactly $S_N(\zeta_1, \ldots, \zeta_N) \leq Q_{S_N}(\alpha)$ for $S_N$ as defined in (8). We omit the detailed arithmetic. □

### 10.15 Proof of Theorem 13

*Proof* Under these assumptions (3) is equal to the optimization problems of Theorem 11 or 12 augmented with the variable $x$ and weak optimization is polynomially reducible to weak separation (see [23]). Tractable weak separation for all constraints except $x \in X$ and (30) is given by the tractable weak optimization over these standard conic-affine constraints. A weak separation oracle is assumed given for $x \in X$. By continuity and given structure of $c(x; \xi)$, we may rewrite (30) as

$$c_i \geq \max_{\xi \in [\xi^{(i-1)}, \xi^{(i)}]} c_k(x; \xi) \quad \forall k = 1, \ldots, K. \tag{55}$$

We polynomially reduce weak $\delta$-separation over the $k^{\text{th}}$ constraint at fixed $c'_i, x'$ to the oracles. We call the $\delta$-optimization oracle to find $\xi' \in [\xi^{(i-1)}, \xi^{(i)}]$ such that $c_k(x'; \xi') \geq \max_{\xi \in [\xi^{(i-1)}, \xi^{(i)}]} c_k(x; \xi) - \delta$. If $c'_i \geq c_k(x'; \xi')$ then we say the constraint and is within $\delta$ of $(c'_i, x')$. If $c'_i < c_k(x'; \xi')$ then we call the subgradient oracle to get $s \in \partial_x c_k(x', \xi')$ with $||s||_\infty \leq 1$ and produce the hyperplane $c_i \geq c(x'; \xi') + s^T(x - x')$ that is violated by $(c'_i, x')$ and for any $(c_i, x)$ satisfying (55) (in particular if it is in the $\delta$-interior) we have $c_i \geq \max_{\xi \in [\xi^{(i-1)}, \xi^{(i)}]} c_k(x; \xi) \geq c_k(x; \xi) \geq c_k(x'; \xi') + s^T(x - x')$ since $s$ is a subgradient. The case for constraints (31) is similar. □

### 10.16 Proof of Proposition 7

*Proof* According to Theorem 11, the observations in Example 2, and by renaming variables, the DRO (3) is given by

$$(P) \quad \min y + \sum_{i=1}^N \left(Q_{D_N}(\alpha) + \frac{i-1}{N}\right) s_i + \sum_{i=1}^N \left(Q_{D_N}(\alpha) - \frac{i}{N}\right) t_i$$

$$\text{s.t. } x \in \mathbb{R}_+, \, y \in \mathbb{R}, \, s \in \mathbb{R}_+^N, \, t \in \mathbb{R}_+^N$$

$$(r-c)x + y + \sum_{i=j}^N (s_i - t_i) \geq (r-c)\xi^{(j)} \qquad\qquad \forall j = 1, \ldots, N+1$$

$$-(c-b)x + y + \sum_{i=j}^N (s_i - t_i) \geq -(c-b)\xi^{(j-1)} \qquad\qquad \forall j = 1, \ldots, N+1.$$

Applying linear optimization duality we get that its dual is

$$(D) \quad \max \ (r-c) \sum_{i=1}^{N+1} \xi^{(i)} p_i - (c-b) \sum_{i=1}^{N+1} \xi^{(i-1)} q_i$$

$$\text{s.t. } p \in \mathbb{R}_+^{N+1}, \ q \in \mathbb{R}_+^{N+1}$$

$$(r-c) \sum_{i=1}^{N+1} p_i - (c-b) \sum_{i=1}^{N+1} q_i \leq 0$$

$$\sum_{i=1}^{N+1} p_i + \sum_{i=1}^{N+1} q_i = 1$$

$$\sum_{i=1}^{j} p_i + \sum_{i=1}^{j} q_i \leq Q_{D_N}(\alpha) + \frac{j-1}{N} \qquad \forall j = 1, \ldots, N$$

$$-\sum_{i=1}^{j} p_i - \sum_{i=1}^{j} q_i \leq Q_{D_N}(\alpha) - \frac{j}{N} \qquad \forall j = 1, \ldots, N.$$

It can be directly verified that the following primal and dual solutions are respectively feasible

$$x = (1-\theta)\xi^{(i_{\text{lo}})} + \theta \xi^{(i_{\text{hi}})},$$

$$y = (r-c)\xi^{(N+1)} - (r-c)x,$$

$$s_i = \begin{cases} (c-b)\left(\xi^{(i)} - \xi^{(i-1)}\right) & i \leq i_{\text{lo}} \\ 0 & \text{otherwise} \end{cases} \qquad \forall i = 1, \ldots, N,$$

$$t_i = \begin{cases} (r-c)\left(\xi^{(i+1)} - \xi^{(i)}\right) & i \geq i_{\text{hi}} \\ 0 & \text{otherwise} \end{cases} \qquad \forall i = 1, \ldots, N$$

$$p_i = \begin{cases} 0 & i \leq i_{\text{hi}} - 1 \\ i/N - \theta - Q_{D_N}(\alpha) & i = i_{\text{hi}} \\ 1/N & N \geq i \geq i_{\text{hi}} + 1 \\ Q_{D_N}(\alpha) & i = N + 1 \end{cases}, \qquad \forall i = 1, \ldots, N$$

$$q_i = \begin{cases} Q_{D_N}(\alpha) & i = 1 \\ 1/N & 2 \leq i \leq i_{\text{lo}} \\ \theta - Q_{D_N}(\alpha) - (i-2)/N & i = i_{\text{lo}} + 1 \\ 0 & i \geq i_{\text{lo}} + 2 \end{cases} \qquad \forall i = 1, \ldots, N$$

and that both have objective cost in their respective programs of

$$z = -(c-b)Q_{D_N}(\alpha)\xi^{(0)} - \frac{c-b}{N} \sum_{i=1}^{i_{\text{lo}}-1} \xi^{(i)} - (c-b)\left(\theta - Q_{D_N}(\alpha) - \frac{i_{\text{lo}} - 1}{N}\right)\xi^{(i_{\text{lo}})}$$

$$+ (r-c)Q_{D_N}(\alpha)\xi^{(N+1)} + \frac{r-c}{N} \sum_{i=i_{\text{hi}}+1}^{N} \xi^{(i)} + (r-c)\left(\frac{i_{\text{hi}}}{N} - Q_{D_N}(\alpha) - \theta\right)\xi^{(i_{\text{hi}})}.$$

This proves optimality of $x$. Adding $0 = (c-b)\theta x - (r-c)(1-\theta)x$ to the above yields the form of the optimal objective given in the statement of the result. □

## 10.17 Proof of Theorem 14

*Proof* Fix $x$. Let $S = \{(a, b) \in \mathbb{R}^{d+1} : ||a||_1 + |b| \leq 1\}$. Using the notation of [45], letting $C$ be the cone of nonnegative measures on $\Xi$ and $C'$ the cone of nonnegative measures on $S$, we write the inner problem as

$$\sup_{F} \ \langle F, c(x; \xi) \rangle_{\Xi}$$

$$\text{s.t. } F \in C, \ \langle F, 1 \rangle_{\Xi} = 1$$

$$\frac{1}{N} \sum_{i=1}^{N} \max\{a^T \xi^i - b, 0\} + Q_{C_N}(\alpha_1) - \langle F, \max\{a^T \xi - b, 0\} \rangle_{\Xi} \geq 0 \quad \forall (a, b) \in S$$

$$\langle F, ||\xi||_2^2 \rangle_{\Xi} \geq Q_{R_N}^{\alpha_2}$$

Invoking Proposition 2.8 of [45] (with the generalized Slater point equal to the empirical distribution), we have that the strongly dual minimization problem is

$$\min_{G, \tau, \theta} \quad \theta + \left\langle G, Q_{C_N}(\alpha_1) + \frac{1}{N} \sum_{i=1}^{N} \max\{a^T \xi^i - b, 0\} \right\rangle_S - Q_{R_N}^{\alpha_2} \tau \tag{56}$$

$$\text{s.t.} \quad G \in C', \tau \in \mathbb{R}_+, \theta \in \mathbb{R}$$

$$\inf_{\xi \in \Xi} \left( \left\langle G, \max\{a^T \xi - b, 0\} \right\rangle_S - c(x; \xi) - \tau \left\| \xi \right\|_2^2 \right) \geq -\theta. \tag{57}$$

If (40) is true, we can infer from constraint (57) that

$$\tau \leq \inf_{\xi \in \Xi} \frac{\theta + \left\langle G, \max\{a^T \xi - b, 0\} \right\rangle_S - c(x; \xi)}{\|\xi\|_2^2} \leq \liminf_{i \to \infty} \frac{\theta + G\{S\} \max\{\|\xi_i'\|_\infty, 1\} - c(x; \xi_i')}{\|\xi_i'\|_2^2} = 0.$$

This shows that $\tau = 0$ is the only feasible choice.

In general, $\tau = 0$ is a feasible choice and fixing it so provides an upper bound on (56) and hence on the original inner problem by weak duality. Moreover, plugging $\xi_0$ in to (57) we conclude that

$$\tau \leq \frac{\theta + G\{S\} \max\{\|\xi_0\|_\infty, 1\} - c(x; \xi_0)}{\|\xi_0\|_2^2} \leq \frac{1}{R^2} \theta + \frac{R+1}{R^2} G\{S\} - \frac{1}{R^2} c(x; \xi_0).$$

Hence, setting $\tau = 0$ in (57) and replacing $\tau$ by the above bound in the objective provides a lower bound on (56) and hence on the original inner problem.

In order to study both cases, and both the upper and lower bounds in the latter case, we consider for the rest of the proof the following general problem given $\eta$ and $\nu$,

$$\min_{G, \theta} \quad \eta \theta + \left\langle G, \nu + \frac{1}{N} \sum_{i=1}^{N} \max\{a^T \xi^i - b, 0\} \right\rangle_S \tag{58}$$

$$\text{s.t.} \quad G \in C', \theta \in \mathbb{R}$$

$$\inf_{\xi \in \Xi} \left( \left\langle G, \max\{a^T \xi - b, 0\} \right\rangle_S - c(x; \xi) \right) \geq -\theta. \tag{59}$$

We first rewrite (59) using the representation (37):

$$\inf_{\xi \in \Xi} \left( \left\langle G, \max\{a^T \xi - b, 0\} \right\rangle_S - c_k(x; \xi) \right) \geq -\theta \quad \forall k = 1, \dots, K.$$

Next, invoking [45], and employing the concave conjugate, we rewrite the $k$th of these constraints as follows:

$$-\theta \leq \inf_{\xi \in \Xi, \, g(\cdot) \geq 0} \sup_{H_k \in C'} \left( \left\langle G, g \right\rangle_S - c_k(x; \xi) + \left\langle H_k, a^T \xi - b - g \right\rangle_S \right)$$

$$= \sup_{H_k \in C', \, G - H_k \in C'} \left( \left\langle H_k, -b \right\rangle_S + c_{k*}(x; \left\langle H_k, a \right\rangle_S) \right).$$

Introducing the variables $H_k$ and this equivalent constraint into the problem (58) and invoking [45] we find that (58) is equal to the dual problem:

$$\max_{p, \, q \, \psi} \quad \sum_{k=1}^{K} p_k c_{k**}(x, \, q_k/p_k)$$

$$\text{s.t.} \quad p_k \in \mathbb{R}_+^K, q \in \mathbb{R}^{K \times d}, \inf_{(a,b) \in S} \psi_k(a, b) \geq 0 \, \forall k = 1, \dots, K$$

$$\sum_{k=1}^{K} p_k = \eta$$

$$\inf_{(a,b) \in S} \left( \nu + \frac{1}{N} \sum_{i=1}^{N} \max\{a^T \xi^i - b, 0\} - \sum_{k=1}^{K} \psi_k(a, b) \right) \geq 0$$

$$\inf_{(a,b) \in S} \left( \psi_k(a, b) - a^T q_k + p_k b \right) \geq 0 \quad \forall k = 1, \dots, K.$$

Since $c_k(x;\xi)$ is closed concave in $\xi$, $c_{k**}(x;\xi) = c_k(x;\xi)$. Moreover, recognizing that $\psi_k = \max\{a^T q_k - p_k b, 0\}$ is optimal, we can rewrite the above problem as:

$$\max_{p,\,q} \quad \sum_{k=1}^{K} p_k c_k(x,\, q_k/p_k) \tag{60}$$

$$\text{s.t.} \quad p_k \in \mathbb{R}_+^K,\, q \in RK \times d$$

$$\sum_{k=1}^{K} p_k = \eta$$

$$\inf_{(a,b)\in S}\left(\nu + \frac{1}{N}\sum_{i=1}^{N}\max\{a^T\xi^i - b, 0\} - \sum_{k=1}^{K}\max\{a^T q_k - p_k b, 0\}\right) \geq 0. \tag{61}$$

Next, we rewrite (61) by splitting it over the different branches on the sum of $K$ maxima, noting that the all-zeros branch is trivial:

$$\nu \geq \sup_{(a,b)\in S}\left(\sum_{k=1}^{k}\gamma_k\left(a^T q_k - p_k b\right) - \frac{1}{N}\sum_{i=1}^{N}\max\{a^T\xi^i - b, 0\}\right) \quad \forall \gamma \in \mathcal{G}.$$

Next, we invoke linear optimization duality to rewrite the $\gamma^{\text{th}}$ constraint as follows:

$$\nu \geq \sup_{a,b,\lambda,t,t'} \quad \left(\sum_{k=1}^{k}\gamma_k q_k\right)^T a - \left(\sum_{k=1}^{k}\gamma_k p_k\right)b - \sum_{i=1}^{N}\lambda_i$$

$$\text{s.t.} \quad \lambda \in \mathbb{R}_+^N,\, t \in \mathbb{R}_+^d,\, t' \in \mathbb{R}_+,\, a \in \mathbb{R}^d,\, b \in \mathbb{R}$$

$$\frac{1}{N}\left(a^T\xi - b\right) - \lambda_i \leq 0 \quad \forall i = 1, \ldots, N$$

$$t' + \sum_{j=1}^{d} t_i \leq 1$$

$$b - t' \leq 0 \quad -b - t' \leq 0$$

$$a - t \leq 0 \quad -a - t \leq 0.$$

$$= \inf_{\mu_\gamma,\rho_\gamma,u_\gamma,v_\gamma,u'_\gamma,v'_\gamma} \quad \rho_\gamma$$

$$\text{s.t.} \quad \mu_\gamma \in \mathbb{R}_+^N,\, \rho_\gamma \in \mathbb{R}_+,\, u_\gamma \in \mathbb{R}_+^d,\, v_\gamma \in \mathbb{R}_+^d,\, u'_\gamma \in \mathbb{R}_+,\, v'_\gamma \in \mathbb{R}_+$$

$$r_\gamma \leq (1, \ldots, 1)$$

$$s_\gamma \geq u'_\gamma + v'_\gamma$$

$$s_\gamma \geq u_{\gamma,j} + v_{\gamma,j} \quad \forall j = 1, \ldots, d$$

$$\sum_{k=1}^{K}\gamma_k q_k - \frac{1}{N}\sum_{i=1}^{N}\mu_{\gamma,i}\xi^i = u_\gamma - v_\gamma$$

$$-\sum_{k=1}^{K}\gamma_k p_k + \frac{1}{N}\sum_{i=1}^{N}\mu_{\gamma,i} = u'_\gamma - v'_\gamma.$$

Introducing the variables $\mu_\gamma, \rho_\gamma, u_\gamma, v_\gamma, u'_\gamma, v'_\gamma$ and this equivalent constraint into (60) and invoking linear optimization duality yields $\mathcal{C}'(x;\nu,\eta)$ and the result. $\qquad\square$

### 10.18 Proof of Theorem 15

*Proof* Under these assumptions, weak optimization is polynomially reducible to weak separation (see [23]) and separation is given for all constraints but (39). We polynomially reduce weak $\delta$-separation over the $k^{\text{th}}$ constraint at fixed $x', h'_k, g'_k$. First we call the concave conjugate oracle to find $\xi'$ such that $h'_k{}^T\xi' - c_k(x';\xi') \leq c_{k*}(x';h'_k) + \delta$. If $g'_k \leq h'_k{}^T\xi - c_k(x';\xi')$ then $x', h'_k, g'_k - \delta$ satisfies the constraint and is within $\delta$ of the given point. Otherwise, we call the subgradient oracle to get $s \in \partial_x c_k(x',\xi')$ with $\|s\|_\infty \leq 1$ and produce the hyperplane $g_k \leq h_k^T\xi' - c_k(x';\xi') - s^T(x - x')$, which is violated by the given point and is valid whenever the original constraint is valid (in particular in its $\delta$-interior). $\qquad\square$